\documentclass[opre,nonblindrev]{informs3} 

\OneAndAHalfSpacedXI 

\usepackage{endnotes}
\let\footnote=\endnote

\usepackage{natbib}
 \bibpunct[, ]{(}{)}{,}{a}{}{,}%
 \def\bibfont{\small}%
\usepackage{booktabs}

\usepackage{float}

\usepackage{bbm}

\newcolumntype{C}[1]{>{\centering\let\newline\\\arraybackslash\hspace{0pt}}m{#1}}

\TheoremsNumberedThrough     
\ECRepeatTheorems

\EquationsNumberedThrough    

\usepackage{setspace}

\usepackage{enumitem}
\usepackage{booktabs}
\usepackage{graphicx}
                 
\usepackage{bm}        

\usepackage[ruled,vlined]{algorithm2e}
\include{pythonlisting}
\usepackage{pdflscape}
\usepackage{mathrsfs}
\usepackage{multirow}

\newcommand{\field}[1]{\ensuremath{\mathbb{#1}}}
\newcommand{\sets}[1]{\ensuremath{\mathcal{#1}}}

\newcommand{\naturals}{\ensuremath{\field{N}}} 
\newcommand{\reals}{\ensuremath{\field{R}}} 
\newcommand{\1}{\ensuremath{{\rm \mathbf e}}} 

\newcommand{{\bu}}{\ensuremath{{\bm u}}}

\newcommand{\I}[1]{\ensuremath{\mathbb{I}\left(#1\right)}} 

\newcommand{\subjectto}{\text{\rm subject to}} 

\newcommand{\cl}{\rm{cl}}

\DeclareMathOperator{\st}{s.t.}

\newcommand{\minimize}{\ensuremath{\mathop{\mathrm{minimize}}\limits}}
\newcommand{\maximize}{\ensuremath{\mathop{\mathrm{maximize}}\limits}}

\definecolor{d}{rgb}{0.0, 0.5, 0.2}

\newcommand{\newpv}[1]{{\color{black}{#1}}}

\usepackage{mdframed}

\newcommand{\Rrisk}{$\mathcal{RP_{\rm{u}}}$}
\newcommand{\Rregret}{$\mathcal{RP_{\rm{r}}}$}

\newcommand{\offlinerisk}[1]{\ensuremath{\mathcal{QP}_{\rm{u}}^{#1}}}
\newcommand{\onlinerisk}[1]{\ensuremath{\mathcal{QP}_{\rm{on,u}}^{#1}}}
\newcommand{\widetildeofflinerisk}[1]{\ensuremath{\widetilde{\mathcal{QP}}_{\rm{u}}^{#1}}}

\newcommand{\offlineregret}[1]{\ensuremath{\mathcal{QP}_{\rm{r}}^{#1}}}

\newcommand{\widetildeofflineregret}[1]{\ensuremath{\widetilde{\mathcal{P}}_{\rm{r}}^{#1}}}

\renewcommand{\arraystretch}{1.5}

\begin{document}


\RUNAUTHOR{Vayanos et al.}

\RUNTITLE{Robust Active Preference Elicitation} 

\TITLE{Robust Active Preference Elicitation\footnote{Previous title: \emph{Active preference elicitation via adjustable robust optimization}}} 
%


%
\ARTICLEAUTHORS{%
\AUTHOR{Phebe Vayanos$^{\dagger}$, Yingxiao Ye, Duncan C. McElfresh, John P. Dickerson, Eric Rice}
\AFF{$^{\dagger}$University of Southern California, Center for Artificial Intelligence in Society, \EMAIL{phebe.vayanos@usc.edu}}
} 

\ABSTRACT{%
\newpv{We study the problem of strategically eliciting the preferences of a decision-maker through a moderate number of pairwise comparison queries with the goal of making them a high quality recommendation for a specific decision-making problem. We are particularly motivated by applications in high stakes domains, such as when choosing a policy for allocating scarce resources to satisfy basic human needs (e.g., kidneys for transplantation or housing for those experiencing homelessness) where a consequential recommendation needs to be made from the (partially) elicited preferences. We model uncertainty in the preferences as being set based and} investigate two settings: \emph{a)} an \emph{offline elicitation} setting, where all queries are made at once, and \emph{b)} an \emph{online elicitation} setting, where queries are selected sequentially over time in an adaptive fashion. We \newpv{propose robust optimization} formulations of these problems which integrate the preference elicitation and recommendation phases \newpv{with aim to either maximize worst-case utility or minimize worst-case regret,} and study their complexity. For the offline case, where active preference elicitation takes the form of a \emph{two and half stage} robust optimization problem with decision-dependent information discovery, \newpv{we provide an equivalent reformulation} in the form of a mixed-binary linear program which we solve via column-and-constraint generation. For the online setting, where active preference learning takes the form of a \emph{multi-stage} robust optimization problem with decision-dependent information discovery, we propose a conservative solution approach. \newpv{Numerical studies on synthetic data demonstrate that our methods outperform state-of-the art approaches from the literature in terms of worst-case rank, regret, and utility.}

\newpv{We showcase how our methodology can be used to assist a homeless services agency in choosing a policy for allocating scarce housing resources of different types to people experiencing homelessness. We use historical data to generate candidate counterfactual housing allocation policies that differ in their fairness, efficiency, and interpretability characteristics and elicit the agency's values over these attributes.}


\newpage

}%


\maketitle


%




\section{Introduction}
\label{sec:introduction}


\newpv{Automated recommendation and decision-support systems are routinely used in many aspects of our lives, e.g., to select adds to display~\citep{Li2010}, to suggest products to purchase~\citep{Gomez16:Netflix}, or a route to take~\citep{Dai2015}. Increasingly, they are used to make recommendations about consequential decisions in high stakes settings, e.g., to suggest a treatment plan to follow~\citep{Blasiak2020} or to prioritize transplant patients for deceased donor kidneys~\citep{bertsimas2013fairness}. Yet, many of these decisions are multifaceted, requiring hard trade-offs to be made. For example, when choosing a treatment, one needs to balance its effectiveness (e.g., reduce pain) with potential risks (e.g., side effects), see~\cite{Johnson2013}. 

Thus, to be able to make a recommendation likely to be adopted, the decision-support system needs to \emph{elicit} the preferences of its user. }
\newpv{Preference elicitation is achieved through either revealed preference data~\citep{Samuelson1948} or stated preference surveys~\citep{Kroes1988}. The most common surveys are \emph{conjoint studies}~\citep{Green1978,Green1990_ConjointAnalysis} which model choices as vectors of attributes (as in multi-attribute utility theory~\citep{Dyer1992}) and characterize preferences in terms of a vector whose elements indicate how much the user values each attribute. The survey consists of several queries that ask the user to either select their preferred item among two alternatives (choice-based conjoint) or to state the strength of their preference between two items (metric-paired comparisons). The former is often preferred as it yields fewer response errors and surveys typically involve a moderate number of queries ($\sim$10) to avoid tiring the user.} 

\newpv{In this paper, we study the problem of eliciting the preferences of a user using a moderate number of pairwise comparison (choice-based conjoint) queries with the goal of making them a high quality recommendation for a \emph{specific} decision-making problem. We call this problem, where elicitation and recommendation phases are inextricably linked, \emph{active preference elicitation}. We are particularly motivated by applications in high stakes low resource domains where revealed preference data is not available, where the number of attributes that characterize a recommendation is often large, and where a consequential recommendation needs to be made from the (partially) elicited preferences. We consider both an \emph{offline} elicitation setting, where all queries are selected in advance, and an \emph{online} elicitation setting, where the choice of query can be adapted to previous questions and answers. The former is relevant when working in low resource settings without access to a computer or when seeking to compare answers across users. The latter should be used whenever feasible to yield better recommendations. Given the critical nature of the problems we focus on, we aim to make \emph{robust recommendations} that either maximize worst-case utility or minimize worst-case regret accounting for information incompleteness, inconsistencies in the responses, and potential model misspecifications. We aim to identify queries that are as informative as possible for the recommendation problem even when answers to queries are ``adversarial'' in the sense that they provide as little information as possible to guide the choice of recommendation. This point of view ensures our system will make recommendations with guaranteed performance. Technically, this problem can be cast as a robust sequential decision-making problem where the answers to the queries are uncertain and will only be revealed if the corresponding query is asked; such problems where information is endogenous are often said to present \emph{decision-dependent information discovery} (DDID). We argue that active preference elicitation in this context is highly relevant to practical problems and requires the development of a new methodological framework.}


\newpv{A concrete motivation for our research is the plight of people experiencing homelessness and the ethical dilemma of deciding how to prioritize them for scarce housing resources. For example, in Los Angeles there are over 66,000 persons experiencing homelessness on any given night and only 22,769 housing units became available over the past year. This shortage makes it hard to decide on an allocation policy as moral trade-offs need to be made between efficiency (e.g., having the least possible returns to homelessness) and equity (e.g., giving people resources commensurate to their needs). These challenges are further complicated by preexisting disparities. For example, in LA, Black people are four times more represented among those experiencing homelessness than in the general population~\citep{BPEH_LAHSA_2018}. Some agencies may prefer allocation rules that do not exacerbate inequalities, while others may prefer giving people equal chances. From our discussions with policy-makers at the Los Angeles Homeless Services Authority (LAHSA), the entity in charge of allocating housing resources in LA, there is consensus that the current matching mechanism is inadequate as it does not align with particular policy objectives, resulting for example in Black people returning to homelessness at almost twice the rate of White individuals~\citep{BPEH_LAHSA_2018}. While selecting a policy among two candidates is feasible, articulating a full preference profile is very challenging. At the same time, the ethical implications of housing allocation require it to closely reflect the interests of its users. Systems that do not match the needs of those they seek to help may deter them from seeking access. Systems that do not serve the interests of their users will result in unintended consequences and may get defunded, exacerbating resource scarcity and causing additional strain. The role of policy-makers in this context is to find ways to finely balance the often conflicting interests of those they serve in a way that aligns with their own values.} 

\newpv{In recent years, several artificial intelligence driven tools have been proposed to design policies for matching homeless persons to housing~\citep{Azizi2018_CPAIOR,Sanmay_AAAI2019_housing}. However, these methods assume that the preferences of policy-makers are perfectly known. For example, \cite{Azizi2018_CPAIOR} proposed a method for designing interpretable policies that minimize returns to homelessness while being, in some predefined sense, fair. To the best of our knowledge, no tool exists for eliciting the value judgements of policy-makers in housing allocation agencies to design policies that align with them. In fact, in private communications, policy-makers at LAHSA have described this as a recurring challenge. This is not surprising considering the many difficulties involved in identifying queries whose responses will be most useful in informing the ultimate choice of policy. First, selecting queries optimally requires solving a sequential decision-making problem with decision-dependent information discovery and combinatorial number of uncertain response scenarios. Second, it requires explicitly accounting for possible inconsistencies in the responses while also allowing for potential model misspecifications.}



\newpv{The problem we study is not unique to the housing allocation system but is rather encountered broadly in the design of policy-making systems, in particular those that decide how to allocate scarce resources needed to satisfy basic needs. Another such system is the US kidney allocation system (KAS), which provides deceased donor kidneys for transplant to those suffering from end stage renal disease. As of 2019, there were approximately 95,000 patients on the national kidney waiting list and only 23,401 transplants took place. At the same time, there are severe inequities in kidney disease, with Black people being five times more likely to be affected than Whites~\citep{USRDS2019}. Thus, when deciding how to prioritize patients, the Organ Procurement and Transplantation Network (OPTN) faces ethical challenges and potential bad implications from inadequate allocation that are not unlike those faced in the housing allocation system.} 

\newpv{Our research objective is to select the best possible queries to ask a user faced with a specific decision-making problem, taking into account any prior information about their preferences. The choice of queries should result in robust recommendations with the best possible guaranteed performance no matter what the user preferences are, and accounting for information incompleteness, possible inconsistencies, and potential model misspecifications, making them suitable to deploy in high stakes settings--we refer to these as \emph{robust queries}. For example, we aim to elicit the moral value judgements of the leadership of a housing allocation agency to be able to recommend a policy that best aligns with them. The idea of examining and comparing sets of policies by looking at their outcomes is common in policy-making settings, e.g., it is used routinely in the KAS.} 

\newpv{The preference elicitation problem has been widely studied in the artificial intelligence (AI), operations research (OR), and marketing literatures, see Section~\ref{sec:related_work}. Several researchers have investigated elicitation strategies and preference models that can be used to make robust recommendations. Yet, this large body of work is not well suited to make recommendations in high stakes low resource settings. Indeed, existing methods focus on obtaining unbiased estimates of preferences or on reducing uncertainty around estimates, rather than on selecting queries that will help explicitly improve the quality of the robust recommendation. This implies they are likely to result in suboptimal recommendations that do not align well with the values of end users. Moreover, none of the preference elicitation methods that can be used to make robust recommendations apply to the offline setting. We consequently utilize techniques from modern adjustable robust optimization that will be able to provide us with optimality guarantees on the quality of the queries selected while also allowing us to effectively cope with model misspecifications, inconsistencies, and information incompleteness.}

\newpv{The problem of selecting comparison queries that will explicitly maximize the quality of the downstream recommendation is a robust optimization problem with DDID \citep{DDI_VKB,vayanos_ROInfoDiscovery}. In particular, active offline (resp.\ online) elicitation can be cast as a \emph{two and a half stage} (resp.\ \emph{multi-stage}) problem with DDID. Our framework does not postulate distributions for uncertain parameters and instead models stochasticity by means of optimization variables that lie in uncertainty sets. This approach builds upon the line of work on polyhedral methods in conjoint analysis, see Section~\ref{sec:related_work}, but provides a different means of modelling inconsistencies that better aligns with popular choice models from the literature. From a modeling standpoint, the key feature that sets us apart from existing work is that our choice of queries is precisely guided by the structure of the recommendation problem. From a solution perspective, this complicates significantly the problem of identifying good queries.}

\newpv{Indeed, when explicitly framed as a multi-stage robust optimization problem with DDID, the preference elicitation problem is particularly challenging: it involves binary adaptive decision variables (which also control the time of information discovery), discrete uncertainty (the responses to each queries) which may also be misspecified, an exponential number of response scenarios, a decision-dependent uncertainty set that involves ``non-equalities'' (to guarantee that responses are aligned with some utility model), and a decision-dependent uncertainty set that involves strict inequalities (to model strict preferences). This type of problem has not been previously studied in the literature on robust optimization, nor in the very nacent work on robust optimization with DDID, see Section~\ref{sec:related_work}, and we show that it is generally NP-hard. By combining techniques from polyhedral theory with methods from robust optimization, we demonstrate that the offline active preference elicitation problems that maximize utility (with arbitrary recommendation set) or minimize regret (with enumerable recommendation set) can be reformulated equivalently as mixed-binary linear optimization problems (MBLP) involving a number of decision variables and constraints that are exponential in the number of queries and propose novel column-and-constraint generation approaches to solve them. For the online setting, we propose conservative solution approaches that solve a series of offline problems in a rolling horizon fashion.}

\newpv{
We demonstrate the performance of our approaches in terms of quality of the recommendation and solution times by conducting extensive numerical studies using simulated data for realistic problem sizes. We also illustrate that changing the recommendation problem from max-min utility to min-max regret impacts the optimal choice of queries to be made, thereby showing the importance of coupling elicitation and recommendation phases.

We showcase how our methods can be used to assist policy-makers at agencies such as LAHSA in choosing a policy for allocating housing resources to people experiencing homelessness. We use historical data to generate candidate counterfactual allocation policies that differ in their fairness, efficiency, and interpretability characteristics and elicit stakeholder values over these attributes.

Our work contributes to the following literature streams. First, to preference elicitation by providing a framework for generating provably optimal robust queries for specific recommendation problems accounting for information incompleteness, model misspecifications, and response inconsistencies, making it suitable to deploy in high stakes settings. Second, to robust optimization with DDID by studying a model with discrete uncertain parameters, exponential number of contingencies, and uncertainty sets involving ``non-equalities'' and strict-inequalities. Finally, to the OR literature by developing the first framework for designing policies for allocating housing to those experiencing homelessness whose characteristics align with human value judgements.
}

\paragraph{Notation.} Throughout this paper, vectors (matrices) are denoted by boldface lowercase (uppercase) letters. The $k$th element of a vector ${\bm x} \in \mathbb R^n$ ($k \leq n$) is denoted by ${\bm x}_k$. We let $\1$ (resp.\ $\1_i$) denote the vector of all ones (resp.\ the $i$th basis vector) of appropriate dimension. With a slight abuse of notation, we may use the maximum and minimum operators even when the optimum may not be attained; in such cases, the operators should be understood as suprema and infima, respectively. Finally, for a logical expression $E$, we define the indicator function $\I{E}$ as $\I{E}:=1$ if $E$ is true and 0 otherwise.

\paragraph{Structure of the Paper.} The remainder of the paper is organized as follows. \newpv{Section~\ref{sec:related_work} introduces related literature.} Section~\ref{sec:model} formalizes the model of our recommendation system and the preference model that underlies our approach. \newpv{Sections~\ref{sec:maxminutility} and~\ref{sec:minimaxregret} study the max-min utility and min-max regret active preference elicitation problems, respectively. Section~\ref{sec:numericals} applies our framework to the problem of designing housing allocation policies that align with human value judgements} and Section~\ref{sec:conclusion} concludes. The proofs of all statements can be found in the Electronic Companion to the paper.


\section{\newpv{Review of Related Work}}
\label{sec:related_work}

\newpv{Our work lies at the intersection of three literature streams: conjoint analysis in marketing, preference elicitation in AI, and stochastic and robust optimization with DDID in OR. It is also connected to preference robust optimization, multi-armed bandits, recommender systems, assortment planning, and dynamic pricing.} 


\textbf{\textsf{\newpv{Conjoint Analysis.}}} 
\newpv{Most of the literature on conjoint analysis is focused on nonadaptive (offline) elicitation, relies on probability models for learning from discrete choices~\citep{McFadden1974,McFadden2001,Stern1990}, and builds on the field of experimental design. Queries are chosen to maximize a proxy of ``precision'' of the parameter estimates produced from the answers to the queries, see e.g.,~\cite{Arora2001,Huber1996,Kuhfeld1994,Kuhfeld2010}.}

\newpv{Several works consider adaptive (online) elicitation and model uncertainty in the utility by means of a polyhedron. This set, that we refer to as \emph{uncertainty set} following robust optimization terminology, is progressively ``shrunk'' as answers to queries are observed. \cite{Toubia_2003} considers metric paired-comparisons while \cite{Toubia_2004,Toubia_2007} and \cite{OHair_LearningPreferences} focus on choice-based conjoint queries. These approaches leverage optimization techniques to select queries \emph{one at a time} (greedily) in a way that will reduce the size of the resulting uncertainty set as much as possible. Inconsistencies are handled either by considering distributions supported by mixtures of polyhedra~\citep{Toubia_2007} or by constraining the number of inconsistencies~\citep{OHair_LearningPreferences} and recommendations using the analytic center of the uncertainty set~\citep{Toubia_2003,Toubia_2004,Toubia_2007} or by maximizing worst-case utility~\citep{OHair_LearningPreferences}.}

\newpv{
Modeling uncertainty using uncertainty sets as in the \emph{polyhedral methods} is particularly attractive when seeking robust solutions--it circumvents the need for estimating distributions and for doing hard Bayesian updates. We thus build upon this line of research and model uncertainty using a polyhedron. However, our work has many differences. First, none of the polyhedral methods apply to the offline setting. Second, our framework explicitly models and solves the elicitation problem as a sequential decision-making problem. Third, our model of inconsistencies is different: it restricts errors in the responses to lie in a polyhedral set. Compared with the approach of~\cite{Toubia_2007}, it does not result in a combinatorial explosion in the number of polyhedra. Relative to~\cite{OHair_LearningPreferences}, it ensures that response errors that are unlikely will be excluded, resulting in a more realistic, less conservative model. These differences ensure that our framework will result in queries with better performance in terms of worst-case utility or worst-case regret.} 

\newpv{Recently, \cite{Vielma2019} proposed a new approach to adaptive choice-based conjoint analysis that is based on Bayesian updating but that preserves some of the attractive interpretability properties of the polyhedral method. As this approach decouples elicitation and recommendation while also taking a probabilistic perspective, it is not well suited to the problems we study.}


\textbf{\textsf{\newpv{Preference Elicitation in AI.}}} \newpv{Online preference elicitation has received considerable attention in AI. Some of this work takes a Bayesian perspective and relies on heuristics to identify queries \citep{Chajewwska00:Making,Boutilier02:POMDP,Eric08:Active,Zhao18:Cost}. Others model uncertainty using polyhedra and focus on gamble queries~\citep{Wang03:Incremental,Boutilier06:Constraint}. Utility is usually elicited using heuristics and recommendations are made using the min-max regret criterion~\citep{Viappiani2013}. In our work, we also investigate recommendations that minimize worst-case regret but focus on choice-based conjoint queries and study also offline elicitation. The AI community has investigated many applications of preference elicitation, see e.g.,~\cite{Sandholm06:Preference,Braziunas10:Assessing,Freedman18:Adapting,Benade17:Preference}.}


\textbf{\textsf{\newpv{Stochastic and Robust Optimization with DDID.}}} 
\newpv{Stochastic programming (SP) with DDID dates back to the works of~\cite{Jonsbraten_DecDepRandElmts} and~\cite{Jonsbraten_thesis}. Most of this literature assumes that the uncertain parameters are discretely distributed and models uncertainty by means of a finite scenario tree whose structure depends on the binary decisions that determine the time of information discovery~\citep{GoelGrossmanGasFields,GoelGrossman_ClassStochastic_DDU,GoelGrossman_NovelBB_GasFields,Colvin_Pharmaceutical,GuptaGrossman_SolutionStrategies}. These methods do not apply to our context which involves also real valued uncertain parameters and we wish to optimize decisions robustly.}



\newpv{Our work most closely relates to RO with DDID~\citep{DDI_VKB,vayanos_ROInfoDiscovery}. \cite{DDI_VKB} propose a decision-rule based approximation approach that relies on a pre-partitioning of the support of the uncertain parameters. \cite{vayanos_ROInfoDiscovery} study a solution method based on the $K$-adaptability approximation. They also investigate an active preference elicitation problem in the context of the US KAS which is however very different from the one we consider in the present paper. Indeed, queries in~\cite{vayanos_ROInfoDiscovery} correspond to a single item and the user is asked to quantify ``how much'' they like the item (on a continuous scale from~0 to~1). Thus, uncertainty is real-valued and moderately sized whereas our problem presents an exponential search space and an exponential number of binary contingencies, requiring a fundamentally different solution approach. In addition, our solution to offline elicitation is exact, whereas their method is approximate.}

\newpv{The problems we study also relate to the works of~\cite{Nohadani2017} and~\cite{Trichakis2018monitoring} which focus on robust monitoring problems. \cite{Nohadani2017} consider problems with time-dependent uncertainty sets that capture increased uncertainty as information ages and observations ``reset uncertainty''. \cite{Trichakis2018monitoring} consider a system with evolving state that can be stopped at any time by the decision-maker to yield a time dependent reward. The decision-maker can choose, in conjunction with the stopping policy, a limited number of monitoring times when to observe the state of the system. Contrary to these problems the state of the system does not change in our case, there is a large number of parameters that we can choose to observe, uncertainty is discrete, and there is a very large number of possible scenarios.} 





\textbf{\textsf{\newpv{Preference Robust Optimization.}}} \newpv{There is very rich literature connecting stochastic and robust optimization and risk preferences~\citep{Dentcheva2004,RiskMeasuresUncertaintySets,Brown2009,Chen2009,Brown2012,Lam2013}. Our paper most closely relates to preference robust optimization. This literature assumes that partial information on preferences is available (e.g., through previous surveys) and focuses on making robust recommendations accounting for ambiguity in risk preferences while capturing sophisticated choice models~\citep{Armbruster_2015,Hu2015,Hu2018,Delage2018,Delage2018b,Haskell2018}. Differing from this work we study the problem of \emph{selecting queries} and focus on learning trade-offs between attributes rather attitudes towards risk.}



\textbf{\textsf{\newpv{Multi-Armed Bandits.}}} \newpv{Multi-armed bandits (MAB) offer a framework for modeling problems where information is acquired over time as a by-product of decisions~\citep{Gittins1989}. Extant models account for capacity constraints~\citep{Wu2015}, contextual information~\citep{pmlr-v9-lu10a,Badanidiyuru2014}, and provide robust formulations~\citep{Auer1995,Auer2003,Syrgkanis2016,Syrgkanis2016b,Lykouris2018}. 
The key feature that differentiates our setting is that we do not obtain a reward as a byproduct of the decisions and that the number of arms is exponential if we model our problem from this lens. Thus MAB are not suited to tackle our problem.}


\textbf{\textsf{\newpv{Recommender Systems.}}} \newpv{Research on recommender systems aims to help users navigate large product spaces and has been motivated from web applications, e.g., to generate playlists for video and music services~\citep{Ladeira2019}, or content for social media platforms~\citep{Gupta2013b}. They aim to \emph{filter} content to provide users with an experience that aligns with their preferences while still being rich/diverse enough, see e.g., \cite{McSherry2002,Reilly2005}. These methods generally rely on a lot more data than we can collect and are not well suited for high stakes settings.}


\textbf{\textsf{\newpv{Assortment Planning.}}} \newpv{Static assortment optimization with known choice behavior has been an active area of research \citep{Ryzin1999,Mahajan2001}. Dynamic assortment optimization which adaptively learns the customers’ choice behavior is receiving increasing attention. Most of this work does not consider the features of the items~\citep{Caro2007,Rusmevichientong_2010,Saure2013,Agrawal2019} or focuses on the pricing problem~\citep{Roth2016,Roth2020}. Our work most closely relates to assortment planning with contextual models~\citep{Cheung2017,Chen2020,Oh2019}.  \cite{DBLP:journals/corr/abs-1910-04183} consider contextual information while also taking a robust viewpoint. However, their work is based on revealed preferences and focused on multiple customers.}


\textbf{\textsf{\newpv{Dynamic Pricing.}}} \newpv{
The literature on dynamic pricing with learning studies pricing in settings where the demand function is unknown and can be learned as a byproduct of decisions~\citep{Kleinberg_IEEESFCS}.
A large body of literature has studied this topic including both parametric \citep{Araman_DP,Broder_OR,Keskin_BDP,Chen_OR, denBoer_MS,Besbes_WP} and non-parametric~\citep{Besbes_DP_OR,Keskin_DP} approaches. The literature also includes models that use a robust optimization approach~\citep{BV_dynamic_pricing,Cohen2020}. Our work is different in that we need to learn preferences from pairwise comparison queries and stated preference surveys.}



\section{Model}
\label{sec:model}


In this section, we define \newpv{the building blocks of our framework that will be used to compute optimal active preference elicitation strategies in Sections~\ref{sec:maxminutility} \newpv{and} \ref{sec:minimaxregret}.}


\subsection{Items, Query and Recommendation Sets, and User Preference Model}
\label{sec:items_preferences_model}

\paragraph{Items.} We assume that, when choosing one product (e.g., policy) over another, a user is basing their decision on the attributes of the options, \newpv{as in conjoint analysis}. Thus, we characterize each item~${\bm x}$ by its~$J$ attributes and model it as a point in a~$J$-dimensional Cartesian space, i.e., ${\bm x} = ({\bm x}_1,\ldots,{\bm x}_J)$, where ${\bm x}_j \in \mathbb R$, $j=1,\ldots,J$, denotes the $j$th attribute of item~${\bm x} \in \reals^J$. In this work, we are motivated by settings where~$J$ is high-dimensional (in the order of 10 or 20 for example). We denote the universe of all item configurations that are possible to produce by~$\sets X \subseteq \mathbb R^J$.




\paragraph{Query Set.} We let $\sets Q \subseteq \sets X$ denote the query set, i.e., the set of items that the recommender system can use to build queries. We 
%
assume that query set $\sets Q$ has finite cardinality $I \geq 2$ and index items in this set by $i \in \sets I := \{1,\ldots,I\}$. We denote the $i$th item by ${\bm x}^i \in \reals^J$, i.e., $\sets Q := \left\{{\bm x}^i  \mid i \in \mathcal I \right\}$.
%
This assumption is very common and also natural\newpv{--e.g., in policy-making settings such as during the OPTN Kidney Transplantation Committee Meetings, there is usually a finite number of candidate policies that are being compared.}

\paragraph{Recommendation Set.} We denote by $\sets R \subseteq \sets X$ the {recommendation set}, i.e., the set of items from which the system can draw to make a recommendation \newpv{or equivalently the feasible set of the optimization problem that the decision-maker is interested in solving}. In general, $\sets Q$ and $\sets R$ need not coincide. For example, we may be able to {design a new item}, that is not currently available or we may be unable to recommend a product that is now out-of-stock. \newpv{In our work, the set $\sets R$ may have a very general structure. For example, it may take the form of the feasible set of a shortest path problem, a knapsack problem, an assignment problem, or a minimum spanning tree problem.}

\paragraph{User Preference Model.} We assume that the user's preferences can be represented with a linear utility function $u : \mathcal X \rightarrow \reals$ defined through $u({\bm x}):= {\bm u}^\top {\bm x} + \gamma({\bm x})$, where ${\bm u} \in \reals^J$ is a vector representing the relative importance they assign to product attributes and $\gamma({\bm x})$ is an idiosyncratic shock to utility. We make the following assumption regarding the vector ${\bm u}$.
\begin{assumption}[Polyhedral Uncertainty]
The vector $\bm{u}$ is an element of the uncertainty set $\mathcal U^0 \subseteq \reals^J$ defined through $\mathcal U^0:=\{ {\bm u} \in \mathbb{R}^J \mid {\bm B} {\bm u} \geq {\bm b} \}$ for some matrix ${\bm B} \in \reals^{M \times J}$ and vector ${\bm b} \in \reals^{M}$. Moreover, $\mathcal U^0$ is non-empty, full-dimensional, and bounded. 
\label{ass:linear_utility}
\end{assumption}
The assumption \newpv{that $u$ is linear} and that prior information on~$\bm u$ can be encoded using linear inequalities is standard in the literature on preference elicitation and polyhedral methods to conjoint analysis, see Section~\ref{sec:related_work}. \newpv{Our utility model is however different than the ones used in this line of work in that we explicitly assume that the utility function is affected by idiosyncratic shocks, in line with several popular choice models from the literature, see e.g., \cite{GreenBook2012}. This modelling choice will allow us to have a polyhedral model for the inconsistencies, see Section~\ref{sec:queries_model}.}

Assumption~\ref{ass:linear_utility} can be relaxed to allow that $\mathcal U^0$ presents equality constraints and that it possesses an inner point (rather than being full-dimensional) at the cost of complicating our proofs.

\begin{example}
If no prior information is available on $\bm u$, since the utility coefficients can be scaled by a constant without affecting the preference ordering over items, one may use $\mathcal U^0= [-1,1]^J$, see e.g.,~\cite{OHair_LearningPreferences}. 
\end{example}

\begin{example}
Alternatively, one may set $\mathcal U^0 = \{ {\bm u} \in \reals_+^J : \1^\top {\bm u}=1 \}$ as in~\cite{Toubia_2003}. In that case, the coefficients ${\bm u}$ can be viewed as partworth utilities, i.e., numerical scores that measure how much each feature influences the user's decision to make that choice.
\label{ex:conjoint_uncertainty_set}
\end{example}


\begin{example}
If item $\overline{\bm x}$ is a benchmark, one may normalize utilities relative to that item, by letting
$\mathcal U^0 = \{ {\bm u} \in [-1,1]^J \; : \; {\bm u}^\top \overline{\bm x} = 1  \}$.
\end{example}



\subsection{Elicitation through Comparison Queries}
\label{sec:queries_model}

Before recommending an item from the set~$\mathcal R$, the recommender system has the opportunity to make~$K$ queries, indexed in the set $\sets K := \{1,\ldots,K\}$, to the user. These queries may enable the system to gain information about~$\bm u$, thereby improving the quality of the recommendation. Each query takes the form of a comparison between two items in the query set, thus being an element of $\sets C := \{ ( i, i' ) \mid i \in \sets I, \; i' \in \sets I, \; i < i' \}$. 
We let ${\bm \iota}^\kappa := ( {\bm \iota}^\kappa_{1} ,{\bm \iota}^\kappa_2 )\in \sets C$ denote the $\kappa$th query, $\kappa \in \sets K$. Thus, ${\bm \iota}^\kappa_1$ and ${\bm \iota}^\kappa_2$ denote the indices of the first and second items in the $\kappa$th query, respectively: query $\kappa$ asks the user to compare ${\bm x}^{{\bm \iota}^\kappa_1}$ and ${\bm x}^{{\bm \iota}^\kappa_2}$. The user must choose one of three possible answers in response to query $\kappa$: \emph{(a)} ``I strictly prefer ${\bm x}^{{\bm \iota}^\kappa_1}$ to ${\bm x}^{{\bm \iota}^\kappa_2}$'' (i.e., ${\bm x}^{{\bm \iota}^\kappa_1} \succ {\bm x}^{{\bm \iota}^\kappa_2}$); \emph{(b)} ``I am indifferent between ${\bm x}^{{\bm \iota}^\kappa_1}$ and ${\bm x}^{{\bm \iota}^\kappa_2}$'' (i.e., ${\bm x}^{{\bm \iota}^\kappa_1} \sim {\bm x}^{{\bm \iota}^\kappa_2}$); or \emph{(c)} ``I strictly prefer ${\bm x}^{{\bm \iota}^\kappa_2}$ to ${\bm x}^{{\bm \iota}^\kappa_1}$''  (i.e., ${\bm x}^{{\bm \iota}^\kappa_1} \prec {\bm x}^{{\bm \iota}^\kappa_2}$). We associate each possible answer to query $\kappa$ to a response scenario ${\bm s}_\kappa \in \sets S:=\{-1,0,1\}$ such that
$$
{\bm s}_\kappa =\begin{cases}
\;\; 1 & \text{ if }  {\bm x}^{{\bm \iota}^\kappa_1} \succ {\bm x}^{{\bm \iota}^\kappa_2} \\ 
\;\; 0 & \text{ if } {\bm x}^{{\bm \iota}^\kappa_1} \sim {\bm x}^{{\bm \iota}^\kappa_2}  \\
-1 &\text{ else}.
\end{cases} 
$$

The information obtained on $\bm u$ depends on the answer to query $\kappa$, i.e., on the response scenario. Since the utility of any given item is affected by idiosyncratic noise, the answers given by the user may be \emph{inconsistent} in the sense that ${\bm x} \succ {\bm y}$ and $u({\bm x}) \leq u({\bm y})$; ${\bm x} \prec {\bm y}$ and $u({\bm x}) \geq u({\bm y})$; or ${\bm x} \sim {\bm y}$ and $u({\bm x}) \neq u({\bm y})$. Naturally, such inconsistencies may not occur if the user is perfectly rational, see e.g., \cite{Arrow_Book1963}. \newpv{Motivated by our utility model, see Section~\ref{sec:items_preferences_model}, we impose a different linear constraint on $\bm u$ depending on the response:
\emph{(a)} If ${\bm s}_\kappa=1$, then $\bm u^\top {\bm x}^{{\bm \iota}^\kappa_1} - \bm u^\top {\bm x}^{{\bm \iota}^\kappa_2} > - {\bm \epsilon}_\kappa$; \emph{(b)} If ${\bm s}_\kappa=0$, then $|\bm u^\top {\bm x}^{{\bm \iota}^\kappa_1} - \bm u^\top {\bm x}^{{\bm \iota}^\kappa_2} | \leq {\bm \epsilon}_\kappa$; and \emph{(c)} If ${\bm s}_\kappa=-1$, then $\bm u^\top {\bm x}^{{\bm \iota}^\kappa_1} -  \bm u^\top {\bm x}^{{\bm \iota}^\kappa_2} < {\bm \epsilon}_\kappa$, where ${\bm \epsilon}:=\{ {\bm \epsilon}_\kappa\}_{\kappa \in \sets K}$ is a vector of inconsistencies. We make the following assumption about ${\bm \epsilon}$.}
\begin{assumption}[Polyhedral Inconsistencies]
The user inconsistencies across $K$ queries are assumed to lie in the uncertainty set $\sets E_\Gamma:=\{ {\bm \epsilon} \in \reals_+^K \; : \; \sum_{\kappa \in \sets K} {\bm \epsilon}_\kappa \leq \Gamma \}$ for some (known) parameter $\Gamma \in \reals_+$, referred to as the budget of uncertainty.
\label{ass:noisy_preferences}
\end{assumption}
\newpv{
\begin{example}[Choice of $\Gamma$]
Standard approaches in RO can be used to choose~$\Gamma$~\citep{BenTal_Book,TheoryApplicationRO}. For example, if the idiosyncratic shocks in utility are assumed to be normal with mean zero and standard deviation $\sigma$, then $\sum_{\kappa \in \sets K} {\bm \epsilon}_\kappa$ is normal with standard deviation $2 \sigma \sqrt{ K }$ and setting $\Gamma = 2\sigma \sqrt{K} \rm{erf}^{-1}(2p-1)$ ensures that ${\bm \epsilon} \in \sets E_\Gamma$ with probability $p$.
\label{ex:choice_of_gamma}
\end{example}
}
\newpv{Modeling inconsistencies as lying in a polyhedron sets us apart from the literature which either uses mixtures of polyhedra for the support of ${\bm u}$~\citep{Toubia_2007}, or integer programming to bound the number of incorrect responses and big-$M$ constants to toggle the sign of the constraints added to $\sets U_0$~\citep{OHair_LearningPreferences}. In particular, our modeling paradigm is more interpretable than the approach of~\cite{Toubia_2007}. In addition, it is less conservative that the framework of~\cite{OHair_LearningPreferences} since our model will only allow inconsistent responses when the user ``marginally'' prefers one item over the other. Moreover, modeling inconsistencies in this way allows us to use duality theory to solve our problems in Sections~\ref{sec:maxminutility} and~\ref{sec:minimaxregret}.

}

\begin{figure}[h!]
    \centering
    \includegraphics[width=0.35\textwidth]{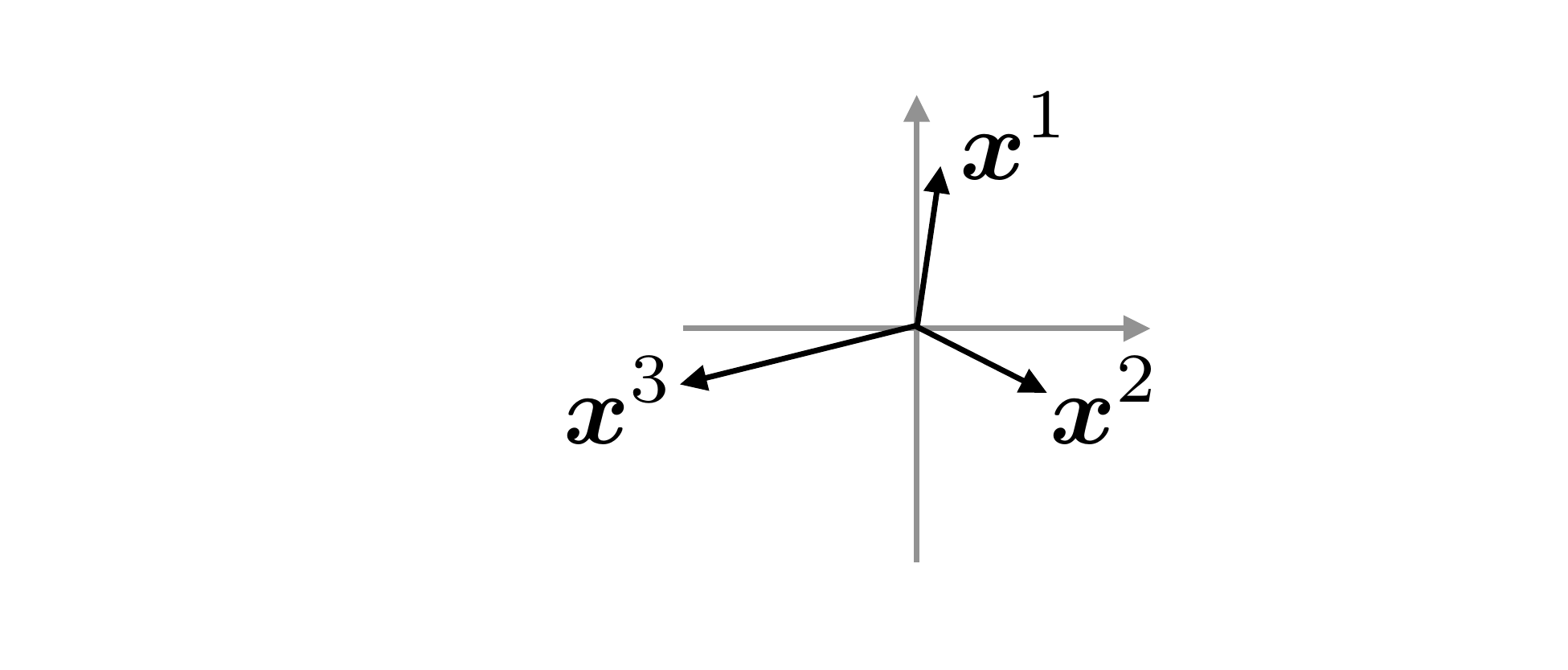} \hspace{0.5cm}
    \includegraphics[width=0.35\textwidth]{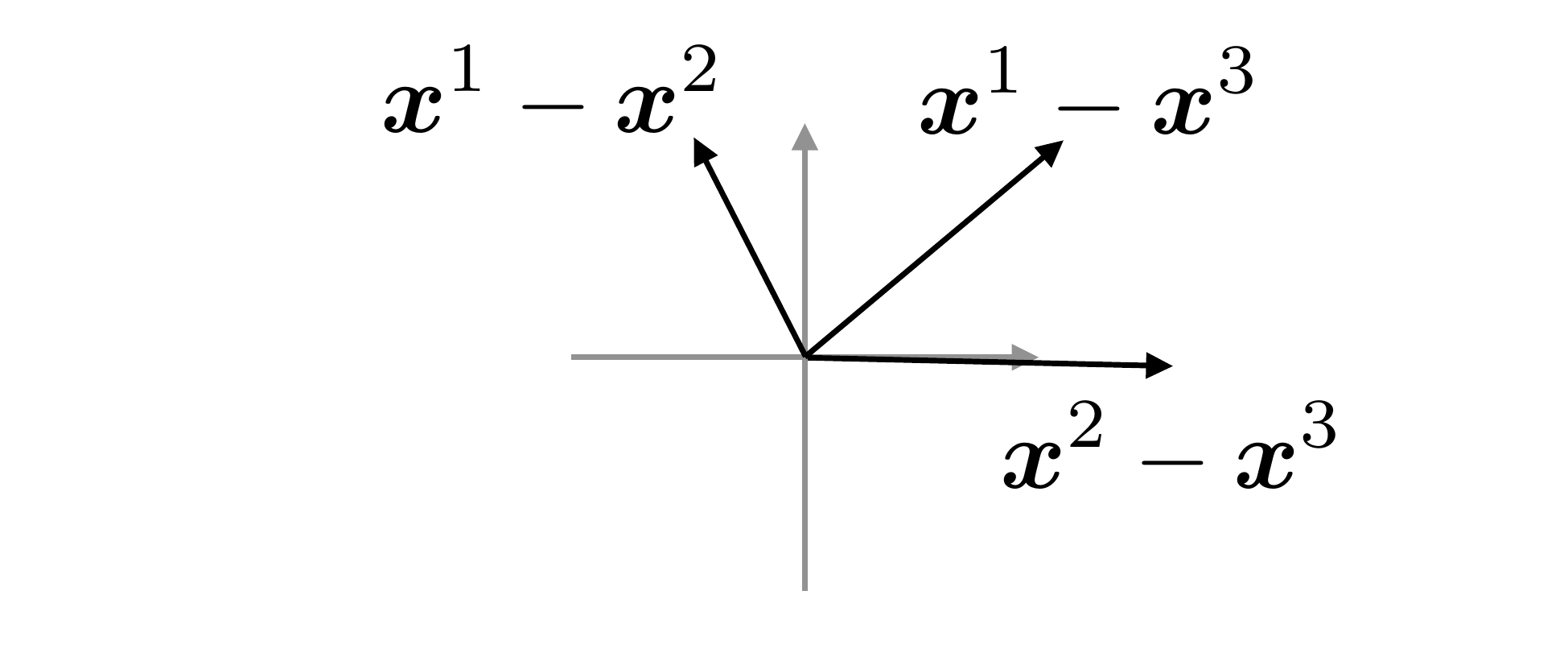}\\
    \includegraphics[width=0.35\textwidth]{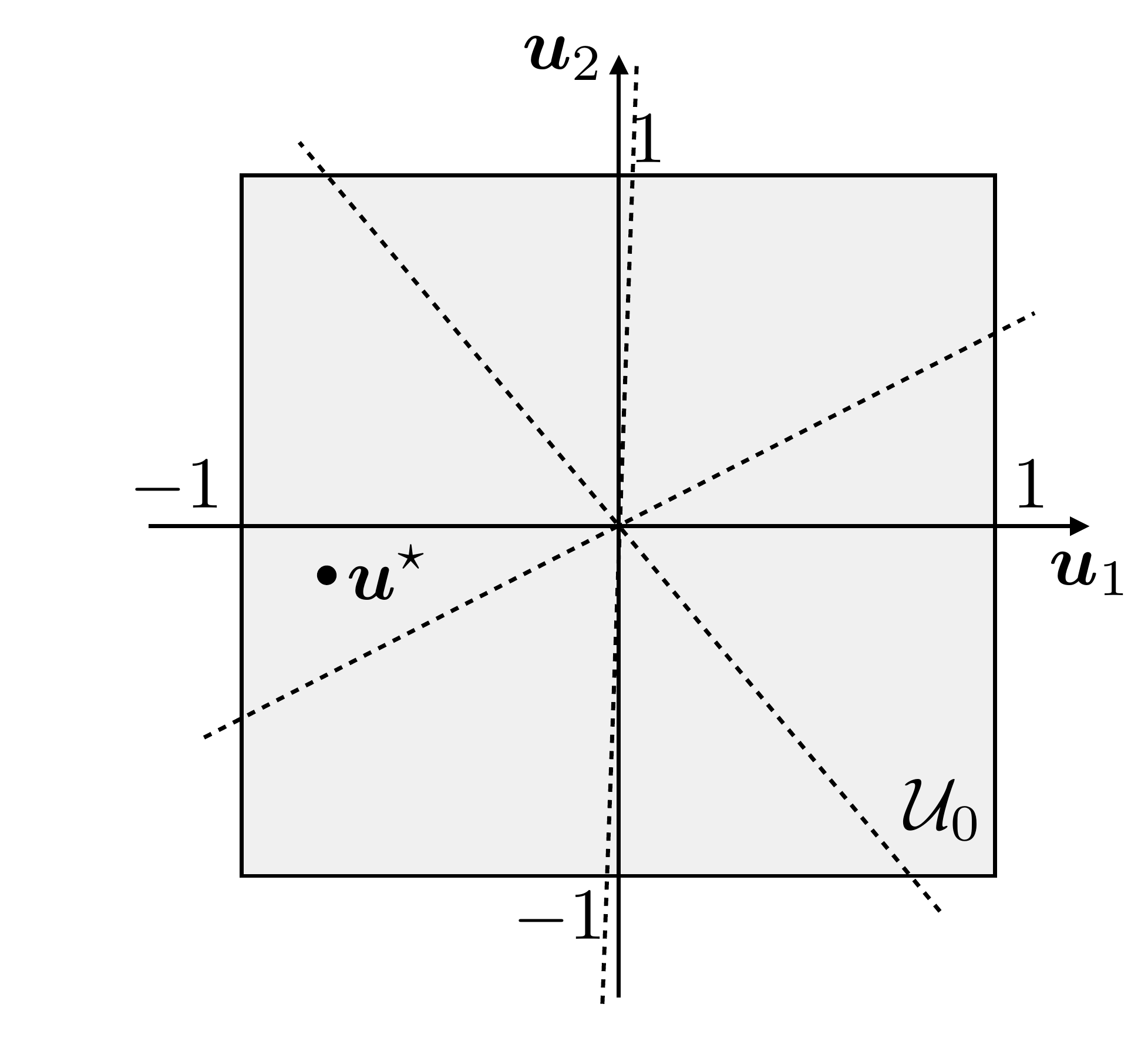}
    \includegraphics[width=0.35\textwidth]{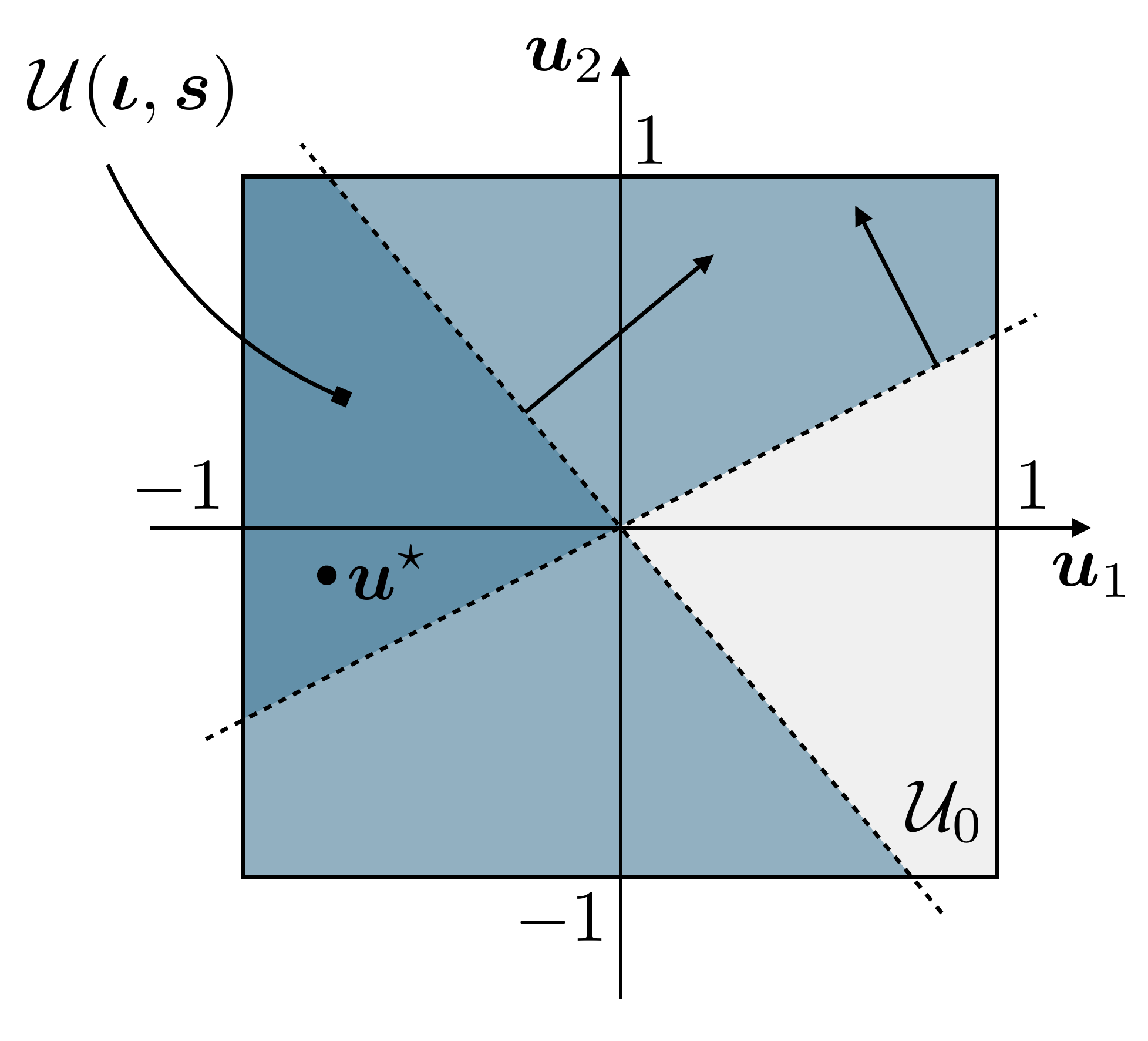}\\ 
    \includegraphics[width=0.35\textwidth]{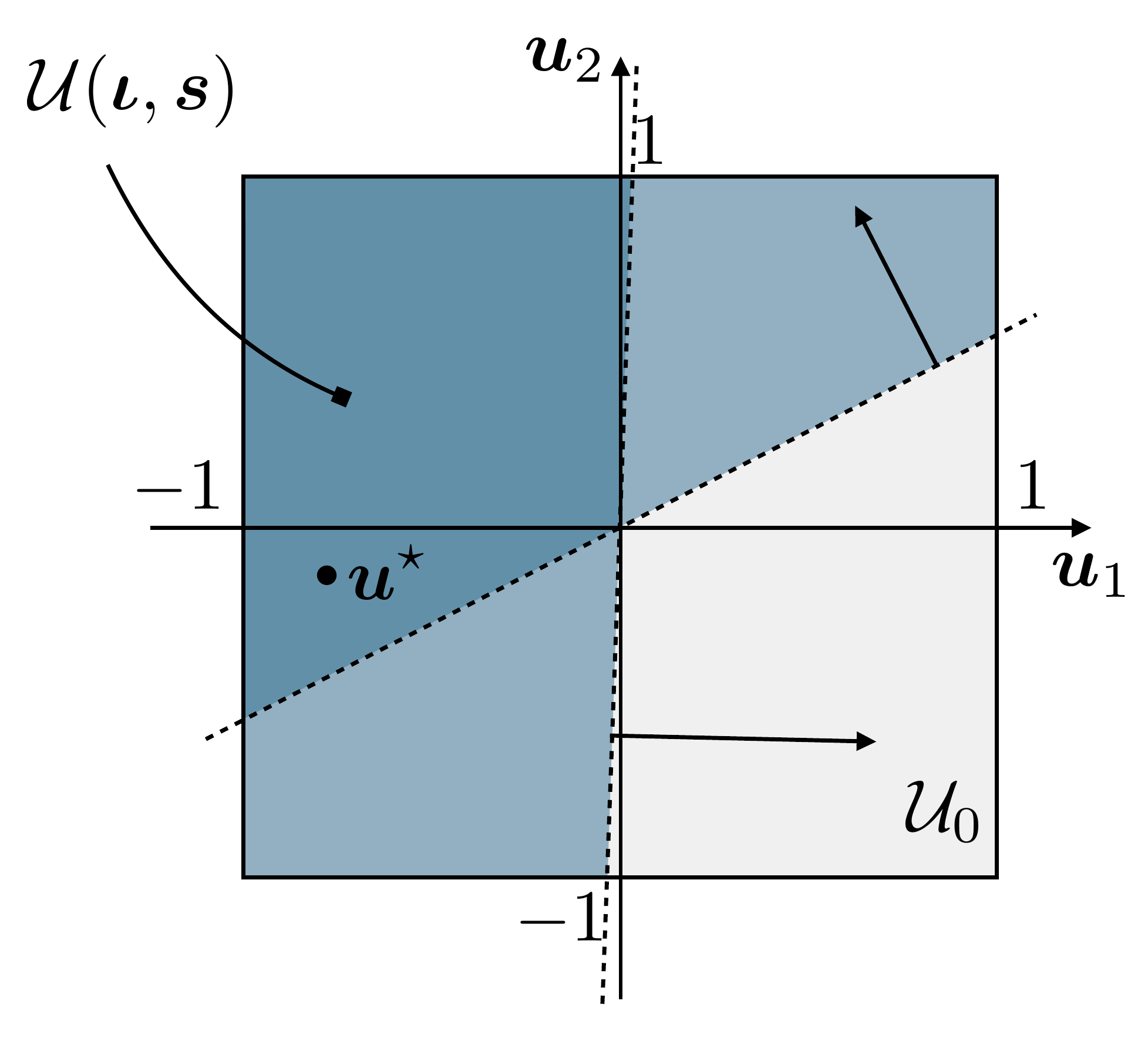}
    \includegraphics[width=0.35\textwidth]{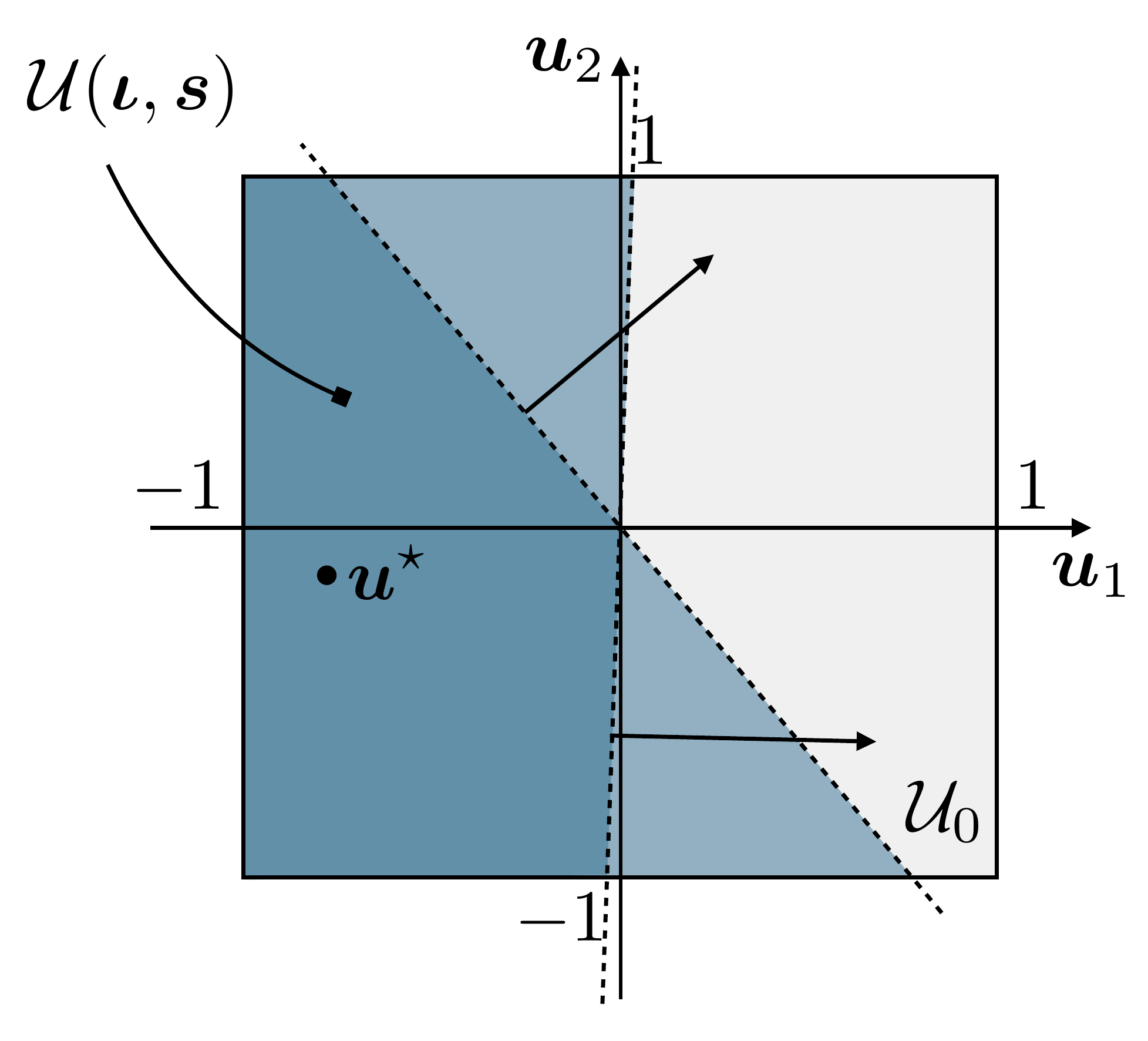}\\ 
    \caption{Illustration of the uncertainty set update procedure when there are $I=3$ items, each with $J=2$ features, $\Gamma=0$, the ``true'' (but unknown) utility vector is ${\bm u}^\star$, and the system is allowed to ask $K=2$ questions before making a recommendation. The first row shows the three items ${\bm x}^1$, ${\bm x}^2$, and ${\bm x}^3$ in $\reals^2$ (L) and the vectors ${\bm x}^1-{\bm x}^2$, ${\bm x}^1-{\bm x}^3$, and ${\bm x}^2-{\bm x}^3$ associated with each of the comparison queries $({\bm x}^1,{\bm x}^2)$, $({\bm x}^1,{\bm x}^3)$, and $({\bm x}^2,{\bm x}^3)$, respectively~(R). The left figure on the second row shows the initial uncertainty set~$\sets U^0$, the vector ${\bm u}^\star$, and the hyperplanes associated with each of the queries. The remaining three figures show the uncertainty set $\sets U({\bm \iota},{\bm s})$ updated in response to the three different pairs of queries $\{ ({\bm x}^1,{\bm x}^2) , ({\bm x}^1,{\bm x}^3) \}$ (row 2, R), $\{ ({\bm x}^1,{\bm x}^2) , ({\bm x}^2,{\bm x}^3) \}$ (row 3, L), and $\{ ({\bm x}^1,{\bm x}^3) , ({\bm x}^2,{\bm x}^3) \}$ (row 3, R). \newpv{On these three figures, the arrows represent the normal vectors to the query hyperplanes}. Note that the uncertainty set changes depending on the queries asked and on the answers given, which in turn depend on the underlying vector ${\bm u}^\star$ (whose value is unknown by the recommender system).}
    \label{fig:preference_elicitation}
\end{figure}

We collect the answers to each query in the vector ${\bm s} := \{ {\bm s}_\kappa \}_{\kappa \in \sets K}$ and accordingly let ${\bm \iota} := \{ ({\bm \iota}^\kappa_1, {\bm \iota}^\kappa_2) \}_{\kappa \in \sets K}$. After all~$K$ queries have been made and the responses to the queries are observed, we can update $\sets U^0$ as:
\newpv{
\begin{equation*}
\mathcal U_\Gamma({\bm \iota},{\bm s}) : = 
\left\{
\begin{array}{ccll}  
{\bm u} \in \mathcal U^0  & : &  \exists {\bm \epsilon} \in \sets E_\Gamma \text{ such that} \\
&& {\bm u}^\top ({\bm x}^{{\bm \iota}^\kappa_1} - {\bm x}^{{\bm \iota}^\kappa_2}) \; > \; -{\bm \epsilon}_\kappa &  \quad  \forall \kappa \in \sets K \; : \; {\bm s}_\kappa=1 \\
&& | {\bm u}^\top ( {\bm x}^{{\bm \iota}^\kappa_1} - {\bm x}^{{\bm \iota}^\kappa_2} ) | \; \leq \; {\bm \epsilon}_\kappa  & \quad  \forall \kappa \in \sets K \; : \; {\bm s}_\kappa=0 \\
&& {\bm u}^\top ( {\bm x}^{{\bm \iota}^\kappa_1} - {\bm x}^{{\bm \iota}^\kappa_2} )\; < \; {\bm \epsilon}_\kappa  &  \quad  \forall \kappa \in \sets K \; : \; {\bm s}_\kappa=-1
\end{array}
\right\}.
\label{eq:uncertainty_set_updated}
\end{equation*}}
We emphasize that after the $K$ comparisons, the set of possible values for ${\bm u}$ depends both on~${\bm \iota}$, the queries made, and on~${\bm s}$, the answers given by the user. \newpv{This idea is illustrated on Figure~\ref{fig:preference_elicitation}.}

\subsection{Robust Recommendations}
\label{sec:recommendation_model}

After~$K$ queries have been made and the answers to these queries have been observed, the recommendation system needs to select an item from the set~$\sets R$ to recommend albeit the coefficients $\bm u$ are only known to belong to the uncertainty set~\newpv{$\sets U_\Gamma({\bm \iota},{\bm s})$}. We thus propose to provide recommendations that are robust to all possible realization of $\bm u$ in the set~\newpv{$\sets U_\Gamma({\bm \iota},{\bm s})$} and investigate two notions of robustness that are popular in the literature: \newpv{max-min utility and min-max regret}.

\paragraph{Maximizing Worst-Case Utility.} Given uncertainty in the utility function coefficients, it is natural to seek recommendations that will maximize worst-case utility. Mathematically, given the sequences, ${\bm \iota}$ and ${\bm s}$, the recommender system offers \newpv{an} item \newpv{${\bm x}$ that} solves the problem
\begin{equation}
\label{eq:riskaverse_recommendation}
\tag{\Rrisk}
\maximize_{{\bm x} \in \sets R} \quad \min_{{\bm u} \in \sets U_\Gamma({\bm \iota},{\bm s})} \quad {\bm u}^\top {\bm x}. 
\end{equation}

\paragraph{Minimizing Worst-Case Regret.} An alternative is to make recommendations according to the worst-case absolute regret criterion which compares the performance of the decision taken relative to the performance of the best decision that should have been taken \emph{in hindsight}, after all uncertain parameters are revealed, see e.g., \cite{Savage51:Theory}. The minimization of worst-case absolute regret is often believed to mitigate the conservatism of classical robust optimization (which maximizes worst-case utility) and is thus attractive in practical applications. According to this criterion and given the sequences ${\bm \iota}$ and ${\bm s}$, the recommender system offers the item \newpv{that} solves the problem
\begin{equation}
\label{eq:regretaverse_recommendation}
\tag{\Rregret}
\minimize_{{\bm x} \in \sets R} \quad \max_{{\bm u} \in \sets U_\Gamma({\bm \iota},{\bm s})} \quad \left\{ \max_{ \; {\bm x}' \in \sets R } \; \; {\bm u}^\top {\bm x}' - {\bm u}^\top {\bm x} \; \right\}. 
\end{equation}


\newpv{We note that for some choices of $\bm \iota$ and $\bm s$, the set \newpv{$\sets U_\Gamma({\bm \iota},{\bm s})$} may be open and even empty and emphasize that t}he \newpv{sets of} optimal solutions to~\eqref{eq:riskaverse_recommendation} and~\eqref{eq:regretaverse_recommendation} depend on both ${\bm \iota}$ and ${\bm s}$. While the system has no control over the responses ${\bm s}$ it \emph{can} select the queries ${\bm \iota} \in \sets C^K$ to improve the quality of the recommendation. \newpv{This observation will motivate our formulations in Sections~\ref{sec:maxminutility} and~\ref{sec:minimaxregret}.} 

\section{Active Preference Elicitation with the Max-Min Utility Decision Criterion}
\label{sec:maxminutility}

\newpv{In this section, we study the offline and online active preference elicitation problems under the max-min utility decision criterion and evaluate their performance on synthetic data.}


\subsection{Offline Active Preference Elicitation with the Max-Min Utility Decision Criterion}
\label{sec:offline_maxminutility}

In the offline problem with the max-min utility decision criterion, all~$K$ comparisons are precommitted \emph{in advance}, before any agent responses are observed \newpv{and the goal is to select robust queries that will result in highest worst-case utility of the recommendation}. \newpv{Offline elicitation is relevant when interactions with the user do not involve a computer (e.g., paper-based questionnaires) or if all users must be presented the same queries (e.g., in a controlled study).} It is also a building block for solving the online preference elicitation problem, see Section~\ref{sec:online_maxminutility}.



\subsubsection{Problem Formulation \& Complexity Analysis}
\label{sec:offline_mmu_problem_formulation}

In the \newpv{offline problem with the max-min utility decision criterion,} the recommender system selects~$K$ queries~${\bm \iota}^\kappa \in \sets C$, $\kappa \in \sets K$, to ask the user. Subsequently, the user, who has a true (but unknown) utility vector~${\bm u}^\star$ from~$\sets U^0$, responds to the queries truthfully \newpv{(albeit potentially inconsistently, see Section~\ref{sec:queries_model})} by selecting answers~${\bm s}_\kappa$ to each query~${\bm \iota}^\kappa$ in a way that complies with their utility vector~${\bm u}^\star$ \newpv{and some choice of ${\bm \epsilon} \in \sets E_\Gamma$}. Note that the utility vector~${\bm u}^\star$ is not observable to the recommender system; in fact, the user themself is not aware of their vector~${\bm u}^\star$ (else, they would directly share it). Once the answers to the queries are observed, the recommender system can certify that~\newpv{${\bm u}^\star \in \sets U_\Gamma({\bm \iota},{\bm s})$} and solves problem~\eqref{eq:riskaverse_recommendation}.

Mathematically, the \newpv{offline problem with the max-min utility decision criterion} is expressible as the two-stage robust optimization problem with decision-dependent information discovery
\begin{equation}
\tag{$\offlinerisk{K}$}
\maximize_{{\bm \iota} \in \sets C^K} \quad \min_{{\bm s}\in \newpv{\sets S_\Gamma}({\bm \iota})} \quad \max_{{\bm x} \in \sets R} \quad \min_{ {\bm u} \in \newpv{\sets U_\Gamma}({\bm \iota} ,{\bm s})} \quad  {\bm u}^\top {\bm x},
\label{eq:offline_mmu}
\end{equation}
where
$
\newpv{\sets S_\Gamma}({\bm \iota}) := \left\{ {\bm s} \in \sets S^K \; \mid \; \newpv{\sets U_\Gamma({\bm \iota},{\bm s})} \neq \emptyset \right\}
$
denotes the set of all answers compatible with some~${\bm u} \in \sets U^0$ \newpv{and $\bm \epsilon \in \sets E_\Gamma$}. Indeed, the set~$\sets U_\Gamma({\bm \iota},{\bm s})$ is empty if and only if \newpv{the answers given by the user are not consistent with any ${\bm u} \in \sets U^0$ and no vector ${\bm \epsilon}\in \sets E_\Gamma$ can resolve the inconsistencies.} Since~$\sets U^0$ is non-empty by Assumption~\ref{ass:linear_utility}, \newpv{$\sets S_\Gamma({\bm \iota}) \neq \emptyset$} for all~${\bm \iota} \in \sets C^K$. Problem~\eqref{eq:offline_mmu} is a ``two-and-a-half'' stage \newpv{dynamic} robust optimization problem with decision-dependent uncertainty set. 


Problem~\eqref{eq:offline_mmu} is difficult to solve for many reasons. First, it is a max-min-max-min problem. Second, the uncertainty sets~\newpv{$\sets S_\Gamma({\bm \iota})$} and \newpv{$\sets U_\Gamma({\bm \iota},{\bm s})$} for the first and second decision-stages are decision-dependent. Third, the uncertainty set~$\sets U_\Gamma({\bm \iota},{\bm s})$ depends upon the first stage uncertainty~${\bm s}$ and is also open, making it difficult to derive computational solution approaches. Fourth, the uncertainty set $\sets S_\Gamma({\bm \iota})$ involves ``non-equalities''. The following lemma shows that Problem~\eqref{eq:offline_mmu} can be considerably simplified by eliminating the dependence of~\newpv{$\sets S_\Gamma({\bm \iota})$} on~${\bm \iota}$, by replacing the strict inequalities in the set~$\sets U_\Gamma({\bm \iota},{\bm s})$ by their loose counterparts, \newpv{and by eliminating indifferent responses from the response set}.

\newpv{
\begin{lemma}
Problem~\eqref{eq:offline_mmu} is equivalent to
\begin{equation}
\tag{$\widetildeofflinerisk{K}$}
\maximize_{{\bm \iota} \in \sets C^K} \quad \min_{{\bm s}\in \widetilde{\sets S}^K} \quad \max_{{\bm x} \in \sets R} \quad \min_{ {\bm u} \in \widetilde{\sets U}_\Gamma({\bm \iota} ,{\bm s}) } \quad  {\bm u}^\top {\bm x},
\label{eq:offline_mmu_4}
\end{equation}
where $\widetilde{\sets S} := \{-1,1\}$ and
\begin{equation*}
\widetilde{\sets U}_\Gamma({\bm \iota},{\bm s}) := 
\left\{
\begin{array}{ccll}  
{\bm u} \in \sets U^0 & : &  \exists {\bm \epsilon} \in \sets E_\Gamma \text{ such that} \\
&& {\bm u}^\top ({\bm x}^{{\bm \iota}^\kappa_1} - {\bm x}^{{\bm \iota}^\kappa_2}) \; \geq \; -{\bm \epsilon}_\kappa &  \quad  \forall \kappa \in \sets K \; : \; {\bm s}_\kappa=1 \\
%
%
&& {\bm u}^\top ( {\bm x}^{{\bm \iota}^\kappa_1} - {\bm x}^{{\bm \iota}^\kappa_2} ) \; \leq \; {\bm \epsilon}_\kappa  &  \quad  \forall \kappa \in \sets K \; : \; {\bm s}_\kappa=-1
\end{array}
\right\},
\end{equation*}
in the sense that the two problems have the same optimal objective value and the same sets of optimal queries.
\label{lem:offline_mmu_simplification}
\end{lemma}
The above result avoids us having to use ``little-$\epsilon$'' variables to convert the strict inequalities in the uncertainty set to loose inequalities as done e.g., by~\cite{OHair_LearningPreferences} in their procedure to calculate the analytic center of the uncertainty set.
}

Next, we show that Problem~\eqref{eq:offline_mmu_4} is $\mathcal{NP}$-hard when the recommendation set is discrete, which motivates the integer optimization based reformulation we provide in the next section.

\begin{theorem}
The following claims about the complexity of Problem~\eqref{eq:offline_mmu_4} hold true.
\begin{enumerate}[label={(\alph*)}]
    \item Suppose the recommendation set $\mathcal{R}$ is convex. Then, Problem~\eqref{eq:offline_mmu_4} is polynomially solvable. Moreover, it holds that Problems $(\widetildeofflinerisk{K})$ and $(\widetildeofflinerisk{0})$ are equivalent for all $K \in \naturals$, i.e., in the case of polyhedral recommendation set, in the worst-case, there is no benefit in asking any queries. In addition, the optimal objective values of Problems~\eqref{eq:offline_mmu_4} and~\eqref{eq:riskaverse_recommendation} coincide.
    %
    %
    \item Evaluating the objective function of Problem~\eqref{eq:offline_mmu_4} is $\mathcal{NP}$-hard even if the recommendation set $\sets R$ consists of only two elements.
\end{enumerate}
\label{thm:offline_mmu_complexity}
\end{theorem}
%
\newpv{The intuition behind claim \emph{(a)} is akin to that of the minimax theorem~\citep{Neumann_1928}. Indeed, under the conditions in the statement, $\sets R$ and $\widetilde{\sets U}_\Gamma({\bm \iota},{\bm s})$ are both convex and so observing~${\bm u} \in \widetilde{\sets U}_\Gamma({\bm \iota},{\bm s})$ before choosing~${\bm x}$ yields no benefit. Since this holds true for any choice of~${\bm \iota}$ and~${\bm s}$, there is no gain (in the worst-case) from asking queries.} Motivated by the complexity results above, we henceforth assume the recommendation set in Problem~\eqref{eq:offline_mmu_4} is discrete and MBLP representable. We provide a general, MBLP reformulation of Problem~\eqref{eq:offline_mmu_4} in Section~\ref{sec:offline_mmu_milp}, a column-and-constraint generation approach in Section~\ref{sec:offline_mmu_ccg}, and additional speed-up strategies in Section~\ref{sec:speedup_strategies}.

\subsubsection{Exact MBLP Reformulation}
\label{sec:offline_mmu_milp}


\newpv{While~\eqref{eq:offline_mmu_4} can in principle be solved by enumeration over ${\bm \iota}$ and ${\bm s}$, such an approach becomes computationally prohibitive as~$K$ grows. We thus propose instead to reformulate~\eqref{eq:offline_mmu_4} as an MBLP by noting that} we can exchange the order of the inner minimization and maximization problems, provided we index the recommendation by~${\bm s}$. \newpv{This observation combined with classical robust optimization techniques yields the following lemma.}

%
%
\begin{lemma}
Problem~\eqref{eq:offline_mmu_4} is equivalent to the following finite program
\newpv{\begin{equation}
\begin{array}{cl}
\maximize & \quad \tau \\
\subjectto  & \quad \tau \in \reals, \; {\bm \iota} \in \sets C^K \\
& \quad \!\!\left. \begin{array}{l}
{\bm x}^{\bm s} \in \sets R, \; {\bm \alpha}^{\bm s} \in \reals^K_+, \; {\bm \beta}^{\bm s} \in \reals_+^M ,\;  \mu^{\bm s} \in \reals_- \\
\tau  \;\; \leq \;\; {\bm b}^\top {\bm \beta}^{\bm s} + \Gamma \mu^{\bm s}  \\
\displaystyle \sum_{\kappa \in \sets K} \; {\bm s}_\kappa \; ({\bm x}^{{\bm \iota}^\kappa_1} - {\bm x}^{{\bm \iota}^\kappa_2}) \; {\bm \alpha}_\kappa^{\bm s} + {\bm B}^\top {\bm \beta}^{\bm s} \; = \; {\bm x}^{\bm s}  \\
 {\bm \alpha}^{\bm s} + \mu^{\bm s} \1 \leq {\bm 0}
\end{array} \quad \quad \right\} \quad \forall {\bm s}\in \widetilde{\sets S}^K,
\end{array}
\label{eq:offline_mmu_4_finite_program}
\end{equation}}
where ${\bm \iota}$ denotes the queries to make and ${\bm x}^{\bm s}$ the items to recommend in response scenario ${\bm s}\in \widetilde{\sets S}^K$.
\label{lem:offline_mmu_4_finite_program}
\end{lemma}
Problem~\eqref{eq:offline_mmu_4_finite_program} can be converted to an MBLP by encoding the choice of ${\bm \iota}^\kappa \in \sets C$, $\kappa \in \sets K$, using binary decisions ${\bm v}^\kappa \in \{0,1\}^I$ and ${\bm w}^\kappa \in \{0,1\}^I$ whose $i$th element is one if and only if item $i$ is the first (resp.\ second) item in query $\kappa$. \newpv{This result is summarized in the following theorem.} 
\newpv{\begin{theorem}
Problem~\eqref{eq:offline_mmu_4_finite_program} is equivalent to the following mixed-binary linear program
\begin{equation}
\begin{array}{cl}
\maximize & \quad \tau \\
\subjectto  & \quad \tau \in \reals, \; {\bm v}^\kappa, \; {\bm w}^\kappa \in \{0,1\}^I , \; \kappa \in \sets K  \\
& \quad {\bm \alpha}^{\bm s} \in \reals^K_+, \; {\bm \beta}^{\bm s} \in \reals_+^M ,\;  \mu^{\bm s} \in \reals_- ,\; {\bm x}^{\bm s} \in \sets R, \; {\bm s}\in \widetilde{\sets S}^K \\
& \quad \overline{\bm v}^{{\bm s \kappa}}, \; \overline{\bm w}^{{\bm s \kappa}} \in \reals^I_+, \;  {\bm s} \in \widetilde{\sets S}^K, \; \kappa \in \sets K \\
& \quad \!\!\left. \begin{array}{l}
\tau  \;\; \leq \;\; {\bm b}^\top {\bm \beta}^{\bm s}  + \Gamma \mu^{\bm s} \\
\displaystyle \sum_{i \in \sets I} {\bm x}^i \; \sum_{\kappa \in \sets K} \; {\bm s}_\kappa \;  ( \overline{\bm v}_i^{{\bm s} \kappa} - \overline{\bm w}_i^{{\bm s} \kappa} ) + {\bm B}^\top {\bm \beta}^{\bm s} \; = \; {\bm x}^{\bm s} \\
{\bm \alpha}^{\bm s} + \mu^{\bm s} \1 \leq {\bm 0} 
%
\end{array} \quad \quad \right\} \quad \forall {\bm s}\in \widetilde{\sets S}^K \\
& 
\quad \1^\top {\bm v}^\kappa = 1, \; \1^\top {\bm w}^\kappa = 1 ,\; 
\displaystyle  1-{\bm w}_i^\kappa \; \geq \; \sum_{i':i'\geq i} {\bm v}_{i'}^\kappa \quad \quad \forall  i \in \sets I  ,\; \forall \kappa \in \sets K \\ 
%
& \quad \!\!\left.  \begin{array}{l}
\overline {\bm v}^{{\bm s}\kappa} \; \leq \; M {\bm v}^\kappa, \;\; \overline {\bm v}^{{\bm s}\kappa} \; \leq \; {\bm \alpha}^{\bm s}_\kappa \1, \;\; \overline {\bm v}^{{\bm s}\kappa} \; \geq \; {\bm \alpha}^{\bm s}_\kappa \1 - M (\1 - {\bm v}^\kappa) \\
\overline {\bm w}^{{\bm s}\kappa} \; \leq \; M {\bm w}^\kappa, \;\; \overline {\bm w}^{{\bm s}\kappa} \; \leq \; {\bm \alpha}^{\bm s}_\kappa \1, \;\; \overline {\bm w}^{{\bm s}\kappa} \; \geq \; {\bm \alpha}^{\bm s}_\kappa \1 - M (\1 - {\bm w}^\kappa)
\end{array} \quad   \right\}  \quad
\begin{array}{c}
\forall {\bm s} \in \widetilde{\sets S}^K \\
\kappa \in \sets K,
\end{array}
\end{array}
\label{eq:offline_mmu_MILP}
\end{equation}
where $M$ is a ``big-$M$'' constant. Given an optimal solution $(\tau,\{{\bm v}^\kappa,{\bm w}^\kappa\}_{\kappa \in \sets K},\{ {\bm \alpha}^{\bm s}, {\bm \beta}^{\bm s},{\bm x}^{\bm s} , \mu^{\bm s} \}_{{\bm s}\in \widetilde{\sets S}^K},  \{ \overline{\bm v}_i^{{\bm s} \kappa}, \overline{\bm w}_i^{{\bm s} \kappa} \}_{ {\bm s}\in \widetilde{\sets S}^K, \kappa \in \sets K}    )$ to Problem~\eqref{eq:offline_mmu_MILP}, an optimal set of queries for~\eqref{eq:offline_mmu} is
$$
{\bm \iota}_1^\kappa = \sum_{i \in \sets I} i \cdot \I{{\bm v}^\kappa_i = 1} \quad \text{ and } \quad {\bm \iota}_2^\kappa = \sum_{i \in \sets I} i \cdot \I{{\bm w}^\kappa_i = 1}, \quad \kappa \in \sets K.
$$
\label{thm:offline_mmu_MILP}
\end{theorem}}
\newpv{The \newpv{first and third} sets of constraints in Problem~\eqref{eq:offline_mmu_MILP} exactly match those in~\eqref{eq:offline_mmu_4_finite_program}. The second set of constraints is almost identical between the two problems with the exception that the choice of queries is now encoded using the ${\bm v}$ and ${\bm w}$ variables.} The \newpv{fourth} set of constraints guarantees that ${\bm \iota}^\kappa_1 < {\bm \iota}^\kappa_2$, i.e., $( {\bm \iota}^\kappa_1, {\bm \iota}^\kappa_2) \in \sets C$. The \newpv{fifth} and \newpv{sixth} sets of constraints linearize the products ${\bm \alpha}_\kappa^{\bm s}{\bm v}^\kappa$ and ${\bm \alpha}_\kappa^{\bm s}{\bm w}^\kappa$. Reformulation~\eqref{eq:offline_mmu_MILP} enables us to solve~\eqref{eq:offline_mmu} as an MBLP using off-the-shelf solvers. For fixed $K$, its size is polynomial in the size of the input. Yet, it is exponential in~$K$. In the next section, we propose a decomposition strategy to speed-up solution as $K$ grows to typical values in the order of 8 to 10. 
Additional strategies are discussed in Section~\ref{sec:speedup_strategies}.


\subsubsection{Column-and-Constraint Generation}
\label{sec:offline_mmu_ccg}

The number of scenarios in Problem~\eqref{eq:offline_mmu_MILP} is exponential in~$K$. At the same time, we expect that only a moderate number of scenarios~${\bm s} \in \widetilde{\sets S}^K$ will be active in the epigraph constraint. \newpv{Column-and-constraint generation (CCG) techniques are a popular approach for addressing problems that possess an exponential number of decision variables and constraints while presenting a decomposable structure, see e.g., \cite{Fischetti1997,Lobel1998,Carvalho1999,Mamer2000,Feillet2010,SADYKOV2011357,Muter2013}, and \cite{Muter2018}. In particular, CCG algorithms have been proposed for solving two-stage robust optimization problems by~\cite{zeng2013solving,Shtern2018}, and~\cite{vayanos_ROInfoDiscovery}. However, these do not apply to problems with a mix of real and discrete uncertainty and/or to problems with two-and-a-half decision stages.} We thus propose a \newpv{new CCG} algorithm to speed-up computation that applies when~$\sets R$ has finite cardinality (in which case~\eqref{eq:offline_mmu} is $\mathcal{NP}$-hard, see Theorem~\ref{thm:offline_mmu_complexity}). We describe our CCG procedure using Problem~\eqref{eq:offline_mmu_4_finite_program}. Naturally, all problems solved would need to be converted to MBLPs using techniques similar to those employed in Theorem~\ref{thm:offline_mmu_MILP}. We omit these to streamline presentation.

The key idea behind our algorithm is to decompose the problem into a relaxed \newpv{main} problem and a series of subproblems indexed by ${\bm s} \in \widetilde{\sets S}^K$. The \newpv{main} problem initially only involves a subset of the constraints (those indexed by ${\bm s} \in \sets S' \subseteq \widetilde{\sets S}^K$) and a \emph{single auxiliary MBLP} is used to iteratively identify indices~${\bm s} \in \widetilde{\sets S}^K$ for which the solution to the \newpv{main} problem becomes infeasible when plugged into subproblem~${\bm s}$. Constraints associated with infeasible subproblems are added to the \newpv{main} problem and the procedure continues until convergence. We now detail this approach.

We define the relaxed \newpv{main} problem parameterized by the index set ${\sets S}'$ as
\newpv{\begin{equation}
\tag{$\mathcal{CCG}_{\text{u}}^{\text{main}}({\sets S}')$}
\begin{array}{cl}
\maximize & \quad \tau \\
\subjectto  & \quad \tau \in \reals, \; {\bm \iota} \in \sets C^K \\
& \quad \!\!\left. \begin{array}{l}
{\bm \alpha}^{\bm s} \in \reals^K_+ , \; {\bm \beta}^{\bm s} \in \reals_+^M ,\; \mu^{\bm s} \in \reals_- ,\; {\bm x}^{\bm s} \in \sets R \\
\tau  \;\; \leq \;\;  {\bm b}^\top {\bm \beta}^{\bm s}  + \Gamma \mu^{\bm s} \\
\displaystyle \sum_{\kappa \in \sets K} {\bm s}_\kappa \; ({\bm x}^{{\bm \iota}^\kappa_1} - {\bm x}^{{\bm \iota}^\kappa_2}) \; {\bm \alpha}_\kappa^{\bm s} + {\bm B}^\top {\bm \beta}^{\bm s}  \; = \; {\bm x}^{\bm s}\\
{\bm \alpha}^{\bm s} + \mu^{\bm s} \1 \; \leq \; {\bm 0}
\end{array} \quad \quad \right\} \quad \forall {\bm s}\in {\sets S}'.
\end{array}
\label{eq:offline_mmu_ccg_rmp}
\end{equation}}
This problem only involves a subset of the decision variables and constraints of Problem~\eqref{eq:offline_mmu_4_finite_program}. Given variables $(\tau,{\bm \iota})$ feasible in the \newpv{main} problem, we define the ${\bm s}$th subproblem, ${\bm s} \in \widetilde{\sets S}^K$, through
\newpv{\begin{equation}
\tag{$\mathcal{CCG}_{\text{u}}^{\text{sub},{\bm s}}(\tau, {\bm \iota})$}
\begin{array}{rl}
\maximize & \quad 0 \\
\subjectto & \quad {\bm \alpha} \in \reals^K_+, \; {\bm \beta} \in \reals_+^M , \; \mu \in \reals_-, \; {\bm x} \in \sets R \\
& \quad \tau \; \leq \;  {\bm b}^\top {\bm \beta} + \Gamma \mu \\
& \quad \displaystyle \sum_{\kappa \in \sets K} {\bm s}_\kappa \; ({\bm x}^{{\bm \iota}^\kappa_1} - {\bm x}^{{\bm \iota}^\kappa_2}) \; {\bm \alpha}_\kappa + {\bm B}^\top {\bm \beta}  \; = \; {\bm x} \\
& \quad {\bm \alpha} + \mu \1 \; \leq \; {\bm 0}.
\end{array}   
\label{eq:offline_mmu_ccg_sp}
\end{equation}}
%
To identify indices of subproblems~\eqref{eq:offline_mmu_ccg_sp} that, given a solution $(\tau,{\bm \iota})$ to the relaxed \newpv{main} problem, are infeasible, we solve a \emph{single} feasibility MBLP defined through
\newpv{\begin{equation}
\tag{$\mathcal{CCG}_{\rm{u}}^{\rm{feas}}({\bm \iota})$}
\begin{array}{cl}
    \minimize & \quad \theta \\
    \subjectto & \quad \theta \in \reals,\; {\bm s} \in \widetilde{\sets S}^K , \; {\bm u}^{\bm x} \in \sets U^0 ,\; {\bm \epsilon}^{\bm x} \in \sets E_\Gamma \quad \forall {\bm x}\in \sets R \\
    & \quad \theta \; \geq \; ({\bm u}^{\bm x})^\top {\bm x} \quad \forall {\bm x} \in \sets R \\
    & \!\!\! \left. \begin{array}{l}
    \quad ({\bm u}^{\bm x})^\top ( {\bm x}^{{\bm \iota}^\kappa_1} - {\bm x}^{{\bm \iota}^\kappa_2})  + {\bm \epsilon}_\kappa^{\newpv{\bm x}} \; \geq \; M ( {\bm s}_\kappa -1) \\
    \quad ({\bm u}^{\bm x})^\top ( {\bm x}^{{\bm \iota}^\kappa_1} - {\bm x}^{{\bm \iota}^\kappa_2})  - {\bm \epsilon}_\kappa^{\newpv{\bm x}} \; \leq \; M( {\bm s}_\kappa + 1) 
    \end{array} \quad \right\} \quad \forall  \kappa \in \sets K, \; {\bm x}\in \sets R .
\end{array}
\label{eq:offline_mmu_ccg_feas}
\end{equation}}
The following proposition enables us to bound the optimality gap associated with a given feasible solution to the relaxed master problem.
\begin{proposition}
Let ${\bm \iota}$ be feasible in the relaxed master problem~\eqref{eq:offline_mmu_ccg_rmp}. Then, ${\bm \iota}$ is feasible in Problem~\eqref{eq:offline_mmu_4} and the objective value of ${\bm \iota}$ in Problem~\eqref{eq:offline_mmu_4} is given by the optimal objective value of Problem~\eqref{eq:offline_mmu_ccg_feas}. If $\sets R$ has fixed finite cardinality, then Problem~\eqref{eq:offline_mmu_ccg_rmp} can be solved as an MBLP using off-the-shelf solvers.
\label{prop:offline_mmu_ccg_algo_correct_1}
\end{proposition}
Proposition~\ref{prop:offline_mmu_ccg_algo_correct_1} implies that, for any ${\bm \iota}$ feasible in the relaxed master problem~\eqref{eq:offline_mmu_ccg_rmp}, the optimal value of~\eqref{eq:offline_mmu_ccg_feas} yields a lower bound to the optimal value of Problem~\eqref{eq:offline_mmu_4}. At the same time, it is evident that for any index set ${\sets S}' \subseteq \widetilde{\sets S}^K$, the optimal value of Problem~\eqref{eq:offline_mmu_ccg_rmp} yields an upper bound to the optimal objective value of Problem~\eqref{eq:offline_mmu_4}. The lemma below is key to identify indices of subproblems~${\bm s} \in \widetilde{ \sets S}^K$ that are infeasible.
\begin{lemma}
The relaxed master problem~\eqref{eq:offline_mmu_ccg_rmp} is always feasible. If~\eqref{eq:offline_mmu_ccg_rmp} is solvable, let \newpv{$(\tau,{\bm \iota}, \{{{\bm \alpha}}^{\bm s}, {\bm \beta}^{\bm s}, \mu^{\bm s},\; {\bm x}^{\bm s} \}_{{\bm s} \in {\sets S}'} )$} be an optimal solution. Else, if~\eqref{eq:offline_mmu_ccg_rmp} is unbounded, set $\tau = \infty$ and let ${\bm \iota} \in \sets C^K$ be such that~\eqref{eq:offline_mmu_ccg_rmp} is unbounded when ${\bm \iota}$ is fixed to that value. Moreover, let \newpv{$(\theta,\{ {\bm u}^{{\bm x}}, {\bm \epsilon}^{\bm x} \}_{{\bm x}\in \sets R}, {\bm s} )$} be optimal in Problem~\eqref{eq:offline_mmu_ccg_feas}. Then, the following hold:
\begin{enumerate}[label=(\roman*)]
    \item $\theta \leq \tau$;
    \item If $\theta = \tau$, then Problem~\eqref{eq:offline_mmu_ccg_sp} is feasible for all ${\bm s} \in \widetilde{\sets S}^K$;
    \item If $\theta < \tau$, then scenario~${\bm s}$ corresponds to an infeasible subproblem, i.e., Problem~\eqref{eq:offline_mmu_ccg_sp} is infeasible. 
\end{enumerate}
\label{lem:offline_mmu_ccg_algo_correct_2}
\end{lemma}
\begin{algorithm}[t!]
\SetAlgoLined
\textbf{Inputs:} Optimality tolerance $\delta$, comparison set $\sets C$, and recommendation set $\sets R$\; Initial uncertainty set $\sets U^0$ and number of queries $K$\;
\textbf{Output:} Query ${\bm \iota}^\star$ from $\sets C^K$, near optimal in Problem~\eqref{eq:offline_mmu_4} with associated objective $\theta$\;
\textbf{Initialization:} ${\bm \iota}^\star \leftarrow \emptyset$; Upper and lower bounds: ${\rm{UB}} \leftarrow +\infty$ and ${\rm{LB}} \leftarrow -\infty$\;
 Initialize index set: ${\sets S}' \leftarrow \emptyset$\;
 \While{${\rm{UB}}-{\rm{LB}} > \delta$}{
  Solve the master problem~\eqref{eq:offline_mmu_ccg_rmp}. If it is solvable, let $(\tau,{\bm \iota}, \{{{\bm \alpha}}^{\bm s}, {\bm \beta}^{\bm s}, {\bm x}^{\bm s} \}_{{\bm s} \in {\sets S}'} )$ be an optimal solution. If it is unbounded, set $\tau = \infty$ and let ${\bm \iota} \in \sets C^K$ be such that~\eqref{eq:offline_mmu_ccg_rmp} is unbounded when ${\bm \iota}$ is fixed to that value\;
  Set ${\rm{UB}} \leftarrow \tau$\;
  Solve the feasibility subproblem~\eqref{eq:offline_mmu_ccg_feas}. Let $(\theta,\{ {\bm u}^{{\bm x}} \}_{{\bm x}\in \sets R}, {\bm s} )$ denote an optimal solution\;
  Set ${\rm{LB}} \leftarrow \theta$\;
    \If{$\theta < \tau$}{
        ${\sets S}' \leftarrow {\sets S}'  \cup \{ {\bm s}\}$\;
    }
 }
 Set ${\bm \iota}^\star \leftarrow {\bm \iota}$\;
\textbf{Result:} Collection of queries ${\bm \iota}^\star$ near-optimal in~\eqref{eq:offline_mmu_4} with objective value $\theta$.
 \caption{Column-and-Constraint Generation procedure for solving Problem~\eqref{eq:offline_mmu_4}.}
 \label{alg:offline_mmu_ccg}
\end{algorithm}
%
%
Propositions~\ref{prop:offline_mmu_ccg_algo_correct_1} and Lemma~\ref{lem:offline_mmu_ccg_algo_correct_2} culminate in Algorithm~\ref{alg:offline_mmu_ccg} whose convergence is guaranteed by the following theorem.
\begin{theorem}
Algorithm~\ref{alg:offline_mmu_ccg} terminates in a final number of steps with a feasible solution to Problem~\eqref{eq:offline_mmu_4}. The objective value attained by this solution is within $\delta$ of the optimal objective value of the problem.
\label{thm:offline_mmu_ccg_algo_converges}
\end{theorem}


\subsubsection{Speed-up Strategies} We propose two additional speed-up strategies.
\label{sec:speedup_strategies}



\textbf{\textsf{Symmetry Breaking. }}
\label{sec:symmetry}
In Problem~\eqref{eq:offline_mmu}, every permutation of the queries ${\bm \iota}$ will yield another solution with the same objective function value. Correspondingly, Problem~\eqref{eq:offline_mmu_MILP} also presents several symmetric solutions: we can permute the indices $\kappa$ in the pairs $\{ ({\bm v}^\kappa,{\bm w}^\kappa) \}_{\kappa \in \sets K}$ to build solutions with the same objective. To speed-up solution time, we propose to augment our MBLP formulation with symmetry breaking constraints which eliminate symmetric solutions from the search space while preserving at least one solution from each equivalence class. 

To build the constraints, note that each comparison is uniquely defined by the vector
\begin{equation}
    {\bm y}^\kappa = {\bm v}^\kappa + {\bm w}^\kappa
    \label{eq:symmetry_breaking_1}
\end{equation}
which has exactly two nonzero elements -- corresponding to the items used in the comparison. To break this symmetry, we require that the binary vectors $\{ {\bm y}^\kappa \}_{\kappa \in \sets K}$ be \emph{lexicographically ordered}: 
if ${\bm y}^\kappa_i = {\bm y}^{\kappa+1}_i$ for all $i<j$, and ${\bm y}^\kappa_j \neq {\bm y}^{\kappa+1}_j$, then ${\bm y}^\kappa_j=0$ and ${\bm y}^{\kappa+1}_j=1$. To enforce lexicographic ordering, we introduce, in addition to~${\bm y}^\kappa$, the binary variables ${\bm z}^{\kappa \kappa'} \in \{0,1\}^I$, which satisfy ${\bm z}^{\kappa \kappa '}_i=1$ if and only if ${\bm y}_{i}^{\kappa} \neq {\bm y}_{i}^{\kappa'}$, $i\in \sets I$. These variables can be uniquely defined using the following linear constraints
%
\begin{equation}
{\bm z}^{\kappa \kappa '} \leq  {\bm y}^\kappa  + {\bm y}^{\kappa '}, \;\;\;
{\bm z}^{\kappa \kappa '} \leq 2 - {\bm y}^\kappa  - {\bm y}^{\kappa '},\;\;\;
{\bm z}^{\kappa \kappa '} \geq {\bm y}^\kappa  - {\bm y}^{\kappa '}, \;\;\;
{\bm z}^{\kappa \kappa '} \geq {\bm y}^{\kappa '} - {\bm y}^\kappa  
\quad  \forall \kappa ,\; \kappa ' \in \sets K \; : \; \kappa  < \kappa '.
\label{eq:symmetry_breaking_2}
\end{equation}
Lexicographic ordering can then be imposed by adding to all MBLP problems the constraints
\begin{equation}
{\bm y}_{i}^{\kappa'} \; \geq \; {\bm y}_{i}^{\kappa} - \sum\limits_{i' \in \sets I : i'< i} {\bm z}^{\kappa \kappa'}_{i'} \quad \forall i\in \sets I,\; \kappa, \kappa' \in \sets K:\kappa<\kappa'.
\label{eq:symmetry_breaking_3}
\end{equation}

Next, note that in the preference elicitation problems, there is no benefit in asking the same query twice: the worst-case utility can only improve by making different queries. Correspondingly, we can add to all MBLP problems the following constraints 
\begin{equation}
\sum\limits_{i\in \sets I} {\bm z}_i^{\kappa,\kappa'} \; \geq \; 1 \quad \forall \kappa, \kappa' \in \sets K \; : \; \kappa < \kappa', 
\label{eq:symmetry_breaking_4}
\end{equation}
which eliminate solutions that involve asking the same queries twice. 

As we will see in Section~\ref{sec:performance}, the constraints above translate to significant speed-ups.



\textbf{\textsf{Warm-Starts.}} Solutions obtained by solving problems with smaller numbers of queries can be used to warm start problems with more queries. In this section, we describe how to build a feasible warm start for Problem~\eqref{eq:offline_mmu_MILP}  with $K+1$ queries from a feasible solution to Problem~\eqref{eq:offline_mmu_MILP} with $K$ queries. Recall that a warm start need only fix solutions for the binary decision variables.

Let $(  \tau, {\bm v} , {\bm w} , {\bm \alpha}, {\bm \beta}, {\bm x}, \overline{\bm v}, \overline{\bm w} )$ be feasible in Problem~\eqref{eq:offline_mmu_MILP} with~$K$ queries. To generate a feasible solution to Problem~\eqref{eq:offline_mmu_MILP} with~$K+1$ queries, we proceed as follows. First, we generate a new query ${\bm \iota}^{K+1}$ randomly from the set
$$
\sets C \; \backslash \; \left\{ \; {\bm \iota} \in \sets C \; : \; \exists \kappa \in \sets K \; : \; {\bm \iota}_1^\kappa = \sum_{i \in \sets I} i \cdot \I{{\bm v}^\kappa_i = 1}  \text{ and }  {\bm \iota}_2^\kappa = \sum_{i \in \sets I} i \cdot \I{{\bm w}^\kappa_i = 1} \; \right\}
$$
and define $\tilde{\bm v}$ and $\tilde{\bm w}$ through $\tilde{\bm v}^\kappa := {\bm v}^\kappa$, $\tilde{\bm w}^\kappa := {\bm w}^\kappa$, $\kappa \in \sets K$, 
$
\tilde{\bm v}_i^{K+1} :=
\I {{\bm \iota}^{K+1}_1 = i}
$, and
$\tilde {\bm w}_i^{K+1} :=
\I{ {\bm \iota}^{K+1}_2 = i }
$
for each $i\in \sets I$. Subsequently, for each ${\bm s} \in \widetilde{\sets S}^{K+1}$, we define 
$
\tilde {\bm x}^{\bm s} \; := \; {\bm x}^{({\bm s}_1,\ldots,{\bm s}_{K})}.
$
Then, $( \tilde{\bm v},\tilde{\bm w}, \tilde{\bm x} )$ constitutes a warm start for Problem~\eqref{eq:offline_mmu_MILP} with $K+1$ queries. The solution hereby constructed asks one additional (random) query and subsequently ignores the answer to the query in the choice of recommendation. It can be shown that it attains an objective value in Problem~\eqref{eq:offline_mmu_MILP} with $K+1$ queries at least as high as that attained by $( \tau, {\bm v} , {\bm w} , {\bm \alpha}, {\bm \beta}, {\bm x}, \overline{\bm v}, \overline{\bm w} )$ in Problem~\eqref{eq:offline_mmu_MILP} with $K$ queries. Note that if warm starts are combined with symmetry breaking constraints, then the solution constructed needs to be permuted so as to satisfy lexicographical ordering to ensure it satisfies the lexicographic constraints. This procedure is detailed in Section~\ref{sec:EC_speedup_strategies}. \newpv{In our numerical experiments, we will use a variant of this warm-start procedure which fixes the queries optimal in the problem with $K$ queries in the problem with $K+1$ queries to speed-up its solution. This approach is detailed in Section~\ref{sec:EC_greedy}.}


\subsection{Online Active Preference Elicitation with the Max-Min Utility Decision Criterion}
\label{sec:online_maxminutility}

\newpv{In Section~\ref{sec:offline_maxminutility}, we assumed that all queries are chosen at once \newpv{before the answers to any of the queries are observed}. However, in many settings of practical interest, queries are made one at a time and the answer to each query is revealed before the next query is selected. In this case, the recommender system has the opportunity to \emph{adjust} their choice of queries, taking into account the information acquired as a byproduct of previous queries and their answers.

After making simplifications akin to those in Lemma~\ref{lem:offline_mmu_simplification} but adapted to this multi-stage setting, we can write the online active preference elicitation problem as
\begin{equation}
\tag{$\onlinerisk{K}$}
\max\limits_{{\bm \iota}^1 \in \mathcal C} \quad \min\limits_{{\bm s}_1 \in \mathcal S} \quad
\max\limits_{{\bm \iota}^2 \in \mathcal C} \quad \min\limits_{{\bm s}_2 \in \mathcal S} \quad
\cdots \quad
\max\limits_{{\bm \iota}^K \in \mathcal C} \quad \min\limits_{{\bm s}_K \in \mathcal S} \quad
\max\limits_{{\bm x} \in \sets R} \quad \min\limits_{ {\bm u} \in  \widetilde{\sets U}_\Gamma({\bm \iota} ,{\bm s})} \quad  {\bm u}^\top {\bm x}.
\label{eq:online_mmu}
\end{equation}
Indeed, for reasons similar to those in the proof of Lemma~\ref{lem:offline_mmu_simplification}, we can allow responses ${\bm s}^\kappa$ to be inconsistent with any ${\bm u} \in \widetilde{\sets U}_\Gamma({\bm \iota} ,{\bm s})$ since this will not affect the optimal value nor the sets of optimal solutions of the query selection problems. Moreover, we can ignore indifferent responses without loss. In spite of these simplifications, Problem~\eqref{eq:online_mmu} is still significantly harder to solve than its offline counterpart~\eqref{eq:offline_mmu} (or even~\eqref{eq:offline_mmu_4}).

We propose to leverage the tools we devised in Section~\ref{sec:offline_maxminutility} to obtain a conservative solution to the online problem. Specifically, we choose, at each step $k \in \sets K$, a query ${\bm \iota}^k$ as if it were the last query to be asked, i.e., in a one-step look-ahead fashion by solving problem~(\offlinerisk{1}) after updating its uncertainty set with the answers $\{{\bm s}_\kappa\}_{\kappa=1}^{k-1}$ obtained to the queries $\{{\bm \iota}^\kappa\}_{\kappa=1}^{k-1}$ actually asked. Concretely, we compute the $k$th query by solving
\begin{equation*}
\maximize_{{\bm \iota}^k \in \sets C} \quad \min_{{\bm s}_k \in \widetilde{\sets S}} \quad \max_{{\bm x} \in \sets R} \quad \min_{ {\bm u} \in \widetilde{\sets U}_{\Gamma(k)}({\bm \iota}^{[k]} ,{\bm s}_{[k]}) } \quad  {\bm u}^\top {\bm x},
\end{equation*}
where, ${\bm \iota}^{[k]}:=\{ {\bm \iota}^{\kappa} \}_{\kappa=1}^k$, ${\bm s}_{[k]}:=\{ {\bm s}_{\kappa} \}_{\kappa=1}^k$, and with a slight abuse of notation
\begin{equation*}
\widetilde{\sets U}_{\Gamma(k)}({\bm \iota}^{[k]},{\bm s}_{[k]}) := 
\left\{
\begin{array}{ccll}  
{\bm u} \in \sets U^0 & : &  \exists {\bm \epsilon} \in \sets E_{\Gamma(k)} \text{ such that} \\
&& {\bm u}^\top ({\bm x}^{{\bm \iota}^\kappa_1} - {\bm x}^{{\bm \iota}^\kappa_2}) \; \geq \; -{\bm \epsilon}_\kappa &  \quad  \forall \kappa \in \{1,\ldots,k\} \; : \; {\bm s}_\kappa=1 \\
&& {\bm u}^\top ( {\bm x}^{{\bm \iota}^\kappa_1} - {\bm x}^{{\bm \iota}^\kappa_2} ) \; \leq \; {\bm \epsilon}_\kappa  &  \quad  \forall \kappa \in \{1,\ldots,k\} \; : \; \tilde{\bm s}_\kappa=-1
\end{array}
\right\}.
\end{equation*}
Note in particular that in the above definition, ${\bm \iota}^{[k-1]}$ and ${\bm s}_{[k-1]}$ are data. We propose to choose $\Gamma(\kappa)$ in a fashion that exactly parallels the choice of~$\Gamma$ in the multi-stage problem.
\begin{example}[Choice of $\Gamma(k)$]
If $\Gamma$ is chosen as in Example~\ref{ex:choice_of_gamma}, then a natural choice is $\Gamma(k) = 2\sigma \sqrt{k} \; \rm{erf}^{-1}(2p-1)$ which ensures that ${\bm \epsilon} \in \sets E_{\Gamma(k)} \subseteq \reals^k$ with probability~$p$.
\label{ex:choice_of_gamma_kappa}
\end{example}

}

\subsection{\newpv{Performance}}
\label{sec:performance}

\newpv{We perform a wide range of experiments on randomly generated datasets of various sizes ($I$, $J$) and varying degree of inconsistencies ($\Gamma$) to evaluate the performance of our approach (\texttt{MMU}). To speed-up computation and to allow the solution of a large number of problems, unless stated otherwise, we employ the CCG algorithm from Section~\ref{sec:offline_mmu_ccg} which we augment with the symmetry breaking constraints from Section~\ref{sec:symmetry}, and a greedy-heuristic that builds upon the warm-start approach, see Section~\ref{sec:EC_speedup_strategies}. To make it easier to evaluate the performance of any given method, we normalize utilities so that they lie in the $[0,1]$ range, where 0 corresponds to the worst-case utility if a recommendation is made before any questions are asked and 1 is the worst-case utility if the recommendation is made with full-knowledge of the utility. We provide here a summary of our results. Details of our experiments can be found in Section~\ref{sec:EC_performance}.

In the offline setting, we compare our approach to random elicitation (\texttt{RAND}) where the performance of a random query is evaluated by plugging it into Problem~\eqref{eq:offline_mmr} (and using the methods in this paper to evaluate the objective of the resulting problem). The results are summarized in Table~\ref{tab:offline_mmu}. From the table, it is apparent that computing an optimal set of queries using \texttt{MMU} is consistently preferable to selecting queries at random, resulting in an improvement in worst-case utility of over 63\% on average across datasets and values of~$K$. Moreover, we observe that while computing the worst-case utility of a set of queries (which is only possible thanks to our approach) is very fast (under 3 seconds on average), evaluating the quality of all sets of queries of cardinality~$K$ is computationally prohibitive (e.g., 2.3E+09 seconds for $K=3$ on average). This motivates the use of our \texttt{MMU} approach that can compute optimal queries in practical times for the application (under 3 hours for realistic problem sizes). From the table it can also be noted that much of the useful information is gathered in the first 5 queries (increase in worst-case utility from 0.62 to 0.72). For completeness, we also conduct a computational study to evaluate the performance of the CCG algorithm and symmetry breaking constraints. The results are summarized in Table~\ref{tab:speedup_summary}. The table shows that relative to solving MBLP~\eqref{eq:offline_mmu_MILP} directly, augmenting it with symmetry breaking constraints and using the CCG algorithm can yield significant speed-ups (up to 135 times on average). The table also shows that the heuristic used to speed-up computations is near-optimal, with an average optimality gap of under~$2.1\%$.
}


\begin{table}[ht]
\caption{\newpv{Performance of our offline elicitation procedure relative to random elicitation. All numbers correspond to averages across randomly generated datasets of various sizes. The numbers in the runtime columns correspond to the time it takes to compute an optimal query and to evaluate the performance of a random query (average over 50 sets of queries) in the case of \texttt{MMU} and \texttt{RAND}, respectively. The numbers in parenthesis in the last column correspond to the expected time it would take to evaluate the performance of all possible sets of~$K$ random queries.}}
\centering
\renewcommand{\arraystretch}{1.1}
\begin{small}
\begin{tabular}{|C{1cm}||C{2.5cm}|C{2.5cm}||C{2.5cm}|C{2.5cm}|}
  \hline
\multirow{2}{*}{K} & \multicolumn{2}{c||}{Worst-Case Utility (Normalized)} & \multicolumn{2}{c|}{Runtime (s)} \\
\cline{2-5}
& \texttt{MMU} & \texttt{RAND} & \texttt{MMU} & \texttt{RAND} \\
  \hline
 1 & \textbf{0.62} & 0.25 & \textbf{35} & 1 (1.3E+03) \\ 
   2 & \textbf{0.62} & 0.32 & \textbf{1668} & 1 (2.1E+06) \\ 
   3 & \textbf{0.66} & 0.37 & \textbf{1854} & 1 (2.3E+09) \\ 
   4 & \textbf{0.69} & 0.41 & \textbf{2346} & 1 (1.9E+12) \\ 
   5 & \textbf{0.72} & 0.44 & \textbf{2936} & 1 (1.3E+15) \\ 
   6 & \textbf{0.75} & 0.46 & \textbf{3370} & 2 (8.1E+17) \\ 
   7 & \textbf{0.75} & 0.49 & \textbf{4435} & 2 (4.8E+20) \\ 
   8 & \textbf{0.77} & 0.52 & \textbf{5759} & 2 (2.5E+23) \\ 
   9 & \textbf{0.78} & 0.54 & \textbf{6637} & 2 (1.1E+26) \\ 
  10 & \textbf{0.78} & 0.56 & \textbf{8121} & 3 (6.0E+28) \\ 
  \hline 
  Avg. & \textbf{0.71} & 0.44 & \textbf{3716} & 1.66 (6.0e+27) \\
   \hline
\end{tabular}
\end{small}
\label{tab:offline_mmu}
\end{table}


\begin{table}[ht]
\caption{\newpv{Performance of our speed-up methods and of the employed heuristic relative to solving the MBLP~\eqref{eq:offline_mmu_MILP} directly. All numbers correspond to averages over randomly generated datasets.}}
\centering
\renewcommand{\arraystretch}{1.1}
\begin{small}
\begin{tabular}{|l|C{2cm}|}
\hline
Avg.\ speed-up of symmetry breaking & $\times65$ \\
Avg.\ speed-up of CCG+symmetry breaking & $\times135$ \\
Avg.\ optimality gap of heuristic & 2.03\% \\
\hline
\end{tabular}
\end{small}
\label{tab:speedup_summary}
\end{table}

\newpv{
In the online setting, we compare to random elicitation, to the polyhedral method of~\cite{Toubia_2004} (\texttt{POLY}), to the probabilistic polyhedral method of~\cite{Toubia_2007} (\texttt{PROB}), to the robust approach of~\cite{OHair_LearningPreferences} (\texttt{ROB}), and to the ellipsoidal method of~\cite{Vielma2019} (\texttt{ELL}). To evaluate the performance of any given method, we conduct simulations, wherein we draw 50 agents at random and simulate their responses (with and without inconsistences). We record statistics on the guaranteed utility and rank of the recommended item which are summarized on Table~\ref{tab:online_mmu}. We omit solver times since all methods identify a query in under a second. From the table, it can be seen that our method consistently outperforms all other approaches in terms of both worst-case and average guaranteed utility (by over $\times2.4$ and $\times1.5$ on average, respectively). It also results in improvements in worst-case and true rank of the recommendation (by over $1.8$ and $1.2$ spots on average, respectively). The table also shows that with \texttt{MMU}, most of the information is gained within the first 4-5 queries and that recommendations with almost optimal worst-case utility (0.91 on average) can be made after just 10 queries.

}

\begin{table}[ht]
\caption{\newpv{Performance of our online elicitation procedure relative to state of the art methods from the literature over 50 randomly generated agents. All numbers correspond to averages across randomly generated datasets.}}
\centering
\renewcommand{\arraystretch}{1.1}
\begin{small}
\begin{tabular}{|C{0.9cm}||C{0.9cm}|C{0.9cm}|C{0.9cm}|C{0.9cm}|C{0.9cm}|C{0.9cm}||C{0.95cm}|C{0.9cm}|C{0.9cm}|C{0.9cm}|C{0.9cm}|C{0.95cm}|}
     \hline
\multirow{2}{*}{K} & \multicolumn{6}{c||}{Worst-Case Guaranteed Utility (Normalized)} & \multicolumn{6}{c|}{Worst-Case True Rank of Recommendation} \\
\cline{2-13}
 & \texttt{MMU} & \texttt{RAND} & \texttt{POLY} & \texttt{PROB} & \texttt{ROB} & \texttt{ELL} & \texttt{MMU} & \texttt{RAND} & \texttt{POLY} & \texttt{PROB} & \texttt{ROB} & \texttt{ELL} \\ 
  \hline
1 & \textbf{0.62} & 0.05 & 0.35 & 0.35 & 0.12 & 0.58 & 27.71 & 35.57 & 32.86 & 32.86 & 37.86 & \textbf{26.71} \\ 
  2 & \textbf{0.74} & 0.15 & 0.16 & 0.19 & 0.20 & 0.26 & \textbf{24.43} & 34.71 & 33.43 & 32.14 & 36.71 & 33.00 \\ 
  3 & \textbf{0.81} & 0.22 & 0.20 & 0.17 & 0.26 & 0.20 & \textbf{21.57} & 32.14 & 32.00 & 33.57 & 35.29 & 32.57 \\ 
  4 & \textbf{0.84} & 0.27 & 0.18 & 0.19 & 0.30 & 0.20 & \textbf{19.71} & 30.86 & 31.71 & 30.00 & 35.00 & 29.71 \\ 
  5 & \textbf{0.87} & 0.36 & 0.15 & 0.16 & 0.34 & 0.26 & \textbf{17.57} & 29.00 & 28.57 & 29.43 & 33.57 & 29.00 \\ 
  6 & \textbf{0.89} & 0.38 & 0.18 & 0.19 & 0.36 & 0.29 & \textbf{16.43} & 29.29 & 28.14 & 29.43 & 30.71 & 26.71 \\ 
  7 & \textbf{0.89} & 0.42 & 0.16 & 0.20 & 0.37 & 0.31 & \textbf{16.57} & 27.43 & 29.14 & 27.43 & 29.43 & 26.43 \\ 
  8 & \textbf{0.90} & 0.50 & 0.16 & 0.18 & 0.36 & 0.37 & \textbf{14.43} & 27.29 & 28.86 & 30.29 & 29.00 & 25.86 \\ 
  9 & \textbf{0.91} & 0.52 & 0.27 & 0.29 & 0.37 & 0.42 & \textbf{13.71} & 25.29 & 27.43 & 27.14 & 28.43 & 22.43 \\ 
  10 & \textbf{0.91} & 0.54 & 0.25 & 0.28 & 0.40 & 0.50 & \textbf{13.86} & 25.43 & 28.43 & 27.29 & 29.86 & 19.43 \\
  \hline
  Avg. & \textbf{0.84} & 0.34 & 0.21 & 0.22 & 0.31 & 0.34 & \textbf{18.60} & 29.70 & 30.06 & 29.96 & 32.59 & 27.19 \\
  \hline \hline
\multirow{2}{*}{K} & \multicolumn{6}{c||}{Average Guaranteed Utility (Normalized)} & \multicolumn{6}{c|}{Average True Rank of Recommendation} \\
\cline{2-13}
 & \texttt{MMU} & \texttt{RAND} & \texttt{POLY} & \texttt{PROB} & \texttt{ROB} & \texttt{ELL} & \texttt{MMU} & \texttt{RAND} & \texttt{POLY} & \texttt{PROB} & \texttt{ROB} & \texttt{ELL} \\ 
  \hline
1 & \textbf{0.62} & 0.23 & 0.36 & 0.36 & 0.17 & 0.59 & 11.59 & 13.72 & 14.11 & 14.13 & 14.43 & \textbf{11.35} \\ 
  2 & \textbf{0.75} & 0.34 & 0.38 & 0.40 & 0.29 & 0.49 & \textbf{9.19} & 12.65 & 12.83 & 11.44 & 13.83 & 11.26 \\ 
  3 & \textbf{0.82} & 0.40 & 0.40 & 0.41 & 0.36 & 0.43 & \textbf{7.17} & 11.88 & 11.25 & 10.58 & 11.87 & 10.86 \\ 
  4 & \textbf{0.86} & 0.47 & 0.43 & 0.47 & 0.41 & 0.44 & \textbf{6.31} & 11.04 & 11.46 & 8.49 & 10.90 & 8.86 \\ 
  5 & \textbf{0.89} & 0.52 & 0.47 & 0.50 & 0.45 & 0.52 & \textbf{5.65} & 10.26 & 7.96 & 8.07 & 10.05 & 7.53 \\ 
  6 & \textbf{0.91} & 0.56 & 0.53 & 0.55 & 0.50 & 0.55 & \textbf{5.23} & 9.20 & 7.81 & 7.21 & 9.36 & 6.58 \\ 
  7 & \textbf{0.92} & 0.61 & 0.54 & 0.56 & 0.53 & 0.59 & \textbf{5.18} & 8.43 & 7.17 & 6.76 & 8.68 & 6.08 \\ 
  8 & \textbf{0.93} & 0.65 & 0.58 & 0.61 & 0.54 & 0.64 & \textbf{4.77} & 7.56 & 6.99 & 6.47 & 8.37 & 5.08 \\ 
  9 & \textbf{0.93} & 0.68 & 0.62 & 0.63 & 0.56 & 0.70 & \textbf{4.39} & 7.18 & 6.28 & 6.26 & 7.85 & 4.62 \\ 
  10 & \textbf{0.94} & 0.70 & 0.64 & 0.66 & 0.59 & 0.73 & 4.32 & 7.05 & 6.71 & 5.73 & 7.82 & \textbf{4.00} \\
  \hline
  Avg. & \textbf{0.86} & 0.54 & 0.49 & 0.51 & 0.44 & 0.57 & \textbf{6.38} & 9.72 & 9.26 & 8.51 & 10.32 & 7.62 \\
  \hline
\end{tabular}
\end{small}
\label{tab:online_mmu}
\end{table}




\section{Active Elicitation with the Min-Max Regret Decision Criterion}
\label{sec:minimaxregret}

\newpv{We now study the offline and online active preference elicitation problems under the min-max regret decision criterion and evaluate their performance on synthetic data.}



\subsection{Active Offline Preference Elicitation with the Min-Max Regret Decision Criterion}
\label{sec:offline_minmaxregret}


\subsubsection{Problem Formulation, Complexity Analysis, \& MBLP Reformulation.}

\newpv{The offline active preference elicitation problem under the min-max regret decision criterion is a variant of its max-min utility counterpart where the objective of the recommender system when selecting queries and making recommendations is to minimize worst-case regret (rather than maximize worst-case utility).} Mathematically, it is expressible as the following two-stage robust optimization problem with decision-dependent information discovery
\begin{equation}
\tag{$\offlineregret{K}$}
\minimize_{{\bm \iota} \in \sets C^K} \quad \max_{{\bm s}\in \newpv{\sets S_\Gamma}({\bm \iota})} \quad \min_{{\bm x} \in \sets R} \quad \max_{ {\bm u} \in \newpv{\sets U_\Gamma}({\bm \iota} ,{\bm s})} \quad \left\{ \;\; \max_{   {\bm x}' \in \sets R } \;\; {\bm u}^\top {\bm x}' - {\bm u}^\top {\bm x} \;\; \right\}.
\label{eq:offline_mmr}
\end{equation}
Problem~\eqref{eq:offline_mmr} is generally $\mathcal{NP}$-hard when the recommendation set is discrete \newpv{as summarized in the following theorem}.

\begin{theorem}
    %
    %
    Evaluating the objective function of Problem~\eqref{eq:offline_mmr} is $\mathcal{NP}$-hard even if the recommendation set $\sets R$ consists of only two elements.
\label{thm:offline_mmr_complexity}
\end{theorem}
\newpv{We note that in contrast to the max-min utility setting, since the objective function of the min-max regret recommendation problem is not concave in~${\bm u}$, the minimax theorem that could be used to simplify the max-min utility problem when $\sets R$ is convex (see Theorem~\ref{thm:offline_mmu_complexity}) does not apply here. This result motivates us to reformulate Problem~\eqref{eq:offline_mmr} as an MBLP as shown in the following theorem}. Throughout the remainder of this section, we assume that $\sets R$ has fixed finite cardinality.

\begin{theorem}
Problem~\eqref{eq:offline_mmr} is equivalent to the following mixed-binary linear program
\newpv{\begin{equation}
\begin{array}{cl}
\minimize & \quad \tau \\
\subjectto  & \quad \tau \in \reals, \; {\bm v}^\kappa, \; {\bm w}^\kappa \in \{0,1\}^I , \; \kappa \in \sets K  \\
& \quad {\bm x}^{\bm s} \in \sets R \quad \forall {\bm s}\in \widetilde{\sets S}^K \\
& \quad {\bm \alpha}^{(\bm x', \bm s)} \in \reals_-^K, \; {\bm \beta}^{(\bm x', \bm s)} \in \reals_-^M , \; \mu^{(\bm x', \bm s)} \in \reals_+ \quad \forall {\bm s}\in \widetilde{\sets S}^K, \; {\bm x}' \in \sets R \\
& \quad \overline{\bm v}^{(\bm x', \bm s)\kappa}, \; \overline{\bm w}^{(\bm x', \bm s)\kappa} \in \reals^I_- \quad \forall {\bm s}\in \widetilde{\sets S}^K, \; {\bm x}' \in \sets R , \; \kappa \in \sets K \\
& \quad \!\!\left. \begin{array}{l}
\tau  \;\; \geq \;\; {\bm b}^\top {\bm \beta}^{(\bm x', \bm s)} + \Gamma \mu^{(\bm x', \bm s)} \\
\displaystyle \sum_{i \in \sets I} {\bm x}^i \; \sum_{\kappa \in \sets K} \; {\bm s}_\kappa \;  \left( \overline{\bm v}_i^{(\bm x', \bm s) \kappa} - \overline{\bm w}_i^{(\bm x', \bm s) \kappa} \right) + {\bm B}^\top {\bm \beta}^{(\bm x', \bm s)} \; = \; {\bm x}' - {\bm x}^{\bm s} \\
{\bm \alpha}^{(\bm x', \bm s)} + \mu^{(\bm x', \bm s)} \1 \geq {\bm 0}
\end{array} \quad \quad \right\} \quad 
\begin{array}{c} \forall {\bm s}\in \widetilde{\sets S}^K, \\ {\bm x}' \in \sets R \end{array} \\
& \quad 
 \1^\top {\bm v}^\kappa = 1, \; \1^\top {\bm w}^\kappa = 1 , \;
\displaystyle  1-{\bm w}_i^\kappa \; \geq \; \sum_{i':i'\geq i} {\bm v}_{i'}^\kappa \quad \forall  i \in \sets I, \;
\forall \kappa \in \sets K \\
& \quad \!\!\left.  \begin{array}{l}
  \overline {\bm v}^{(\bm x', \bm s)\kappa} \; \geq \; -M {\bm v}^\kappa ,\; \overline {\bm v}^{(\bm x', \bm s)\kappa} \; \geq \; {\bm \alpha}^{(\bm x', \bm s)}_\kappa \1 \\
 \overline {\bm v}^{(\bm x', \bm s)\kappa} \; \leq \; {\bm \alpha}^{(\bm x', \bm s)}_\kappa \1 + M (\1 - {\bm v}^\kappa) \\
 \overline {\bm w}^{(\bm x', \bm s)\kappa} \; \geq \; -M {\bm w}^\kappa ,\; \overline {\bm w}^{(\bm x', \bm s)\kappa} \; \geq \;  {\bm \alpha}^{(\bm x', \bm s)}_\kappa \1 \\
 \overline {\bm w}^{(\bm x', \bm s)\kappa} \; \leq \; {\bm \alpha}^{(\bm x', \bm s)}_\kappa \1 + M (\1 - {\bm w}^\kappa)
\end{array} \quad   \right\}  \quad
\begin{array}{c}
\forall {\bm s} \in \widetilde{\sets S}^K, \; {\bm x}' \in \sets R, \\
\kappa \in \sets K,
\end{array}
\end{array}
\label{eq:offline_mmr_MILP}
\end{equation}}
where $M$ is a ``big-$M$'' constant. In particular, given an optimal solution $(\tau,{\bm v},{\bm w}, {\bm \alpha}, {\bm \beta},\overline{\bm v},\overline{\bm w})$ to Problem~\eqref{eq:offline_mmr_MILP}, an optimal set of queries for Problem~\eqref{eq:offline_mmr} is given by
$$
{\bm \iota}_1^\kappa = \sum_{i \in \sets I} i \cdot \I{{\bm v}^\kappa_i = 1} \quad \text{ and } \quad {\bm \iota}_2^\kappa = \sum_{i \in \sets I} i \cdot \I{{\bm w}^\kappa_i = 1}, \quad \kappa \in \sets K.
$$
\label{thm:offline_mmr_MILP}
\end{theorem}
%


\subsubsection{Column-and-Constraint Generation}
\label{sec:offline_mmr_ccg}

For fixed $K$, Problem~\eqref{eq:offline_mmr_MILP} is polynomial in the size of the input. Yet, it presents a number of decision variables and constraints that are exponential in~$K$. We thus propose a column-and-constraint generation approach, similar to the one proposed in Section~\ref{sec:offline_mmu_ccg} for the max-min utility case, that enables us to scale to practical values of~$K$. In the interest of space, we defer the details to the electronic companion, see Section~\ref{sec:EC_CCG_MMR}.

\newpv{Naturally, the speed-up strategies from Section~\ref{sec:speedup_strategies} directly apply to the min-max regret setting.}



\subsection{Active Online Elicitation with the Min-Max Regret Decision Criterion}
\label{sec:online_minmaxregret}

\newpv{Similar to the max-min utility case, we can formulate the online version of the active preference elicitation problem with the min-max regret decision criterion as a min-max-\ldots-min-max-max problem with decision-dependent uncertainty sets. We propose to solve this problem approximately in a one-step look-ahead fashion in a way that exactly parallels the approach from Section~\ref{sec:online_maxminutility}}.

\subsection{\newpv{Performance}}
\label{sec:performance_mmr}

\newpv{We perform a wide set of experiments on randomly generated data to evaluate the performance of our min-max regret elicitation procedure (\texttt{MMR}). Our results are summarized in Table~\ref{tab:offline_online_mmr}. Since these are qualitatively very similar to those reported in the max-min utility case, see Section~\ref{sec:performance}, we only report averages across datasets and values of $K$ and refer the reader to Section~\ref{sec:EC_performance} for more information. We note that since the regret problem is harder than the utility elicitation problem, the instances considered here were smaller than those in Section~\ref{sec:performance} to be able to report results on a variety of problems. From the table, we observe that our proposed approach outperforms existing methods on all metrics on average. Most notably, our offline method reduces the worst-case regret by a factor of 2 on average. Although the reported numbers are averages, a similar behavior is exhibited by individual realizations.}

\begin{table}[ht]
\caption{\newpv{Performance of our min-max regret offline (top two rows) and online (bottom four rows) elicitation procedures relative to state of the art methods from the literature. All numbers correspond to averages across $K$ and randomly generated datasets; for the online elicitation case, we use 50 randomly generated agents.}}
\centering
\renewcommand{\arraystretch}{1.1}
\begin{small}
\begin{tabular}{|C{7cm}||C{0.9cm}|C{2cm}|C{0.9cm}|C{0.9cm}|C{0.9cm}|C{0.9cm}|}
     \hline
 & \texttt{MMR} & \texttt{RAND} & \texttt{POLY} & \texttt{PROB} & \texttt{ROB} & \texttt{ELL} \\ 
  \hline
Worst-Case Regret (Normalized) & \textbf{0.51} & 0.66 & -- & -- & -- & -- \\
Runtime (s) & \textbf{1540} & 2 (1E+20)  & -- & -- & -- & --  \\
\hline 
Worst-Case Guaranteed Regret (Normalized) & \textbf{0.38}  &  0.77 & 0.82  & 0.83  &  0.72  &   0.58  \\
Worst-Case True Rank of Recommendation & \textbf{4.37} & 7.74 & 8.04  &  8.33 & 7.69   &  5.57   \\
Average Guaranteed Regret (Normalized) & \textbf{0.31} & 0.54 & 0.59  & 0.62  &  0.53  & 0.41    \\
Average True Rank of Recommendation & \textbf{1.72} & 2.56 & 2.65  & 2.66  & 2.74   & 1.95    \\
  \hline 
\end{tabular}
\end{small}
\label{tab:offline_online_mmr}
\end{table}

\newpv{In our last set of experiments on synthetic data, we investigate whether there are benefits in employing the min-max regret solution relative the max-min utility solution -- in particular, whether it results in solutions that have lower regret than solutions to the problem that maximizes worst-case utility. The results are summarized on Table~\ref{tab:mmu_vs_mmr_summary}, where it can be seen that the worst-case regret queries (worst-case utility) have lower regret (higher utility) than their max-min utility (min-max regret) counterparts. This implies that although computing min-max regret queries comes at a higher computational cost, they may be worthwhile to compute if one is keen on having low regret.

\begin{table}[ht]
\caption{Comparison between max-min utility and min-max regret based queries on a randomly generated synthetic dataset. All approaches were given a one hour time limit. The numbers in the table correspond to averages across values of $K$.}
\renewcommand{\arraystretch}{1.1}
\centering
\begin{small}
\begin{tabular}{|C{3cm}|C{3cm}||C{3cm}|C{3cm}|} 
     \hline
\multicolumn{2}{|c||}{Max-Min Utility Solution} &	\multicolumn{2}{c|}{Min-Max Regret Solution} \\
\hline 
 Normalized Worst-Case Utility & Normalized Worst-Case Regret &	Normalized Worst-Case Utility &		Normalized Worst-Case Regret  \\ 
 \hline
\textbf{0.84} &  0.58 & 0.75 & \textbf{0.53} \\
\hline
    \end{tabular}
\end{small}
    \label{tab:mmu_vs_mmr_summary}
\end{table}

}


\section{\newpv{Designing Allocation Policies that Align with Human Value Judgements}}
\label{sec:numericals}

\newpv{In this section, we investigate a real-world application of our preference elicitation framework. In particular, we consider the problem of eliciting the preferences of policy-makers in charge of allocating scarce housing resources to people experiencing homelessness. We envision our methodology will enable homeless services organizations such as LAHSA to allocate resources in a way that aligns with the values of its policy-makers to adequately balance the interests of those they serve. We first describe the L.A.\ housing allocation system and motivate various attributes of policies that are relevant to be considered when evaluating them against one another. We then discuss how the dataset we have at our disposal can be used to estimate the performance of counterfactual housing allocation policies in terms of these metrics. We conclude by performing a numerical case study based on historical data that showcases the promise of the proposed approach.}

\subsection{\newpv{The L.A.\ Housing Allocation System}}
\label{sec:LA_housing_system}


\newpv{In LA County, LAHSA coordinates housing and services for homeless families and individuals. It employs a \emph{Coordinated Entry System (CES)} which facilitates the coordination and management of resources and services. The CES enables organizations within the same community to pool both their housing resources and clients. The County is divided into eight geographic areas designated as Service Planning Areas (SPAs). Each SPA is expected to have a balance of homeless services used to serve clients in its area and consists of three systems that operate mostly independently of one another: the Youth system that serves those under 24 years old, the Single Adult system, and the Family system. Each SPA System has a ``matcher'' in charge of working with LAHSA and with case workers that know the clients to match resources as they become available to waitlisted clients. For simplicity, and to align with the data we have at our disposal, we focus here on the Single Adult system and note that the other systems operate in a similar fashion. 

When an individual interacts with the housing allocation system for the first time, e.g., through street outreach or by staying at a shelter, they typically answer a survey, the Vulnerability Index -- Service Prioritization Decision Assistance Tool (VI-SPDAT), which will help determine risk and prioritization when providing them with assistance \citep{VISPDAT}. The VI-SPDAT consists of~34 predominantly yes-or-no items intended to measure an individual’s level of vulnerability across four domains: their history of housing and homelessness, individual risk factors, socialization and daily functions, and wellness. The answers to the VI-SPDAT survey are also used to compute a \emph{risk score} ranging from~0 to~17 with higher scores indicating greater vulnerability. This information is entered into the Homeless Management Information System (HMIS) together with client demographic information (e.g., age, race/ethnicity). At this point, they are waitlisted for housing in their SPA and matchers can see them in their system, in what is called the \emph{community queue}.}

\newpv{Housing resources can broadly be categorized in two groups: \emph{permanent supportive housing (PSH)}, which combines low-barrier affordable housing, health care, and other supportive services; and \emph{rapid rehousing (RRH)}, which provides short-term rental assistance and services. Individuals who do not receive a housing resource will still receive services depending on their needs. 

To help decide how to match people to resources, a housing allocation policy is in place which, roughly speaking, instructs matchers to direct clients towards more supportive resources the more vulnerable they are (according to the score). In LA, the allocation policy recommends that those who score 12 or higher receive PSH, those who score 8-11 receive either PSH or RRH, and those who score 7 or lower receive RRH. Different communities employ variants of this policy with different ``cut-scores'' and potentially having a bucket of (low) scores for which no resource is intended. A salient advantage of this policy is its \emph{interpretability} which is very important when operating in high-stakes settings~\citep{Rudin2019}. Yet, this policy is not tied to outcomes in the sense that it has not been designed with e.g., the goal of being effective nor, in some sense, fair.}

\newpv{In practice, this recommended policy is not fully adhered to due to several reasons of which we discuss the main ones. First, it is not implementable as there are insufficient resources to match the demand in each score ``bucket''. Second, the matcher and case worker, when discussing a particular case during \emph{case conferencing}, may deem that an individual's vulnerability score is not representative of their true vulnerability. Third, resources are matched as they arrive and would never ``wait'' for a more vulnerable individual. Finally, to date matchers are only responsible for matching PSH resources whereas RRH resources are matched by referral through a mechanism that is more ad hoc (this is in the process of changing).

Each time a resource is allocated to an individual, the match is reported in the HMIS together with all interactions that the individual had with the system. In particular, in addition to PSH/RRH matches, other exits from homelessness are recorded, such as if an individual self-resolved (independently obtained private housing) or went to live with family. 


Being sensitive to the fact that it supports communities
disproportionately affected by systemic racial, gender, and other injustices, LAHSA instigated a number of committees in charge of evaluating their system from an equity lens. Most notably, it assembled the Ad Hoc Committee on Women and Homelessness and the Ad Hoc Committee on Black People Experiencing Homelessness~\citep{WAH_LAHSA_2017,BPEH_LAHSA_2018}. An important finding of the latter committee is that in the current system, although resources are allocated at equal rates among Black and White people, Black people return to homelessness at twice the rates of Whites. Similar findings have been reported in other communities~\citep{Petry_2021,Hill2021,Hsu2021}. Motivated by these findings and by the desire to make more effective use of its resources, LAHSA is interested in building a more evidence based system that will effectively allocate resources while better aligning with the values of policy-makers in serving its clients with a focus on equity.



}

\subsection{\newpv{Trade-Offs and Policy Attributes}}
\label{sec:policy_attributes}

\newpv{A housing allocation policy is a function $\pi : \sets X \rightarrow \Delta^{|\sets T|}$ that maps an individual's characteristics $x \in \sets X$ to their probability of getting each intervention in the set $\sets{T}:=\{\text{PSH}, \text{RRH}, \text{SO}\}$, where SO stands for ``Services Only'', i.e., no housing intervention. We focus our attention on policies that are \emph{feasible}, in the sense that the expected number of clients we assign to PSH and RRH does not exceed the number of these resources that are available (according to this definition, the current recommended policy is not feasible). From the above discussion, the key characteristics that drive the choice of a policy are its interpretability, effectiveness, and fairness properties. We now formalize these mathematically.}

\begin{figure}
    \begin{minipage}[c]{0.5\textwidth}
    \caption{\newpv{Sample prescriptive tree to allocate housing resources based on answers to the VI-SPDAT survey questions. The numbers at each leaf correspond to the probability of getting each intervention for all individuals that reach that leaf.}}\label{fig:prescriptive_tree}
    \end{minipage}\hfill
  \begin{minipage}[c]{0.4\textwidth}
    \includegraphics[width=0.9\textwidth]{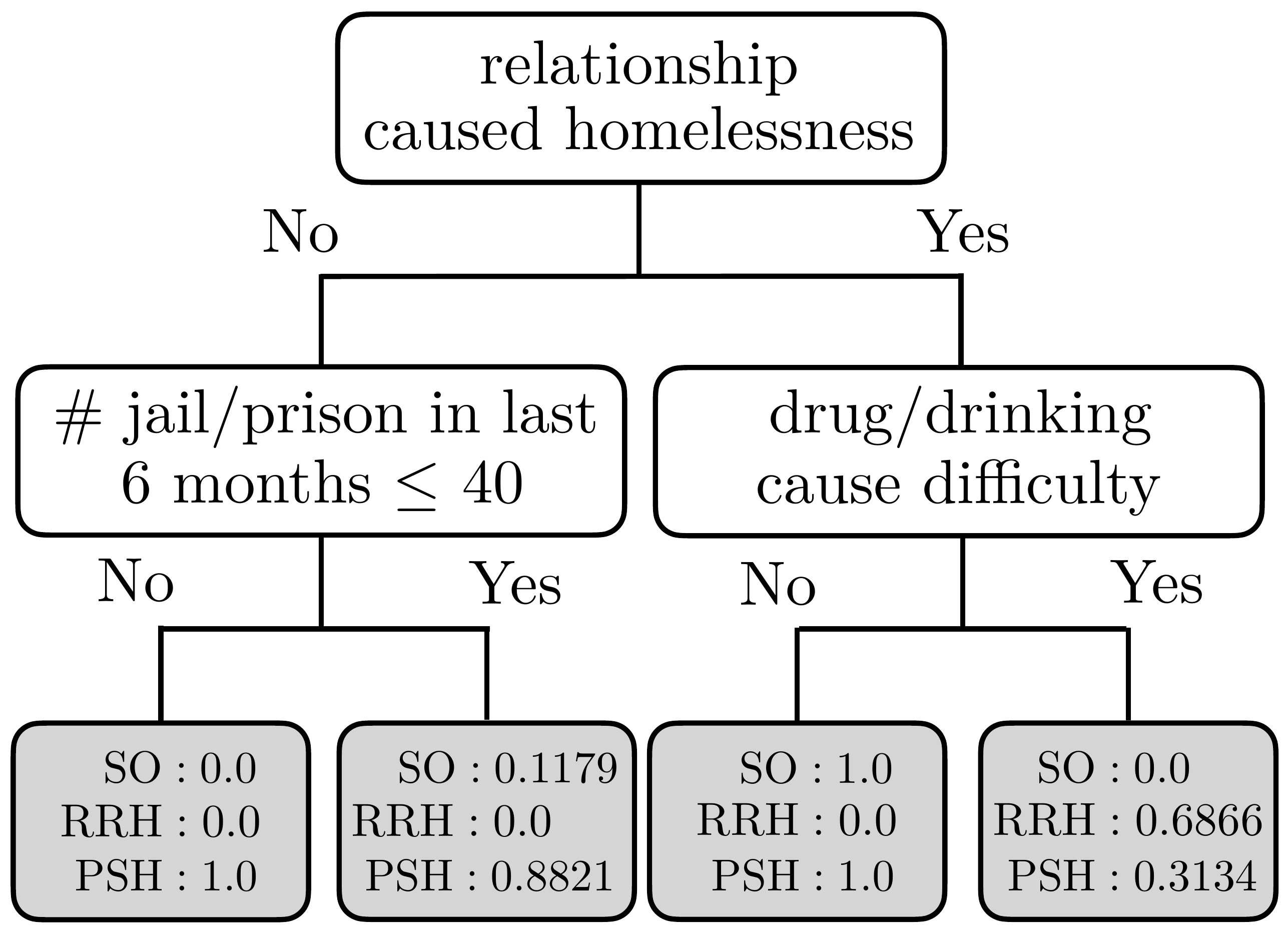}
  \end{minipage}
\end{figure}

\paragraph{Prescriptive Trees and Measures of Interpretability.} \newpv{The first step in designing a policy is agreeing on one (or several) policy classes  to consider \citep{Azizi2018_CPAIOR}. Motivated by the desire to have interpretable policies, we focus here on prescriptive trees~\citep{Kallus2017RecursiveData,Bertsimas2019OptimalTrees,Jo_2021}. These take the form of a binary tree. In each branching node of the tree, a binary test is performed on a specific feature. Two branches emanate from each branching node, with each branch representing the outcome of the test. If a datapoint passes (resp.\ fails) the test, it is directed to the right (resp.\ left) branch. At each leaf node, each treatment is assigned with a certain probability. Thus, each path from root to leaf represents a rule that assigns each treatment with the same probability to all individuals that reach that leaf. An example of a prescriptive tree can be seen in Figure~\ref{fig:prescriptive_tree}. Prescriptive trees are among the most interpretable treatment rules and are thus particularly attractive when designing social and public health interventions. We measure interpretability of a prescriptive tree using two attributes: the number of branching nodes and the number of features it uses (the greater these are, the less interpretable the policy).
}

\paragraph{\newpv{Measure of Effectiveness.}} \newpv{To measure effectiveness of a policy, we must use outcomes that policy-makers care about and that are also observable in the HMIS. According to the U.S.\ Department of Housing and Urban Development (HUD), an important outcome to mitigate is \emph{chronic homelessness}. While this is not directly observable in the data, a proxy that can be observed is \emph{housing (in)stability}. Using potential outcomes notations, see~\cite{hernan2019causal}, we let $Y(t) \in \{0,1\}$ denote the outcome under treatment $t\in \sets T$--thus $Y(t)=1$ iff the individual will still be housed after one year of receiving treatment $t$. Note that individuals who do not receive a PSH nor RRH resource may still have a positive outcome if e.g., they self-resolved. Mathematically, we can express effectiveness as the expected number of people that are still housed under the given policy, i.e.,
\begin{equation}
    Q(\pi) \; := \; \mathbb E\left[  \sum_{t \in \sets T} \pi( t | X ) Y(t) \right],
    \label{eq:policy_efficiency}
\end{equation}
where the expectation is taken over the joint distribution of $X,Y(1),\ldots,Y(T)$.
}

\paragraph{Measures of Fairness.} \newpv{To measure fairness, we focus on metrics that have been investigated by the ad hoc committees charged with evaluating equity in the system, see Section~\ref{sec:LA_housing_system}. These look at outcomes by group with respect to various protected characteristics. Letting $\sets P$ denote the set of protected characteristics and $\sets X_p$ the set of values that this covariate can take, we can express the conditional average outcome for each group as
\begin{equation}
    Q_p(\pi,x_p)\; := \; \mathbb E\left[ \left. \sum_{t \in \sets T} \pi( t | X ) Y(t) \right| X_p = x_p \right] \text{ for all } p \in \sets P \text{ and } x_p \in \sets X_p,
    \label{eq:policy_fairness}
\end{equation}
where $X_p \subseteq X$ represents the value of protected characteristic~$p$. An associated metric that we also consider an important attribute in the choice of policy is the smallest conditional average outcome for each protected characteristic, i.e., $\min_{x_p \in \sets X_p} Q_p(\pi,x_p)$ for all $p \in \sets P$. This is motivated by the notion of max-min fairness that seeks to provide maximum benefit to the least advantaged members of society~\citep{Rawls1991}. The last attribute related to fairness that we consider is the number of protected features we branch on (with smaller numbers being preferred).
} 


\subsection{Using our Preference Elicitation Framework}
\label{sec:using_our_framwork}

\newpv{To use our framework in practice several candidate policies must be identified and their fairness, effectiveness, and interpretability attributes must either be calculated (e.g., number of branching nodes) or estimated from data (e.g., outcomes by race). We now describe in more detail the data we have and the process we use to estimate the various fairness and effectiveness metrics.}

\paragraph{Raw Data Description.} \newpv{The raw dataset we use features HMIS data from 16 Continuums of Care that were sourced, accessed, anonymized, and provided to us by OrgCode Consulting, Inc. The dataset includes the demographic characteristics and VI-SPDAT responses for $I=22,165$ unsheltered single adults assessed between February 2015 and April 2018. Demographic characteristics include individuals’ self-reported age, gender, and race or ethnicity. We collect demographic information and VI-SPDAT responses in the feature vector $X_i \in \sets X \subseteq \mathbb R^F$, $i\in \sets I:=\{1,\ldots,I\}$. The data also includes information on the intervention $T_i \in \sets T$ received by individual $i \in \sets I$. Housing interventions include exits from homelessness to either PSH, RRH, family, or self-resolve. In our analysis, we group family and self-resolve interventions under the SO category. Finally, the data reports the observed outcome $Y_i$, which indicated if the client is still housed 365 days after their initial exit from homelessness, see~\cite{Petry_2021} for details.}. 

\paragraph{\newpv{Protected Characteristics and Groups.}} \newpv{We measure fairness with respect to race/ethnicity, age, and gender. We group individuals into three racial/ethnic groups, namely White, Black, and Other, where the Other category represents 5.9\% of the population and includes those self reported as Asian, Hispanic, Native American, and Pacific Islander. Similarly, we group individuals into four age groups, those aged [25,41], (41,48], (48,54], and (54,84], such that each group has approximately the same size. Finally, there were three self reported gender groups--male, female, and transgender. Thus, $\sets P := \{ \text{race}, \text{gender}, \text{age}\}$ and e.g., $\sets X_{\text{race}}:=\{ \text{White}, \text{Black}, \text{Other}\}$.}

\paragraph{\newpv{Policy Attributes.}} \newpv{This approach of grouping individuals gives us a total of $J=17$ features for each policy: 1 related to efficiency, 2 to interpretability, and the remaining to fairness.}

\paragraph{\newpv{Potential Outcomes and Assumptions.}} \newpv{The HMIS data is an observational dataset where housing interventions are not allocated at random (as in randomized clinical trials). We view the HMIS data in a potential outcomes framework. In particular, the joint distribution of $X$ and $\{Y(t)\}_{t \in \sets T}$ is unknown and we can only observe, for each point $i$ in the data the outcome $Y_i = Y_i(T_i)$ under the treatment received. Critically, we cannot control (and do not know) the historical treatment assignment policy and the outcomes $Y_i(t)$ for $t \neq T_i$ remain \emph{unobserved}. We assume however, that the conditional probability of receiving a treatment depends only on the covariates $X$ and that the probability of receiving any given treatment conditional on $X$ is positive, i.e., $\mu(t,x) :=\mathbb P(T=t|X=x)>0$ almost surely for all $t$. These two assumptions, standard in causal inference, will allow us to use the available data to estimate the performance of counterfactual housing allocation policies. While not possible to verify them from the data, these assumptions mostly hold based on what we know about the system (indeed, while some data from case conferencing that is not in HMIS may have been used to assign resources, this seems to be a rare occurence).}

\paragraph{Estimating the Performance of Counterfactual Policies.}  \newpv{Evaluating several of the fairness and efficiency metrics such as~\eqref{eq:policy_efficiency} and~\eqref{eq:policy_fairness} requires knowledge of the joint distribution of $X$ and $\{Y(t)\}_{t \in \sets T}$. While this distribution is unknown, we can leverage the assumptions above to estimate them from the historical observations $\{(X_i, T_i, Y_i)\}_{i \in \mathcal{I}}$. In our work, we propose to use a doubly robust estimator, see~\cite{Dudik2011DoublyLearning}. Specifically, we estimate $Q(\pi)$ an $Q_p(\pi,x_p)$ through
\begin{equation}
%
\hat Q(\pi) := \frac{1}{I}\sum_{i \in \mathcal I} \left(\hat{\nu}_{\pi(X_i)}(X_i)+ (Y_i - \hat{\nu}_{T_i}(X_i)) \frac{ \pi(T_i | X_i) }{\hat \mu(T_i, X_i)}\right)
\label{eq:policy_value_estimate}
\end{equation}
and
\begin{equation}
\hat Q_p(\pi,x_p) := \frac{1}{|i \in \mathcal I : X_{i,p} = x_p|}\sum_{i \in \mathcal I : X_{i,p} = x_p} \left(\hat{\nu}_{\pi(X_i)}(X_i)+ (Y_i - \hat{\nu}_{T_i}(X_i)) \frac{ \pi(T_i | X_i) }{\hat \mu(T_i, X_i)}\right),
\label{eq:conditional_policy_value_estimate}
\end{equation}
respectively. Here, $\hat \mu$ is an estimator of $\mu$, which is obtained for instance using machine learning, by fitting a model to $\{(X_i,T_i)\}_{i \in \sets I}$, and $\hat{\nu}_{t}(x)$ is and estimator of $\mathbb{E}(Y|T=t, X=x)$ obtained by fitting a model on the subpopulation that was assigned treatment $t$. Provided at least one of $\hat \mu$ or $\hat \nu$ converges almost surely to $\mu$ or $\nu$, respectively, these doubly robust estimators of $Q(\pi)$ and $\hat Q_p(\pi,x_p)$ will be asymptotically unbiased. All fairness and efficiency metrics of a policy $\pi$ can then be estimated using the two estimators above. 
}

\subsection{Numerical Case Study}

\newpv{We now summarize our findings. We first discuss our models of $\hat \nu_t(x)$, $t\in \sets T$, and $\hat \mu(t,x)$ and then present the results of the offline and online elicitation procedures.}

\newpv{Since we use estimation procedures to evaluate the counterfactual performance of a policy and in particular its fairness properties, it is imperative that our learned models of $\hat \nu_t(x)$ and $\hat \mu(t,x)$ in \eqref{eq:policy_value_estimate} and~\eqref{eq:conditional_policy_value_estimate} be fair. In particular, they should have the same true positive rates (TPR) and same false positive rates (FPR) across groups for each protected group and feature. This notion of fairness is known as \emph{equalized odds}, see~\citet{Hardt:2016}. We calculate the equalized odds fairness metric by taking the maximum, across all protected characteristics, of the greatest difference in both TPR and FPR between all groups. Thus, this is a number between 0 and 1 and we would like it to be as small as possible to indicate a fair classifier. We use the approach from~\cite{Aghaei_StrongCT} to learn optimal fair classification trees--we denote this by \texttt{OFT}. This framework has the advantage that it performs close to the best unfair models (e.g., CART, random forests) in terms of accuracy when no fairness constraints are imposed, that it applies to multi-class classification problems, and to problems with multiple protected groups and multiple levels for each group.}

\begin{figure}[h!]
    \centering
    \includegraphics[height=0.3\textwidth]{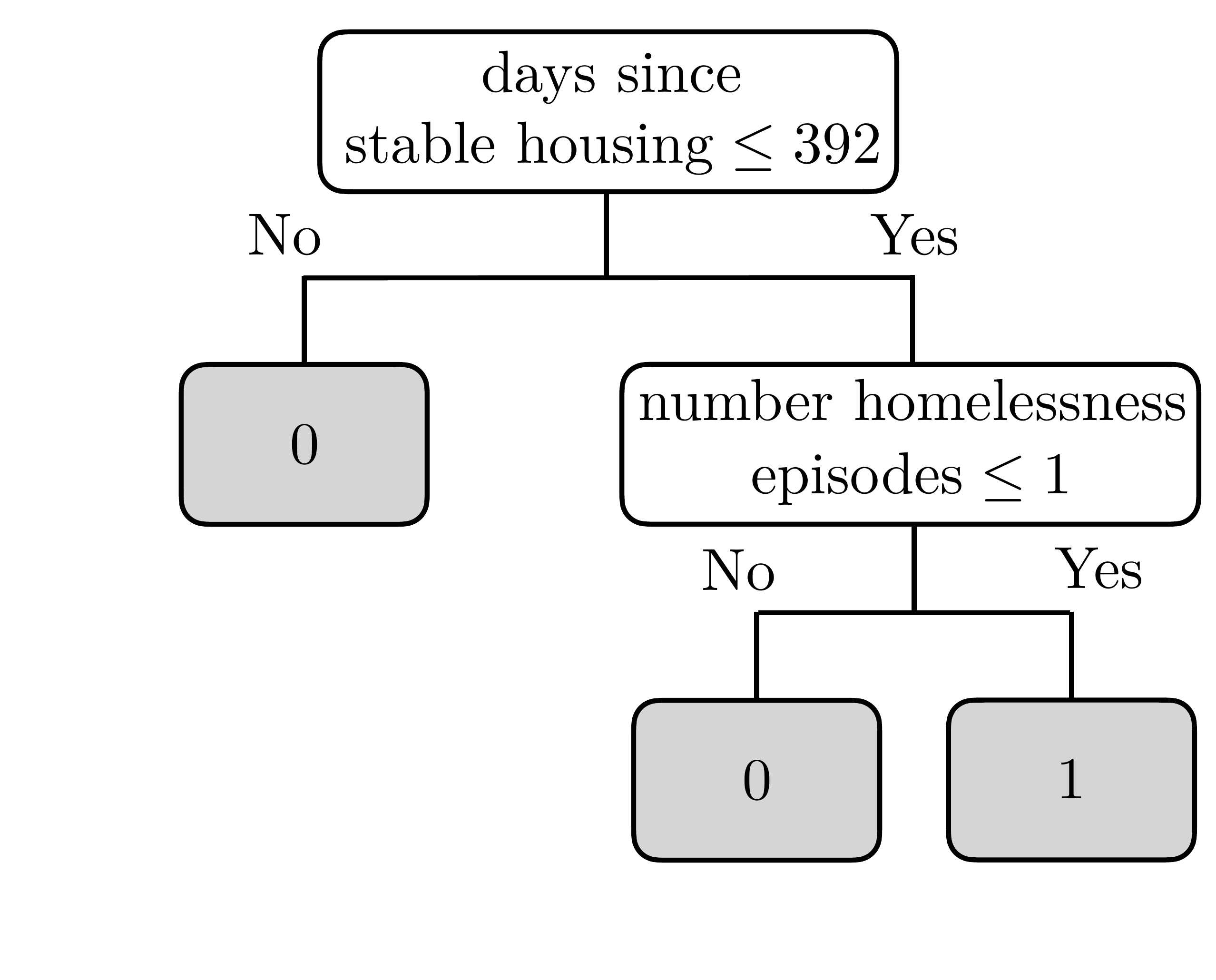}
    \hspace{1cm}
    \includegraphics[width=0.45\textwidth]{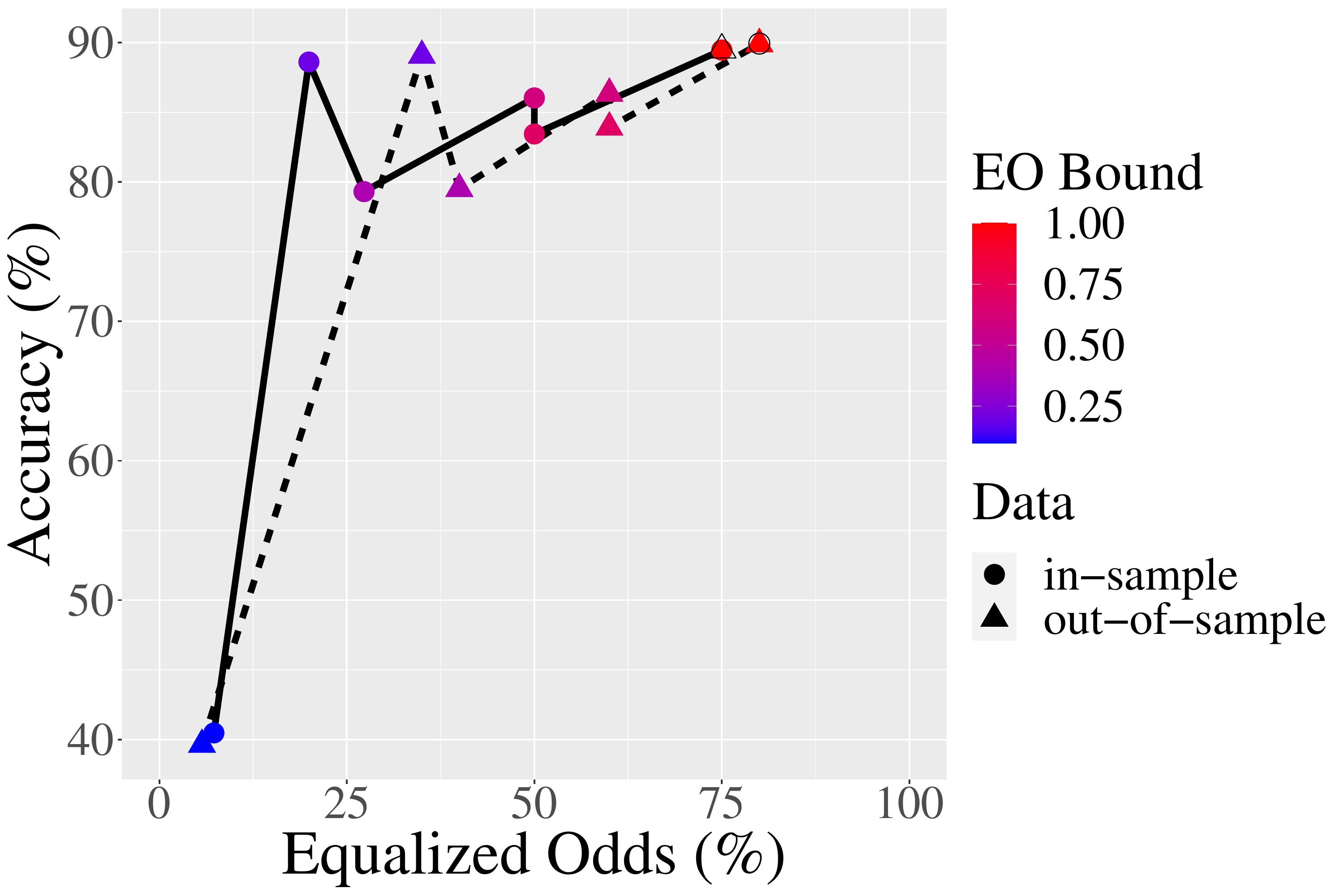} \\ \vspace{0.5cm}
    \includegraphics[height=0.3\textwidth]{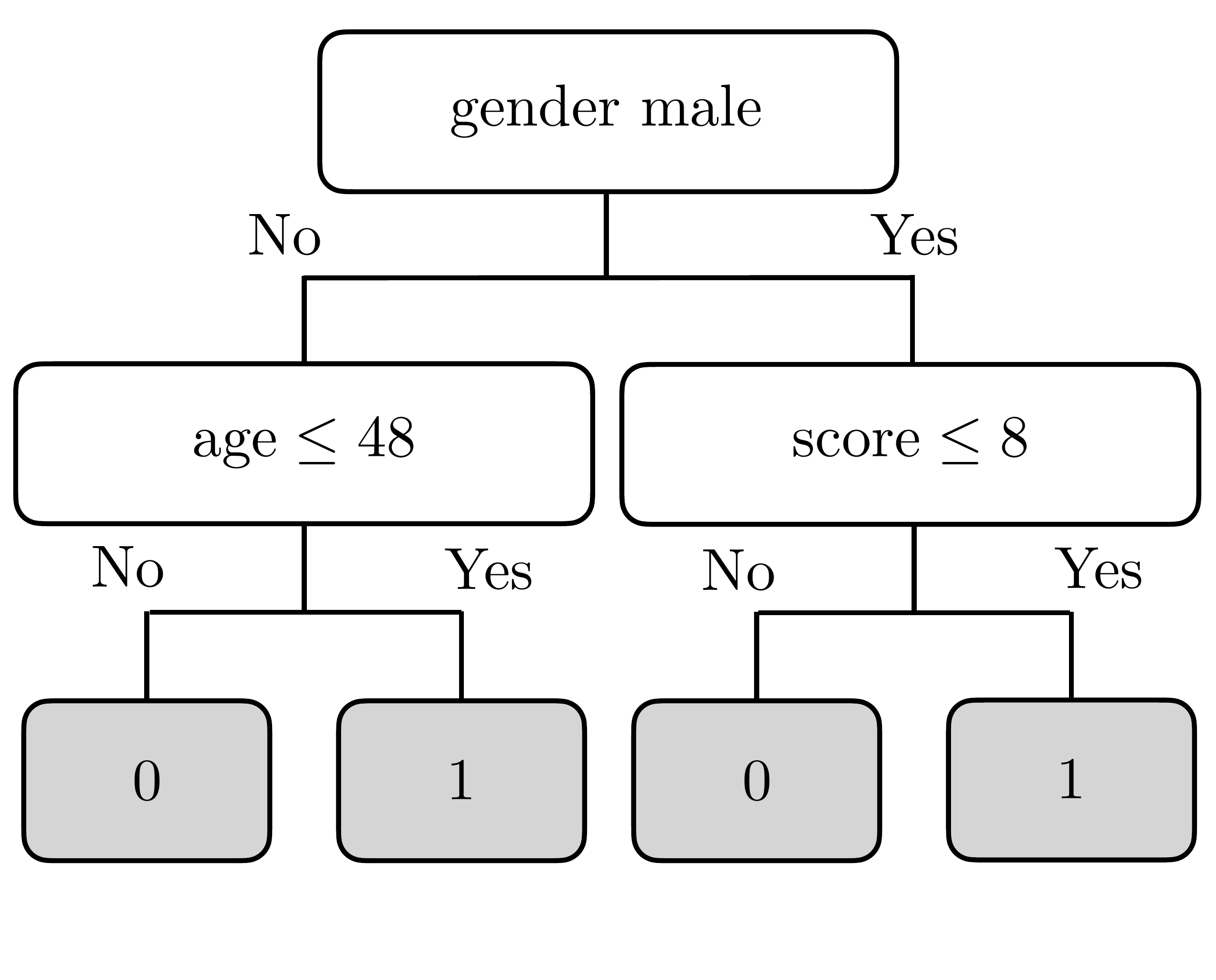} \hspace{1cm}
    \includegraphics[width=0.45\textwidth]{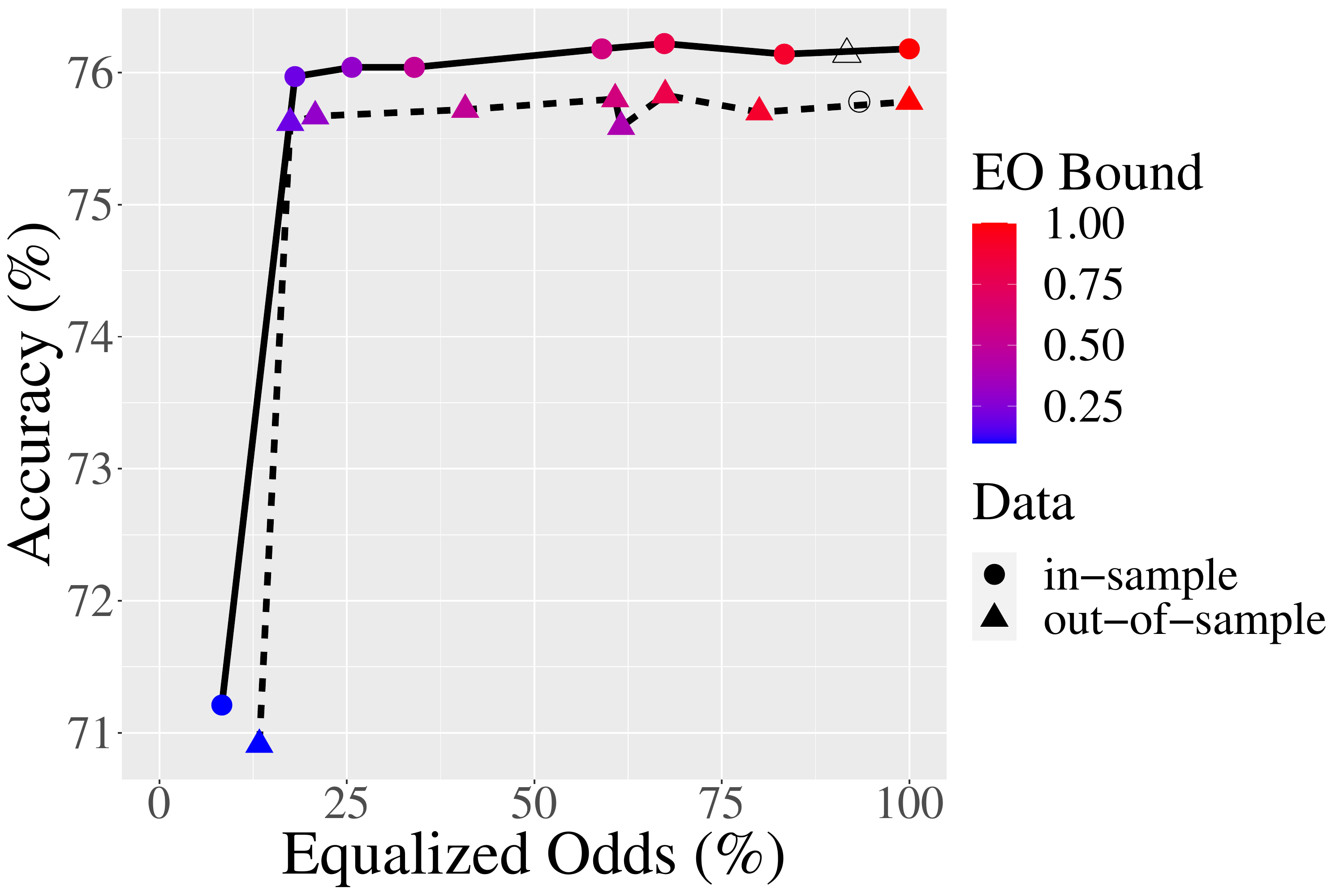} \\
    \vspace{0.5cm}
    \includegraphics[height=0.3\textwidth]{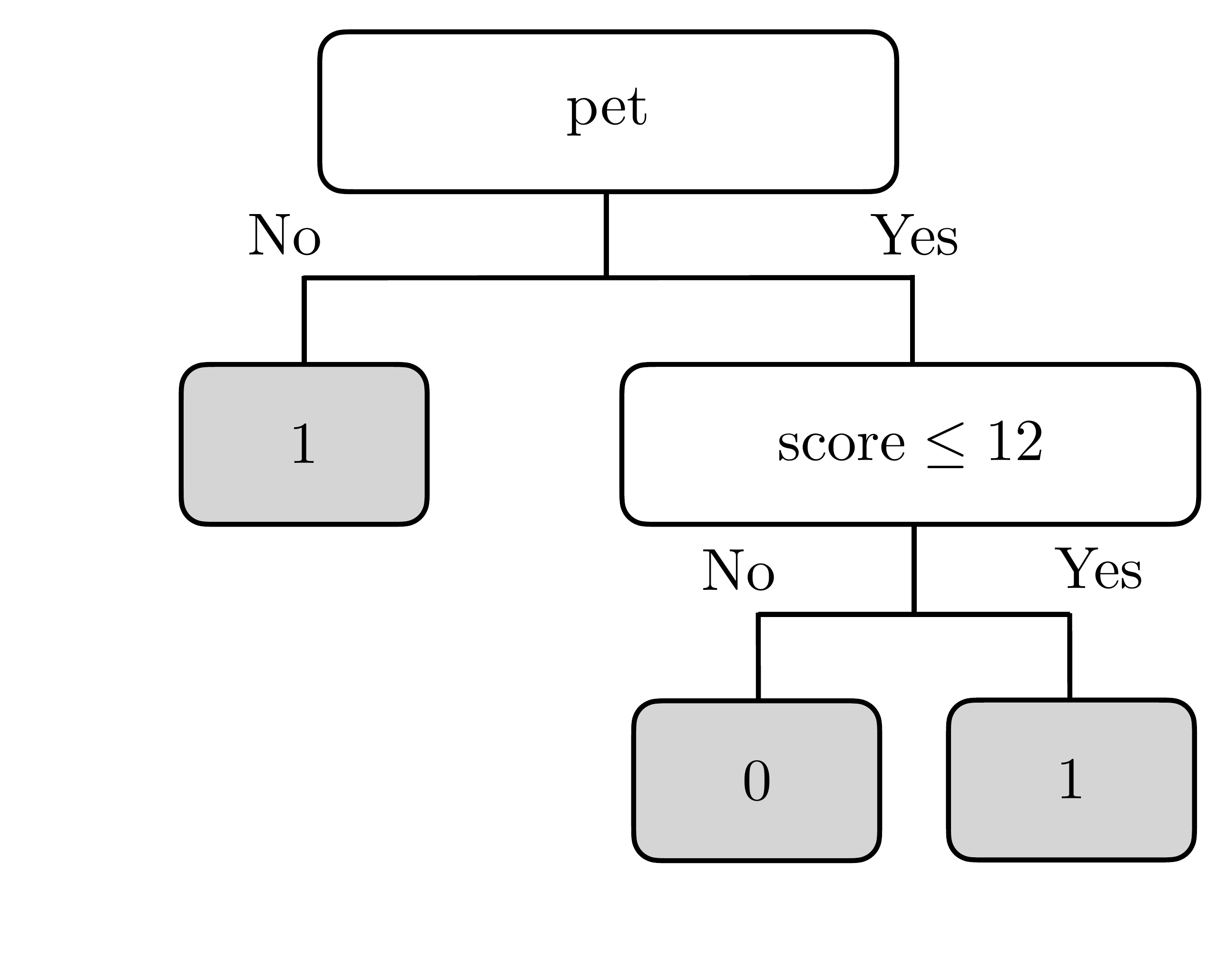}
    \hspace{1cm}
    \includegraphics[width=0.45\textwidth]{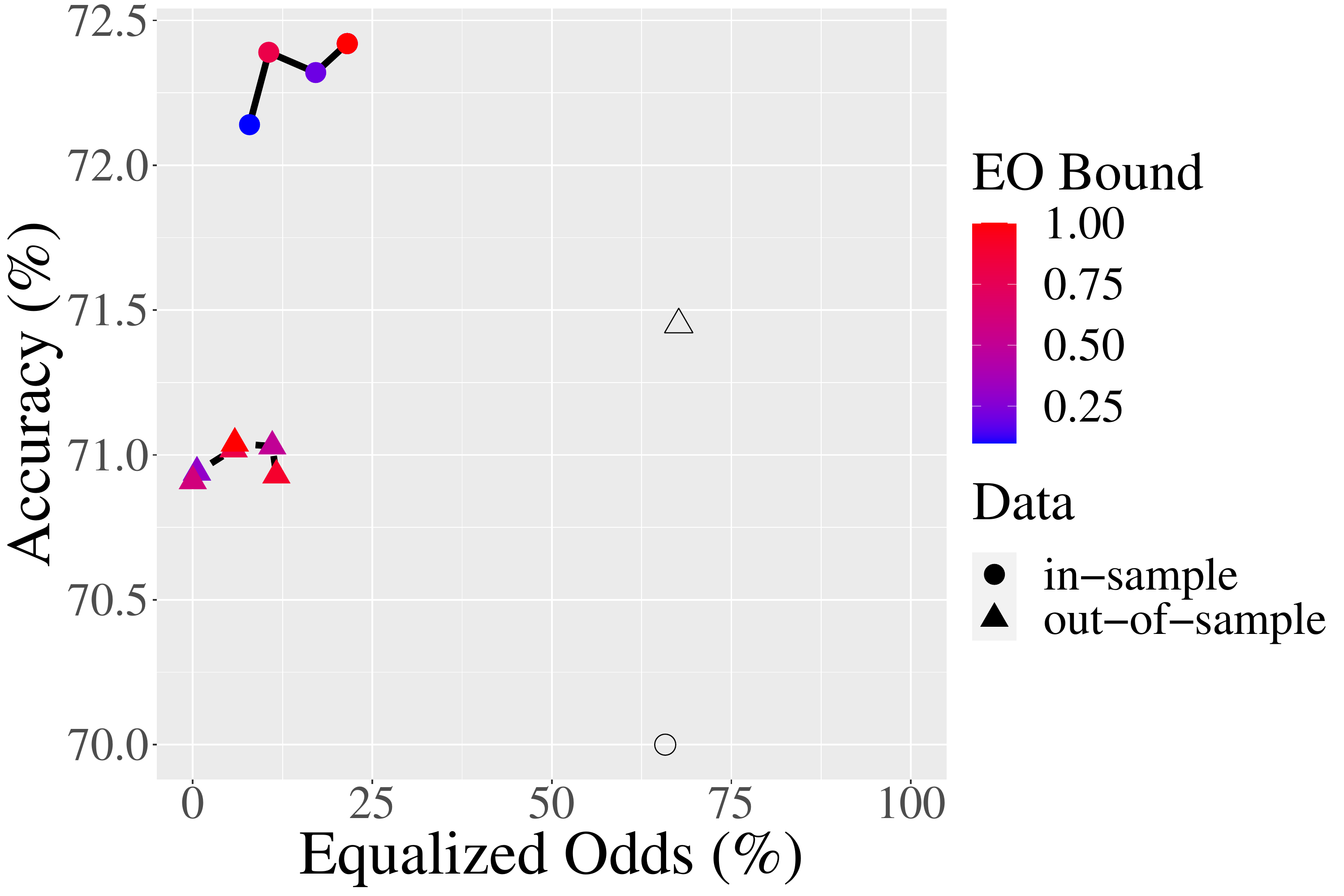}
    \caption{\newpv{Optimal and fair decision trees for estimating $\nu_t(x) = \mathbb E(Y | X = x, T=t)$ for $t\in \{ \text{SO}, \text{RRH}, \text{PSH}\}$ (top, middle, and bottom, respectively). The figures on the right show the fairness/accuracy trade-off both in-sample and out-of-sample for CART (empty circle and triangle) and for \texttt{OFT} (filled circles and triangles). The figures on the left show the optimal fair classification trees that we chose to use; these corresponded to the trees that were the most fair among those that resulted in minimal loss in accuracy.}} 
    \label{fig:opt_fair_DT_outcomes}
\end{figure}

\newpv{\paragraph{Fair Learning of $\hat \nu_t(x)$.} We used several models to learn outcomes $\hat \nu_t(x) = \mathbb E(Y|X=x, T=t)$ for $t\in \sets T$. Results of the best models are shown on Figure~\ref{fig:opt_fair_DT_outcomes}. As shown in the figure, in both training and testing the best CART model had equalized odds of over 75\% showcasing the importance of imposing fairness directly in the model. In contrast, by using \texttt{OFT}, we can see that in all cases we could significantly decrease equalized odds of the learned models compared to CART at no cost to accuracy. For $t \in \{\text{SO}, \text{RRH}\}$, decreasing equalized odds beyond about 20\% results in significant drop in accuracy. For our predictions, we chose the model associated with the lowest possible equalized odds that did not result in significant drop in accuracy. For the case of the SO model, rather than optimizing accuracy, we optimized \emph{balanced} accuracy since this dataset was significantly imbalanced. From the figure, we can see that the learned models are interpretable and also make a lot of sense. For example, in the case of the SO intervention, only those recently stably housed with no more than one episode of homelessness are predicted to have a positive outcome.

}


\begin{figure}[h!]
    \centering
    \includegraphics[height=0.3\textwidth]{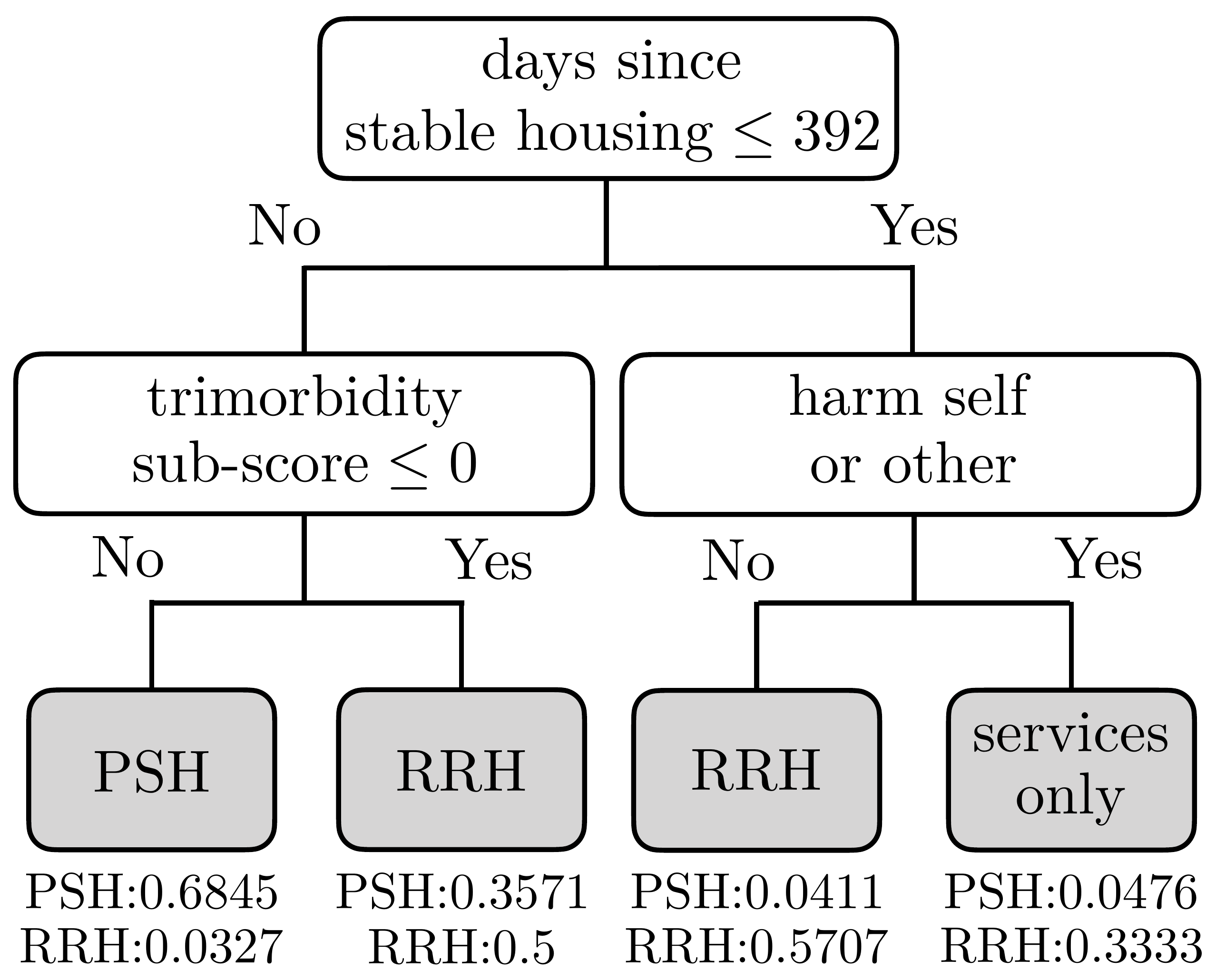} \hspace{1cm}
    \includegraphics[width=0.45\textwidth]{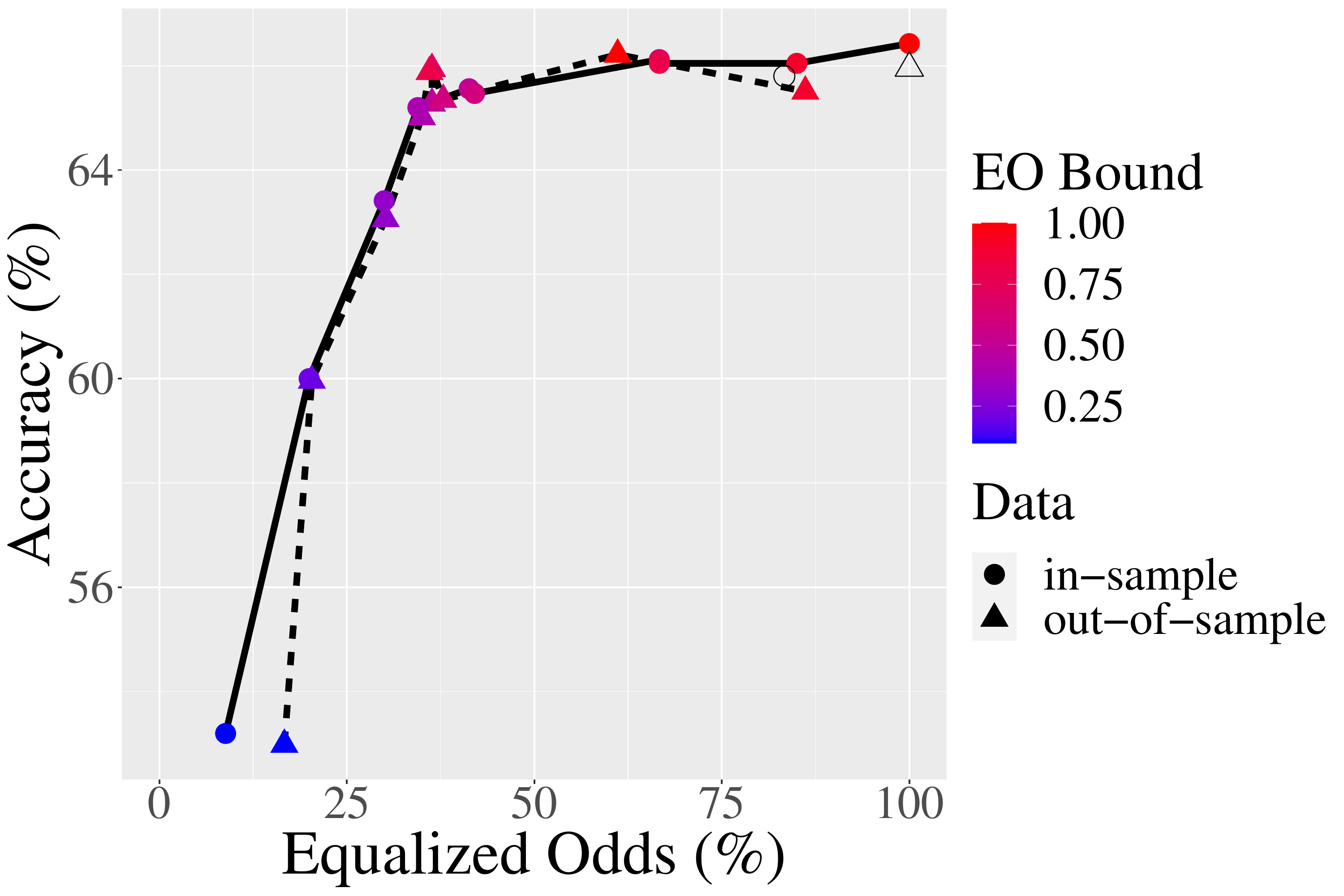}
    \caption{\newpv{Optimal and fair decision-tree for predicting $\hat \mu(t,x)$. The meaning of the circles and triangles is the same as in Figure~\ref{fig:opt_fair_DT_outcomes}. The tree chosen to predict propensity scores is shown on the left.}}
    \label{fig:opt_fair_DT_ps}
\end{figure}

\paragraph{Fair Learning of $\hat \mu(t,x)$.} \newpv{Learning the propensity score is a multi-class classification problem. We generalize the equalized odds fairness metric to multiple classes as the greatest difference in error rates of any type across groups. Our results on learning $\hat \mu(t,x)$ are summarized in Figure~\ref{fig:opt_fair_DT_ps}. From the figure, we see that the best CART model results in equalized odds of over 80\%, whereas using \texttt{OFT} allows us to carefully tune the trade-off between fairness and accuracy and in particular to reduce equalized odds to under 38\% at no cost to accuracy. Moreover, we see that the model obtained to fairly predict propensity score is very interpretable and natural: the more the indications of vulnerability an individual has, the more likely they are to get more supportive interventions.}

\paragraph{Generating Feasible Counterfactual Prescriptive Tree-Based Policies.} \newpv{To generate feasible counterfactual prescriptive trees with different fairness, efficiency, and interpretability properties, we proceed as follows. For each policy we wish to design (25 in total), we sample a utility vector ${\bm u}$ uniformly at random from the set $\sets U = \{ {\bm u} \in \mathbb R^{17}_+ : \textbf{\textrm{e}}^\top {\bm u} = 1\}$. We also draw the depth~$d$ of the tree uniformly from the set $\{1,2,3\}$. Then, we use the method from~\cite{Jo_2021} to design an optimal and feasible tree based policy of depth $d$ that optimizes the combination of the aforementioned metrics of fairness, efficiency, and interpretability with weight vector ${\bm u}$. Finally, we record the relevant attributes of the policy, resulting in a dataset with $I=25$ policies and $J=17$ attributes.
}

\begin{figure}
    \begin{minipage}[c]{0.45\textwidth}
    \caption{\newpv{Optimality (top) and scalability (bottom) results for the offline min-max regret preference elicitation problem~\eqref{eq:offline_mmr} on the real data when $\sigma=0$ (left) and $\sigma=0.025$ (right). Approach \texttt{MMR} is shown with red dots. The median performance of \texttt{RAND} across 50 sets of $K$ random queries is shown with the blue dashed line. In particular the run time for \texttt{RAND} is the time it takes to evaluate the performance of one query. The blue shaded region shows the range of performance of \texttt{RAND} across these 50 sets of queries.}} \label{fig:offline_mmr_data}
    \end{minipage}\hfill
  \begin{minipage}[c]{0.5\textwidth}
    \includegraphics[width=\textwidth]{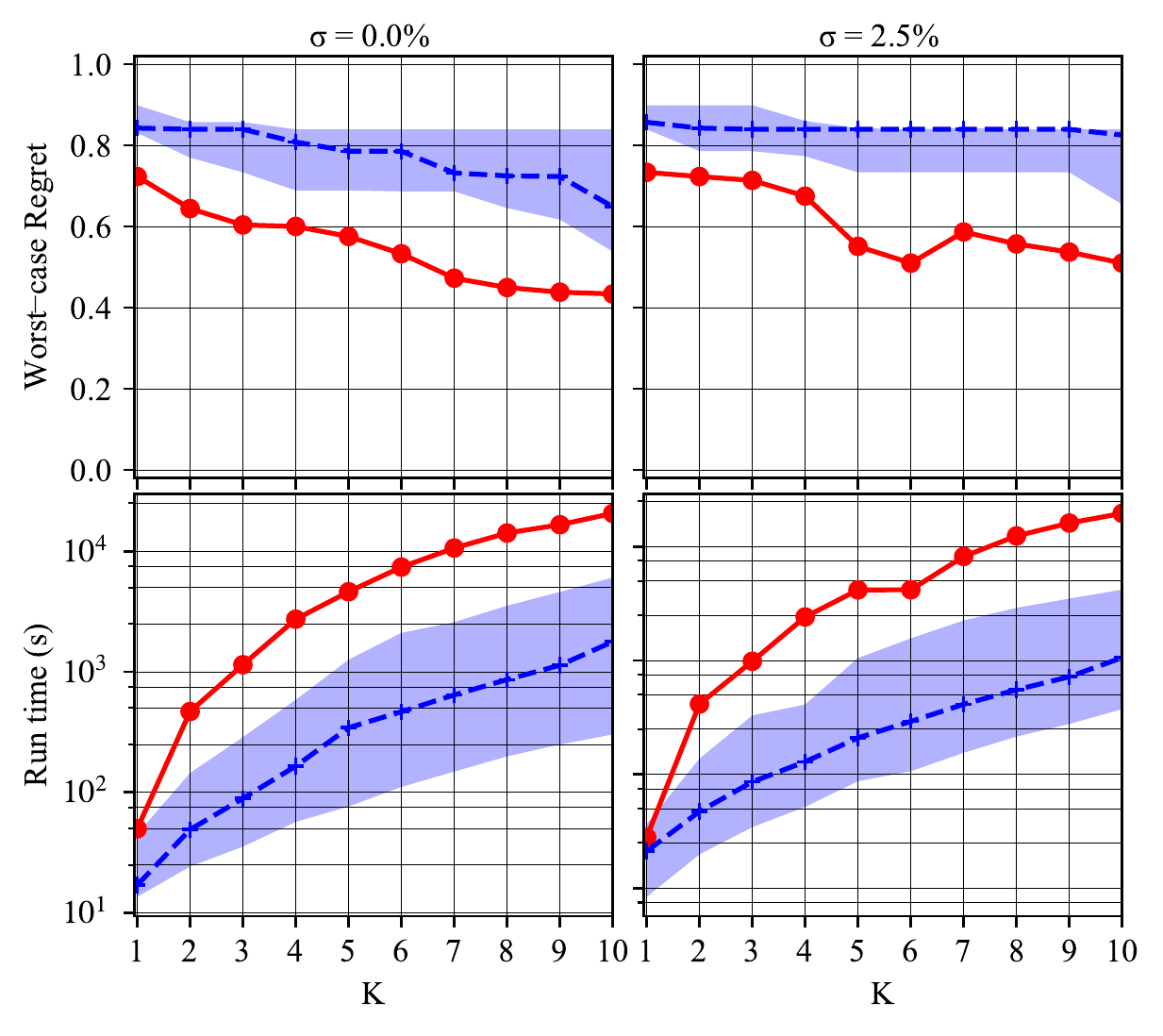}
  \end{minipage}
\end{figure}

\paragraph{\newpv{Results of Offline \texttt{MMR} Elicitation.}} \newpv{We focus on the min-max regret case as this usually results in less conservative solutions. The results of offline \texttt{MMR} on this dataset of candidate prescriptive tree policies are shown on Figure~\ref{fig:offline_mmr_data}. From the figure, we see that our approach consistently outperforms random elicitation across all numbers of queries in terms of worst-case regret. In particular, with just a single query, our approach is able to outperform the average performance of random elicitation with 10 queries. This is true whether responses given by the user are inconsistent or not. In addition, optimizing over queries is faster than evaluating the performance of all possible queries.}

\begin{figure}
    \begin{minipage}[c]{0.45\textwidth}
    \caption{\newpv{Worst-case regret (top), worst-case true regret (middle), and worst-case rank (bottom) for the online min-max regret preference elicitation problem on real data when $\sigma=0$ (left) and $\sigma=0.025$ (right). The median performance of \texttt{MMR} (resp.\ \texttt{RAND}, \texttt{POLY}, \texttt{PROB}, \texttt{ROB}, \texttt{ELL}) across 50 random utility vectors ${\bm u}$ and inconsistencies ${\bm \epsilon}$ is shown with red dots {(resp.\ blue plus, black triangles, black inverted triangles, green crosses, and purple squares).}}}\label{fig:online_mmr_data}
    \end{minipage}\hfill
  \begin{minipage}[c]{0.5\textwidth}
    \includegraphics[width=\textwidth]{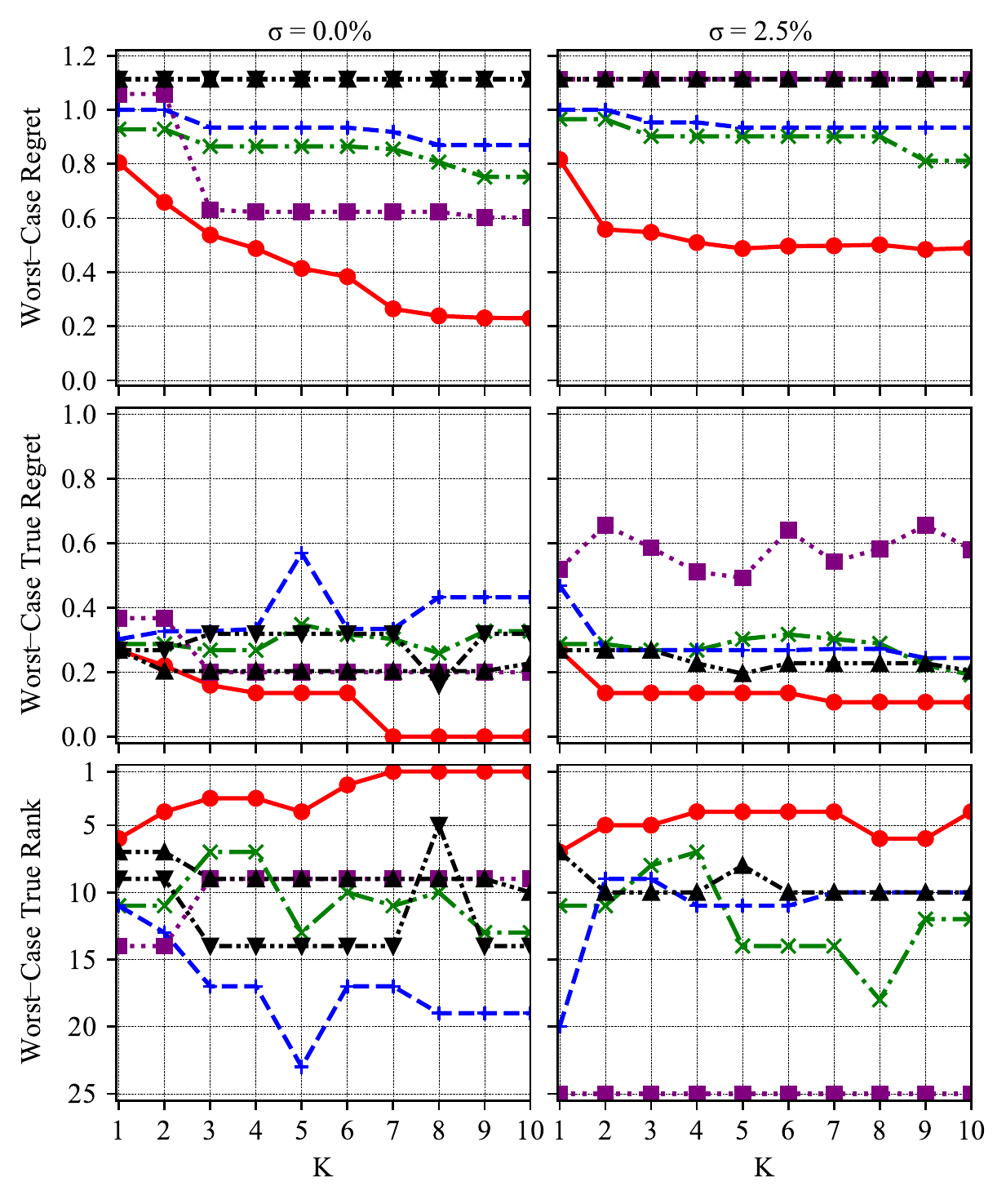}
  \end{minipage}
\end{figure}

\paragraph{\newpv{Results of Online \texttt{MMR} Elicitation.}} \newpv{Since in the online case it is not possible to compute the objective value of the problem, we evaluate our method through simulation by sampling 50 agent utilities from the uncertainty set and allowing potential inconsistencies in their responses. The results are summarized in Figure~\ref{fig:online_mmr_data}. From the figure, it can be seen that our approach consistently outperforms all other methods from the literature in terms of worst-case regret, worst-case true regret of the recommendation, and worst-case true rank of the recommendation. In particular, after just 7 queries, our method is able to consistently recommend the best item if there are no inconsistencies and an item ranked 6 or higher if there are inconsistences.}

\paragraph{\newpv{Robustness to Choice of $\Gamma$.}} \newpv{A challenge when working with real data is that our model may be misspecified, e.g., if the standard deviation $\sigma$ used to select queries is different than the true standard deviation $\sigma^\star$ of the shocks in the user utility, see Section~\ref{sec:model}. We perform an evaluation of our model on a range of such misspecifications. The results are summarized in Table~\ref{tab:sigma_comparison}. From the table, we see that our query selection methods are quite robust to our modeling assumptions. In particular, worst-case regret, worst-case true regret, and worst-case true rank only mildly deteriorate if the model is different than assumed.}


\begin{table}[ht]
\caption{\newpv{Performance of our approach when the standard deviation~$\sigma$ of the inconsistencies used to select $\Gamma$ in the optimization is different than the true standard deviation $\sigma^\star$ of the inconsistencies. Numbers are averages across values of~$K$, $\sigma$, and $\sigma^\star$. The detailed results can be found in Tables~\ref{tab:offline_sigma_comparison_detail} (offline) and~\ref{tab:online_sigma_comparison_detail} (online).}}
\centering
\renewcommand{\arraystretch}{1.1}
\begin{small}
\begin{tabular}{|C{2.5cm}||C{2.5cm}|C{2.5cm}|C{2.5cm}|}
  \hline
Problem Type & Avg.\ Difference in Normalized Worst-Case Regret & Avg.\ Difference in Normalized Worst-Case True Regret & Avg.\ Difference in Worst-Case True Rank \\
  \hline
     Offline & 0.03 & 0.03 & 2.44 \\ 
     Online & 0.00 & 0.01 & 0.48 \\ 
    
   \hline
\end{tabular}
\end{small}
\label{tab:sigma_comparison}
\end{table}

\section{Conclusion}
\label{sec:conclusion}

In this paper, \newpv{motivated by the problem of making recommendations in high stakes settings where hard trade-offs need to be made,} we proposed novel formulations for the offline and online \newpv{max-min utility and min-max regret} active preference elicitation problems. These take the form of robust optimization problems with decision-dependent information discovery. We studied the complexity of these problems and provided exact reformulations and conservative approximations of the offline and online problems, respectively. We provided efficient solutions procedures and performed extensive computational experiments that showed the superiority of our approach on both synthetic data and real data from the homeless management information system. In the future, we plan to deploy this algorithm on Amazon Mechanical Turk\footnote{\url{https://www.mturk.com}} and on policy-makers at the Los Angeles Homeless Services Authority to be able to elicit their preferences and recommend housing allocation policies that best meet their needs to help mitigate homelessness.






\theendnotes


\ACKNOWLEDGMENT{\newpv{P.\ Vayanos and E.\ Rice thank Dr.\ Iain DeJong for sharing the HMIS data used our numerical study. They also gratefully acknowledge support from Schmidt Futures and from the Hilton C.\ Foundation, the Homeless Policy Research Institute, and the Home for Good foundation under the ``C.E.S.\ Triage Tool Research \& Refinement'' grant. P. Vayanos is also grateful for the support from the James H. Zumberge Faculty Research and Innovation Fund at the University of Southern California, from the U.S. Department of Transportation METRANS University Transportation Center and National Center for Sustainable Transportation under Award No.\ DTRT13-G-UTC57, and from the National Science Foundation under CAREER award number 2046230. The authors are grateful to Marina Genchev, Sr.\ Manager for the Adult Coordinated Entry System at LAHSA, Erin Cox, former Policy Manager at LAHSA, and Maggie Potthoff, former Policy Supervisor at LAHSA, for their support of this work and for the very interesting discussions that helped shape this research. They would also like to thank Andr\'{e}s G\'{o}mez for vibrant conversations on polyhedral theory and Sina Aghaei and Nathan Jo for guidance on using and adapting their optimal classification trees and optimal prescriptive trees codes, respectively.}}



\bibliographystyle{informs2014} 
\bibliography{Mendeley.bib}

\newpage



\ECSwitch


\ECHead{E-Companion}


\section{Illustration of the Uncertainty Set Update Procedure} 
\label{sec:EC_polyhedral_method}

\newpage


\section{Generating Warm-Starts Satisfying Lexicographic Constraints} 
\label{sec:EC_speedup_strategies}


Algorithm~\ref{alg:warm_start} details the procedure for generating warm-starts proposed in Section~\ref{sec:speedup_strategies} to ensure that the lexicographic constraints from Section~\ref{sec:symmetry} are satisfied.


\begin{algorithm}[ht!]
\SetAlgoLined
\textbf{Inputs:} A feasible solution $( \tau, {\bm v} , {\bm w} , {\bm \alpha}, {\bm \beta}, {\bm x}, \overline{\bm v}, \overline{\bm w} )$ to Problem~\eqref{eq:offline_mmu_MILP} with $K$ queries with objective $\tau$\;
\textbf{Output:} A warm start $( \tilde{\bm v},\tilde{\bm w}, \tilde{\bm x} )$ to Problem~\eqref{eq:offline_mmu_MILP} with $K+1$ queries with objective $\geq \tau$\;
\textbf{Initialization:}\\
\For{$\kappa \in \sets K$}{
    ${\bm \iota}_1^\kappa \gets \sum_{i \in \sets I} i \cdot \I{{\bm v}^\kappa_i = 1}$\;
    ${\bm \iota}_2^\kappa \gets \sum_{i \in \sets I} i \cdot \I{{\bm w}^\kappa_i = 1}$\;
}
$
\widetilde{\sets C} \leftarrow  \left\{ \; {\bm \iota} \in \sets C \; : \; {\bm \iota} \neq {\bm \iota}^\kappa, \; \kappa \in \sets K \; \right\}
$\;

Select ${\bm \iota}^{K+1}$ at random from $\widetilde{\sets C}$\;
\For{${\bm s} \in \sets S^{K+1}$}{
    $\tilde {\bm x}^{\bm s} \; := \; {\bm x}^{({\bm s}_1,\ldots,{\bm s}_{K})}$\;
}
$( \tilde {\bm \iota} , \bm \ell) \gets \text{lexicographic order}({\bm \iota})$\;
\For{${\bm s} \in \sets S^{K+1}$}{
    $\tilde {\bm x}^{\bm s} \gets \tilde {\bm x}^{ {\bm s}(\bm \ell) }$\;
}
\For{$\kappa \in \sets K$, $i \in \sets I$}{
$\tilde{\bm v}^{\kappa}_i \leftarrow \I{ {\bm \iota}^{\kappa}_1 = i }$\; 
$\tilde{\bm w}^{\kappa}_i \leftarrow \I{ {\bm \iota}^{\kappa}_2 = i }$\;
}
\textbf{Result:} $( \tilde{\bm v},\tilde{\bm w}, \tilde{\bm x} )$

\textit{Note:} The lexicographic order function takes as input a collection of vectors and returns two elements; the first element corresponds to the lexicographically ordered collection of vectors and the second element is a vector of the same dimension as the number of vectors input whose $i$th element denotes the position of the $i$th output vector in the input\;
 \caption{An algorithm for building warm-starts}
 \label{alg:warm_start}
\end{algorithm}


\section{Greedy Solution Approach to Speed-Up Computations}
\label{sec:EC_greedy}

To speed-up solution of our problems so as to be able to showcase performance on a variety of instances, we employ a heuristic approach in all our experiments (unless explicitly stated otherwise). \newpv{This heuristic is a variant of the warm-start approach, see Section~\ref{sec:speedup_strategies} and is detailed in Algorithm~\ref{alg:heuristic}.} This algorithm returns a feasible but potentially suboptimal solution to the MBLP counterpart of the offline preference elicitation problem to be solved. \newpv{The idea of iteratively fixing some of the decision variables in the problem is not new and has a long standing history in operations research, see e.g., \cite{nemhauser1978analysis,BoydBook,vayanos_ROInfoDiscovery,Subrahmanyan_DP_JR}}.

\begin{algorithm}[h!]
\SetAlgoLined
\textbf{Inputs:} 
    Instance of \newpv{Problem~\eqref{eq:offline_mmu} (resp.\ \eqref{eq:offline_mmr})}\; 
\textbf{Output:} Conservative (suboptimal) set of $K$ queries $\{{\bm \iota}^{\kappa}\}_{\kappa \in \sets K}$ to ask the user\;

 \For{$\kappa \in \{1,\ldots,K\}$}{
 
    \eIf{$\kappa=1$}
    {
        Solve the MBLP reformulation of \newpv{Problem~$(\offlinerisk{1})$, (resp.\ $(\offlineregret{1})$) either directly or using the CCG algorithm}\;
        
        Let ${\bm \iota}^{\star,1}$ denote an optimal query\;
    }
    {
        Solve the MBLP reformulation of \newpv{Problem~$(\offlinerisk{\kappa})$ (resp.\ $(\offlineregret{\kappa})$)} with the added constraints that ${\bm \iota}^k={\bm \iota}^{\star,k}$ for all $k \in \{1,\ldots,\kappa-1\}$\;
        
        Let $\{{\bm \iota}^{\star,k}\}_{k=1}^\kappa$ denote an optimal query\;
    }
 }
 \textbf{Result:} Return $\{{\bm \iota}^{\star,\kappa}\}_{\kappa\in \sets K}$.
 \caption{Heuristic algorithm for \newpv{solving Problem~\eqref{eq:offline_mmu} (resp.\ \eqref{eq:offline_mmr})}.} 
 \label{alg:heuristic}
\end{algorithm}


\section{\newpv{Column-and-Constraint Generation for Min-Max Regret Active Preference Elicitation}}
\label{sec:EC_CCG_MMR}

We describe our column-and-constraint generation procedure using finite programs. Naturally, all problems solved would need to be converted to MBLPs first, using techniques similar to those employed in the proof of Theorem~\ref{thm:offline_mmr_MILP}. We omit these conversions to streamline presentation.

We define the relaxed master problem parameterized by the index set ${\sets S}' \subseteq \widetilde{\sets S}^K$
\newpv{\begin{equation}
\tag{$\mathcal{CCG}_{\text{r}}^{\text{main}}({\sets S}')$}
\begin{array}{cl}
\minimize & \quad \tau \\
\subjectto  & \quad \tau \in \reals, \; {\bm \iota} \in \sets C^K ,\; {\bm x}^{\bm s} \in \sets R, \; {\bm s} \in \widetilde{\sets S}^K \\
& \quad \!\!\left. \begin{array}{l}
{\bm \alpha}^{(\bm x', \bm s)} \in \reals_-^K, \; {\bm \beta}^{(\bm x', \bm s)} \in \reals_-^M , \; \mu^{(\bm x', \bm s)} \in \reals_+ \\
\tau  \;\; \geq \;\; {\bm b}^\top {\bm \beta}^{(\bm x', \bm s)} + \Gamma  \mu^{(\bm x', \bm s)} \\
\displaystyle \sum_{\kappa \in \sets K} {\bm s}_\kappa \; ({\bm x}^{{\bm \iota}^\kappa_1} - {\bm x}^{{\bm \iota}^\kappa_2}) \; {\bm \alpha}_\kappa^{(\bm x', \bm s)} + {\bm B}^\top {\bm \beta}^{(\bm x', \bm s)} \; =  \; {\bm x}' - {\bm x}^{\bm s} \\
{\bm \alpha}^{(\bm x', \bm s)} + \mu^{(\bm x', \bm s)} \1 \geq {\bm 0}
\end{array} \quad \quad \right\} \quad 
\begin{array}{c}
\forall {\bm s}\in {\sets S}', \\  \bm x' \in \sets R.
\end{array}
\end{array}
\label{eq:offline_mmr_ccg_rmp}
\end{equation}}
This problem only involves a subset of the decision variables and constraints of Problem~\eqref{eq:offline_mmr_MILP} (those indexed by ${\bm s} \in {\sets S}' \subseteq \widetilde{\sets S}^K$).
To identify \newpv{scenarios that,} given a solution $(\tau,{\bm \iota})$ to the relaxed \newpv{main} problem, \newpv{must be added to~\eqref{eq:offline_mmr_ccg_rmp}}, we solve a \emph{single} feasibility MBLP defined through 
\newpv{\begin{equation}
\tag{$\mathcal{CCG}^{\rm{feas}}_{\rm{r}}({\bm \iota})$}
\begin{array}{cl}
    \maximize & \quad \theta \\
    \subjectto & \quad \theta\in \reals,\; {\bm s} \in \widetilde{\sets S}^K \\
    & \quad {{\bm x}'}^{,\bm x} \in \sets R,\; {\bm u}^{\bm x} \in \sets U_0, \; {\bm \epsilon}^{\bm x} \in \sets E_\Gamma \quad \forall {\bm x} \in \sets R \\
    & \quad \theta \; \leq \; ({\bm u}^{\bm x})^\top ( {{\bm x}'}^{,\bm x} - {\bm x}) \quad \forall {\bm x} \in \sets R \\
    & \!\!\! \left. \begin{array}{l}
    \quad ({\bm u}^{\bm x})^\top ( {\bm x}^{{\bm \iota}^\kappa_1} - {\bm x}^{{\bm \iota}^\kappa_2})  + {\bm \epsilon}_\kappa^{\newpv{\bm x}} \; \geq \; M ( {\bm s}_\kappa -1) \\
    \quad ({\bm u}^{\bm x})^\top ( {\bm x}^{{\bm \iota}^\kappa_1} - {\bm x}^{{\bm \iota}^\kappa_2})  - {\bm \epsilon}_\kappa^{\newpv{\bm x}} \; \leq \; M( {\bm s}_\kappa + 1) 
    \end{array} \quad \right\} \quad \forall  \kappa \in \sets K, \; {\bm x}\in \sets R .
\end{array}
\label{eq:offline_mmr_ccg_feas}
\end{equation}}
%

%
Problems~\eqref{eq:offline_mmr_ccg_rmp} and~\eqref{eq:offline_mmr_ccg_feas} can be used in Algorithm~\ref{alg:offline_mmr_ccg} \newpv{to solve Problem~\eqref{eq:offline_mmr}. Convergence of the algorithm} is guaranteed by the following theorem.
\begin{theorem}
Algorithm~\ref{alg:offline_mmr_ccg} terminates in a final number of steps with a feasible solution to the regret averse active preference elicitation Problem~\eqref{eq:offline_mmr}. The objective value attained by this solution is within $\delta$ of the optimal objective value of the problem.
\label{thm:offline_mmr_ccg_algo_converges}
\end{theorem}

\begin{algorithm}[t!]
\SetAlgoLined
\textbf{Inputs:} Optimality tolerance $\delta$, comparison set $\sets C$, and recommendation set $\sets R$\; Initial uncertainty set $\sets U^0$ and number of queries $K$\;
\textbf{Output:} Query ${\bm \iota}^\star$ from $\sets C^K$, near optimal in Problem~\eqref{eq:offline_mmr_2} with associated objective $\theta$\;
\textbf{Initialization:} ${\bm \iota}^\star \leftarrow \emptyset$; Upper and lower bounds: ${\rm{UB}} \leftarrow +\infty$ and ${\rm{LB}} \leftarrow -\infty$\;
Initialize index set: ${\sets S}' \leftarrow \emptyset$\;
 \While{${\rm{UB}}-{\rm{LB}} > \delta$}{
  Solve the master problem~\eqref{eq:offline_mmr_ccg_rmp}. If it is solvable, let $(\tau,{\bm \iota}, \{{{\bm \alpha}}^{(\bm x', \bm s)}, {\bm \beta}^{(\bm x', \bm s)} \}_{\bm x' \in {\sets R}, {\bm s} \in {\sets S}'}, \{ {\bm x}^{\bm s} \}_{{\bm s} \in {\sets S}'} )$ be an optimal solution. If it is unbounded, set $\tau = -\infty$ and let ${\bm \iota} \in \sets C^K$ be such that~\eqref{eq:offline_mmr_ccg_rmp} is unbounded when ${\bm \iota}$ is fixed to that value\; 
  Set ${\rm{LB}} \leftarrow \tau$\;
  Solve feasibility subproblem~\eqref{eq:offline_mmr_ccg_feas}. Let $(\theta,\{ {\bm u}^{{\bm x}}, {{\bm x}'}^{,{\bm x}} \}_{{\bm x}\in \sets R}, {\bm s} )$ denote an optimal solution\;
  Set ${\rm{UB}} \leftarrow \theta$\;
    \If{$\theta > \tau$}{
        $\sets S' \leftarrow \sets S' \cup \{ {\bm s}  \}$
    }
 }
 Set ${\bm \iota}^\star \leftarrow {\bm \iota}$\;
\textbf{Result:} Collection of queries ${\bm \iota}^\star$ near-optimal in~\eqref{eq:offline_mmr_2} with objective value $\theta$.
 \caption{Column-and-Constraint Generation procedure for solving Problem~\eqref{eq:offline_mmr_2}.}
 \label{alg:offline_mmr_ccg}
\end{algorithm}


\section{\newpv{Detail of Numerical Experiments on Synthetic Instances}}
\label{sec:EC_performance}



\newpv{We performed four set of experiments. In the first (resp.\ second) set of experiments, we evaluate the performance of our offline (resp.\ online) elicitation procedure for max-min utility (\texttt{MMU}) and min-max regret (\texttt{MMR}). In the third set of experiments, we investigate the computational performance of the CCG algorithm and the symmetry breaking constraints. In the last set of experiments, we investigate the relative performance of max-min utility and min-max regret solutions. Unless explicitly stated otherwise, we solve the offline elicitation problems~\eqref{eq:offline_mmu} and \eqref{eq:offline_mmr} using the column-and-constraint generation procedures discussed in Sections~\ref{sec:offline_mmu_ccg} and \ref{sec:offline_mmr_ccg}, respectively, with tolerance~$\delta = 1 \times 10^{-3}$.  we augment all our formulations with the symmetry breaking constraints proposed in Section~\ref{sec:symmetry}. To speed-up computation further, we also employ the greedy heuristic from Section~\ref{sec:EC_speedup_strategies}. Without this conservative heuristic, the quality of the solutions obtained by our approach would be higher than those reported in these experiments.

For the first two sets of experiments, we investigate 7 settings of $(I,J,\sigma)$ for each of \texttt{MMU} and \texttt{MMR}, where $\sigma$ denotes the standard deviation of the inconsistencies. For each setting, we generate the features ${\bm x}^i$, $i\in \sets I$, uniformly at random from the~$J$-dimensional sphere of radius ten. Throughout these experiments, we assume $\sets U^0= [-1,1]^J$ and set $\Gamma= 2 \sigma \sqrt{K} \text{erf}^{-1}(2p-1)$, where $p=90\%$ is a confidence level and this definition is compatible with the elements of ${\bm \epsilon}$ being independent and normally distributed with standard deviation $\sigma$. Throughout our experiments, we use slightly larger datasets for \texttt{MMU} than for \texttt{MMR}. The reason is that we benchmark against random elicitation and evaluating the performance of a large number of random queries becomes computationally prohibitive for large datasets in the \texttt{MMR} case.

In order to facilitate an interpretation of the performance of our approach, we standardize the utilities (regrets) of the recommended items. In the case of max-min utility, we standardize worst-case utilities between zero and one. \newpv{Mathematically, the normalized worst-case utility associated with a choice of queries ${\bm \iota}$ is calculated as
$$
u_{\rm{nwc}}({\bm \iota}) \; = \; \frac{ u_{\rm{wc}}({\bm \iota}) - u_{\rm{wc}}^0 }{u_{\rm{wc}}^{\rm full} - u_{\rm{wc}}^0 },
$$
where $u_{\rm{wc}}^0$ (resp.\ $u_{\rm{wc}}^{\rm full}$) denote the worst-case utility of the item recommended if no queries are asked (resp.\ if the utility is known).} Thus, if after asking queries, recommendations are made optimally (to maximize worst case utility), then their normalized worst-case utility will be in the range $[0,1]$. On the other hand, if recommendations are made in a way that does not maximize worst-case utility (e.g., according to the analytic center of the uncertainty set as in~\cite{Toubia_2004}), then they may have a normalized worst-case utility that is negative. We can then interpret the normalized worst-case utility as a percentage relative to the full-information case (with 0\% representing no information and 100\% representing full-information). \newpv{In the case of min-max regret, the normalized worst-case regret is calculated as
$$
r_{\rm{nwc}}({\bm \iota}) \; = \; \frac{ r_{\rm{wc}}({\bm \iota}) - r_{\rm{wc}}^0 }{r_{\rm{wc}}^{\rm full} - r_{\rm{wc}}^0 },
$$
where $r_{\rm{wc}}^0$ (resp.\ $r_{\rm{wc}}^{\rm full}$) denote the worst-case regret of the item recommended if no queries are asked (resp.\ if the utility is fully known).} Thus, if after asking queries, recommendations are made to minimize worst-case regret, then their normalized worst-case regret will be in the range $[0,1]$. On the other hand, if recommendations are made in a way that does not minimize worst-case regret, then they may have a normalized worst-case regret that is greater than one. We can then interpret the normalized worst-case regret as a percentage relative to the no-information case (with 100\% representing no information and 0\% representing full-information).

All of our experiments were performed on a Linux virtual machine (Ubuntu 18.04.1 LTS), with 8GB RAM and two Intel Xeon 2.6GHz cores. All linear optimization problems were solved using Gurobi version~9.0.0. We use Mosek version~9.1.l0 to compute the analytic centers of polyhedra as required by~\texttt{POLY}, \texttt{PROB}, and \texttt{ROB}.

\paragraph{Offline Elicitation.} In our offline experiments, we compare \texttt{MMU} (resp.\ \texttt{MMR}) to random elicitation (\texttt{RAND}) for $K \in \{1,\ldots,10\}$. For any fixed~$K$, \texttt{RAND} selects~$K$ queries at random. To the best of our knowledge, no other solution approach exists in the literature for identifying~$K$ queries to ask at once. For \texttt{RAND}, we sample 50 sets of $K$ random queries. We evaluate the worst-case utility (resp.\ worst-case regret) of any given choice of queries~${\bm \iota}^\star$, which we denote by $u_{\rm{wc}}({\bm \iota}^\star)$ (resp.\ $r_{\rm{wc}}({\bm \iota}^\star)$), by fixing ${\bm \iota}={\bm \iota}^\star$ in Problem~\eqref{eq:offline_mmu_MILP} (resp.\ \eqref{eq:offline_mmr_MILP}). To speed-up this objective function evaluation step, we employ a variant of the column-and-constraint generation procedure presented in this paper. We omit the details in the interest of space. The optimality and scalability results for \texttt{MMU} (resp.\ \texttt{MMR}) are provided in Figure~\ref{fig:offline_mmu} (resp.\ \ref{fig:offline_mmr}).

\begin{figure}[ht!]
\centering
    \includegraphics[width=0.8\textwidth]{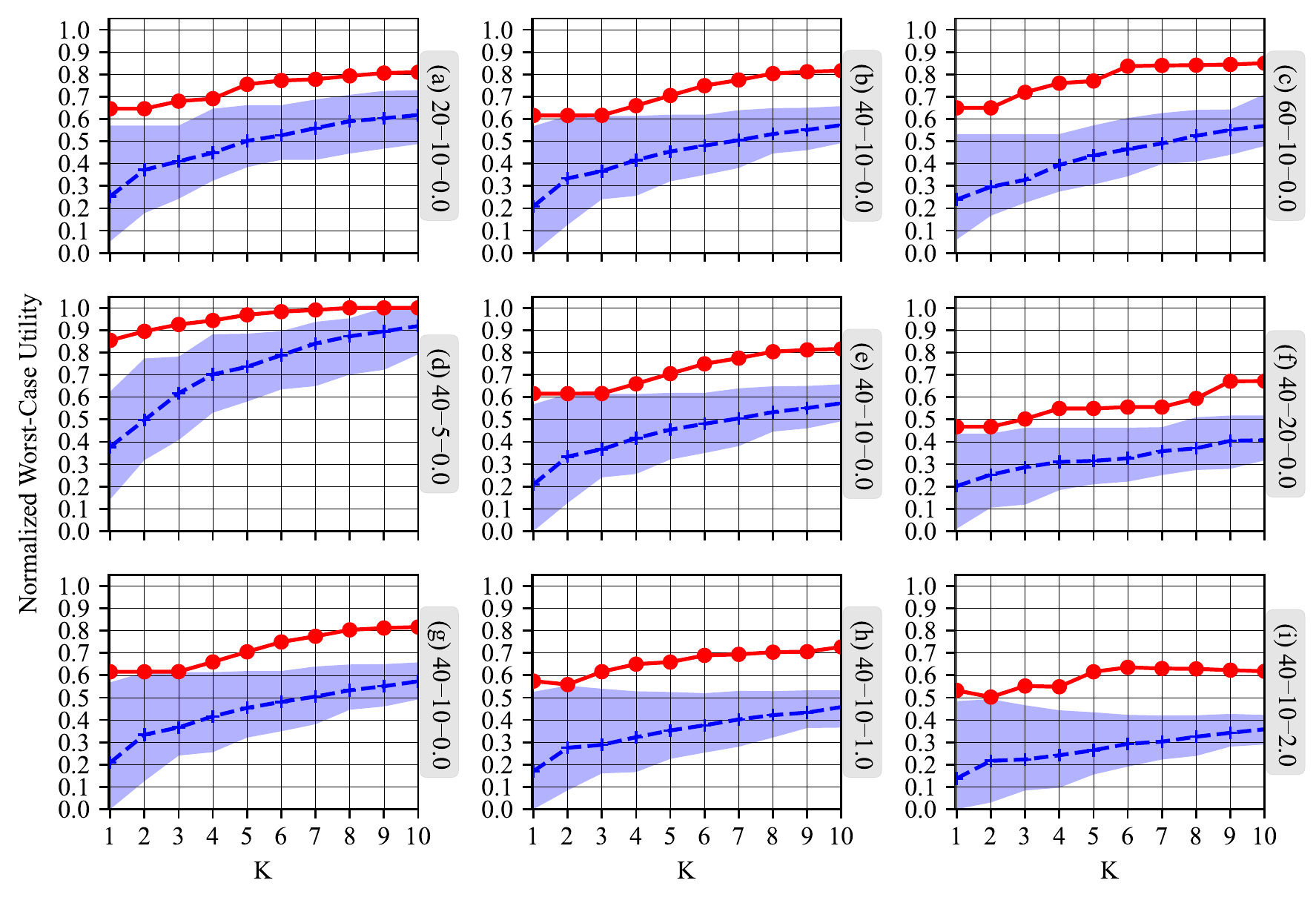}
    \includegraphics[width=0.8\textwidth]{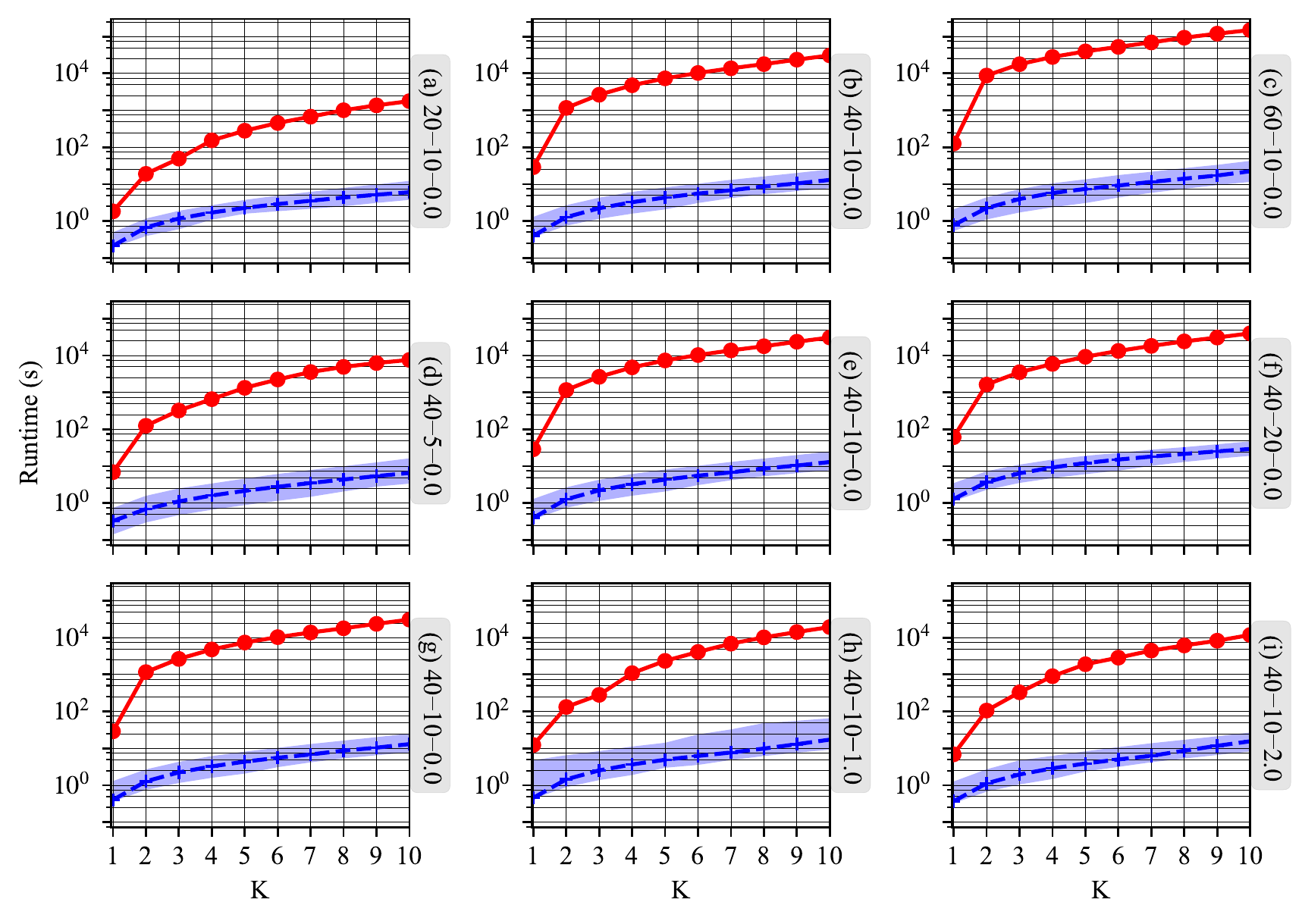}
    \caption{\newpv{Optimality (top) and scalability (bottom)} results for the offline \newpv{max-min utility} preference elicitation problem~\eqref{eq:offline_mmu} on synthetic data. Each label on the right of each facet corresponds to the characteristics of the instance solved ($I-J-\sigma$). Approach \texttt{MMU} is shown with red dots. The median performance of \texttt{RAND} across 50 sets of $K$ random queries is shown with blue crosses. \newpv{In particular the run time for \texttt{RAND} is the time it takes to evaluate the performance of one query.} The blue shaded region shows the range of performance of \texttt{RAND} across the 50 sets of queries.}
\label{fig:offline_mmu}
\end{figure}


\begin{figure}[h!]
\centering
    \includegraphics[width=0.8\textwidth]{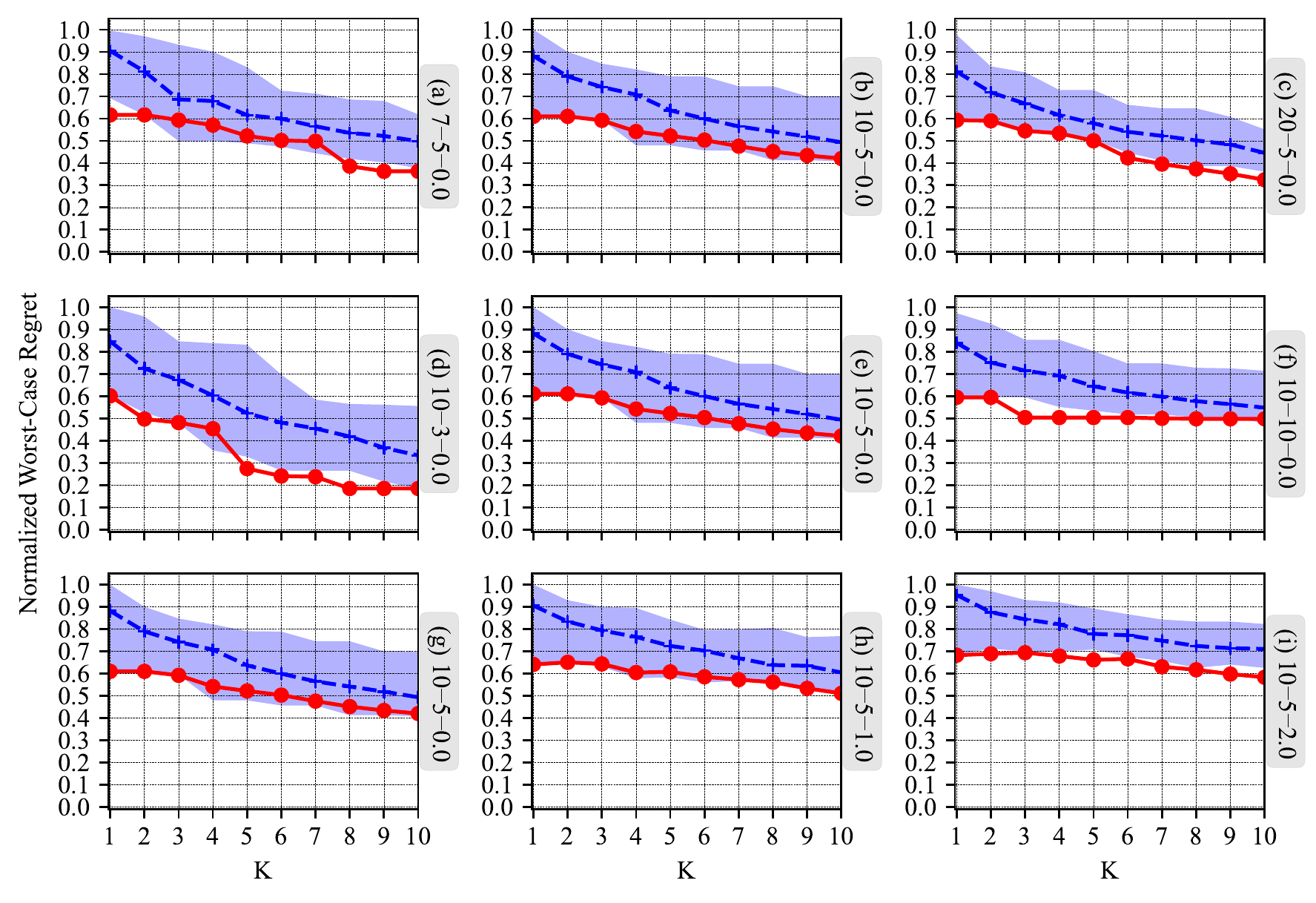}
    \includegraphics[width=0.8\textwidth]{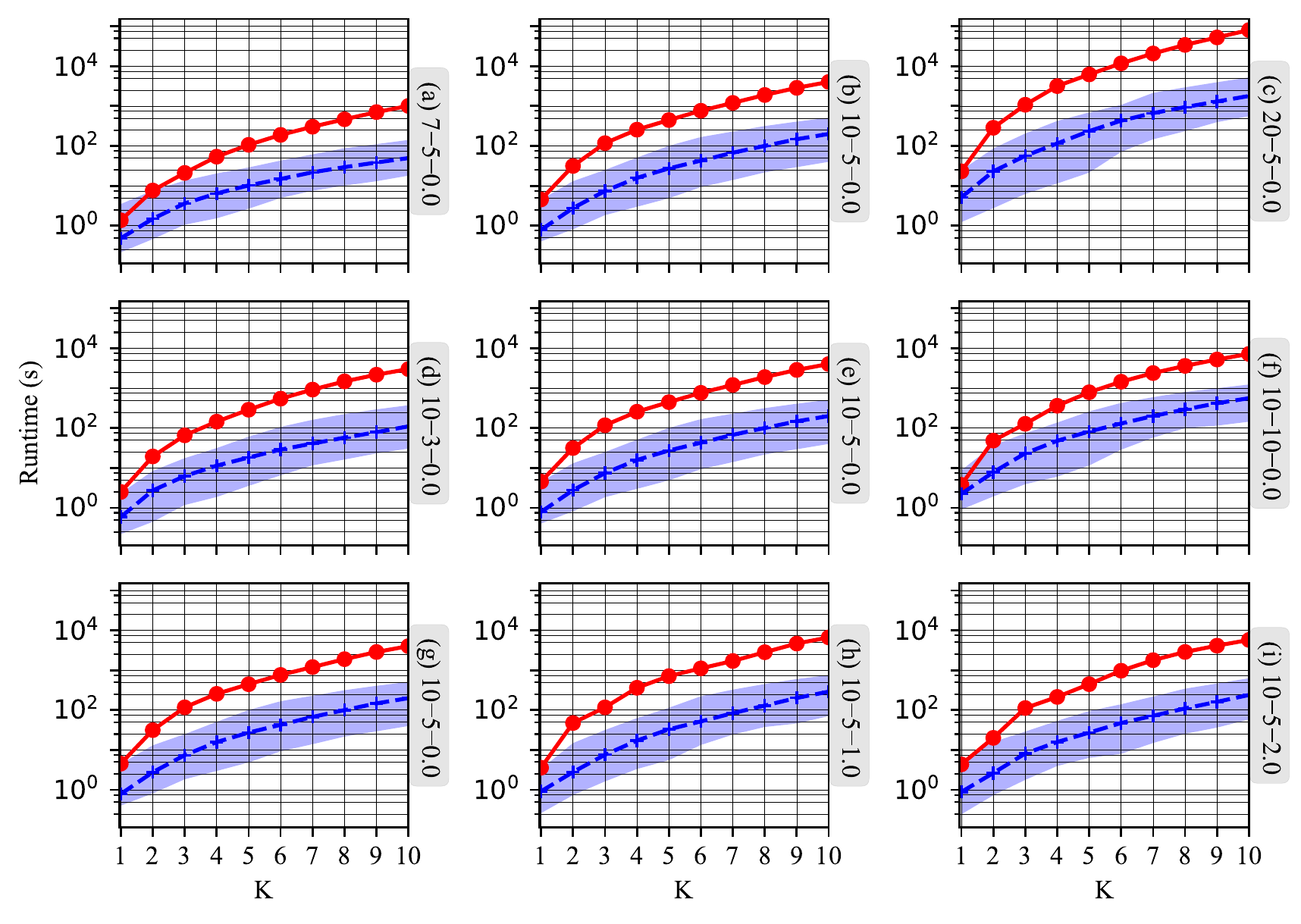}
    \caption{\newpv{Optimality (top) and scalability (bottom)} results for the offline \newpv{min-max regret} preference elicitation problem~\eqref{eq:offline_mmr} on synthetic data. The facet labels, graphs, shapes, lines, and colors have similar interpretation as in Figure~\ref{fig:offline_mmu}.}
    \label{fig:offline_mmr}
\end{figure}

}

\newpv{

\paragraph{Online Elicitation.} In the online setting, we compare our conservative solution approaches (\texttt{MMU} and \texttt{MMR}) to random elicitation, to the polyhedral method of~\cite{Toubia_2004} (\texttt{POLY}), to the probabilistic polyhedral method of~\cite{Toubia_2007} (\texttt{PROB}), to the robust approach of~\cite{OHair_LearningPreferences} (\texttt{ROB}), and to the ellipsoidal method of~\cite{Vielma2019} (\texttt{ELL}). For the online setting, we evaluate the performance of all approaches through simulation using 50 pairs (${\bm u},{\bm \epsilon}$). The utility vectors are drawn uniformly at random from the $J$-dimensional sphere of radius one. The inconsistencies are drawn from the normal distribution with standard deviation $\sigma$. For each method, we generate queries according to its elicitation approach, generate answers to the queries using ${\bm u}$ and ${\bm \epsilon}$, and subsequently make a recommendation according to the recommendation approach used by that method. We record the worst-case utility, worst-case regret, and rank of the recommended item. These are summarized in Figure~\ref{fig:online_mmu} (resp.\ \ref{fig:online_mmr}) for \texttt{MMU} (resp.\ \texttt{MMR}).


\begin{figure}[h!]
\centering
    \includegraphics[width=0.8\textwidth]{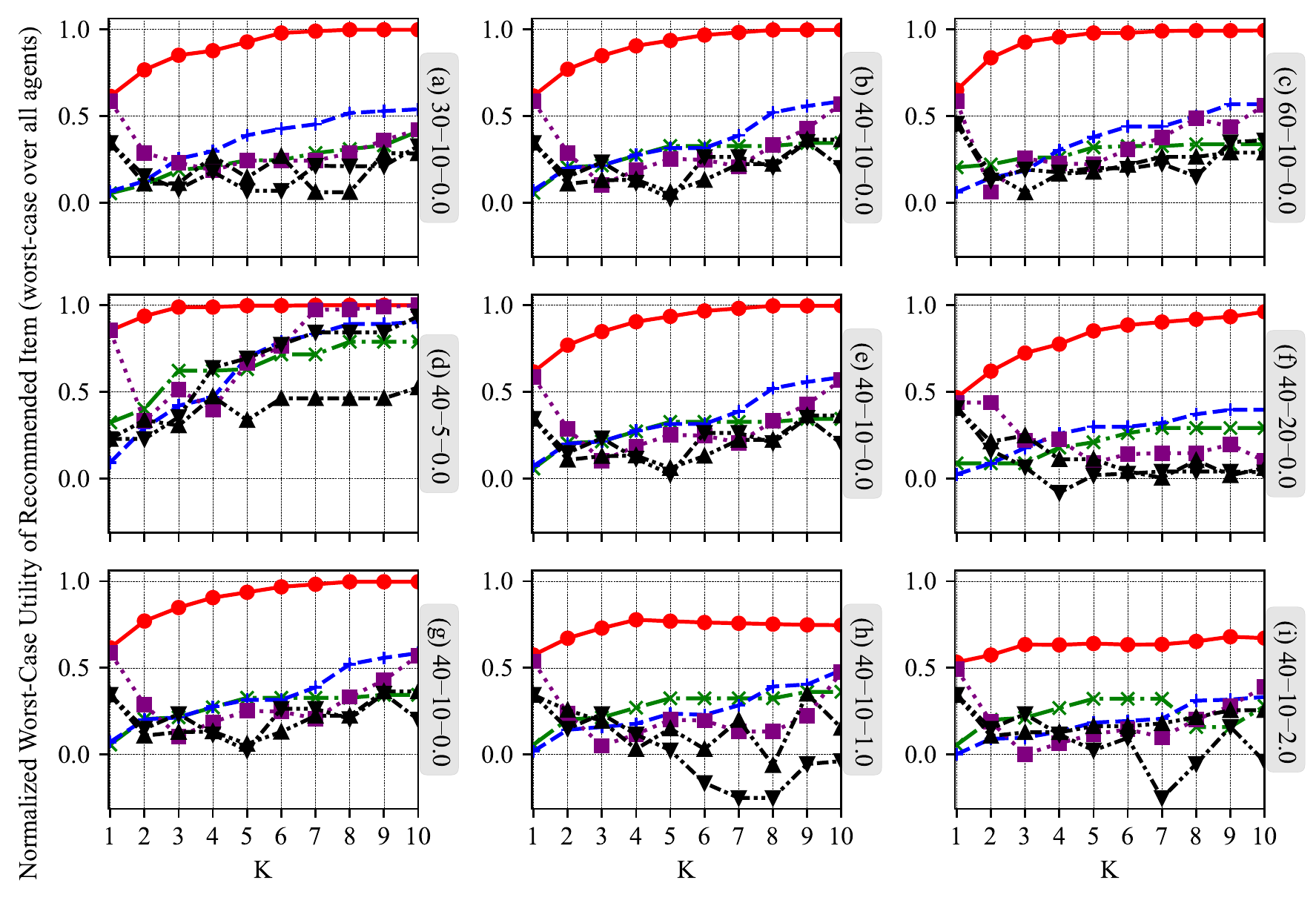}
    \includegraphics[width=0.8\textwidth]{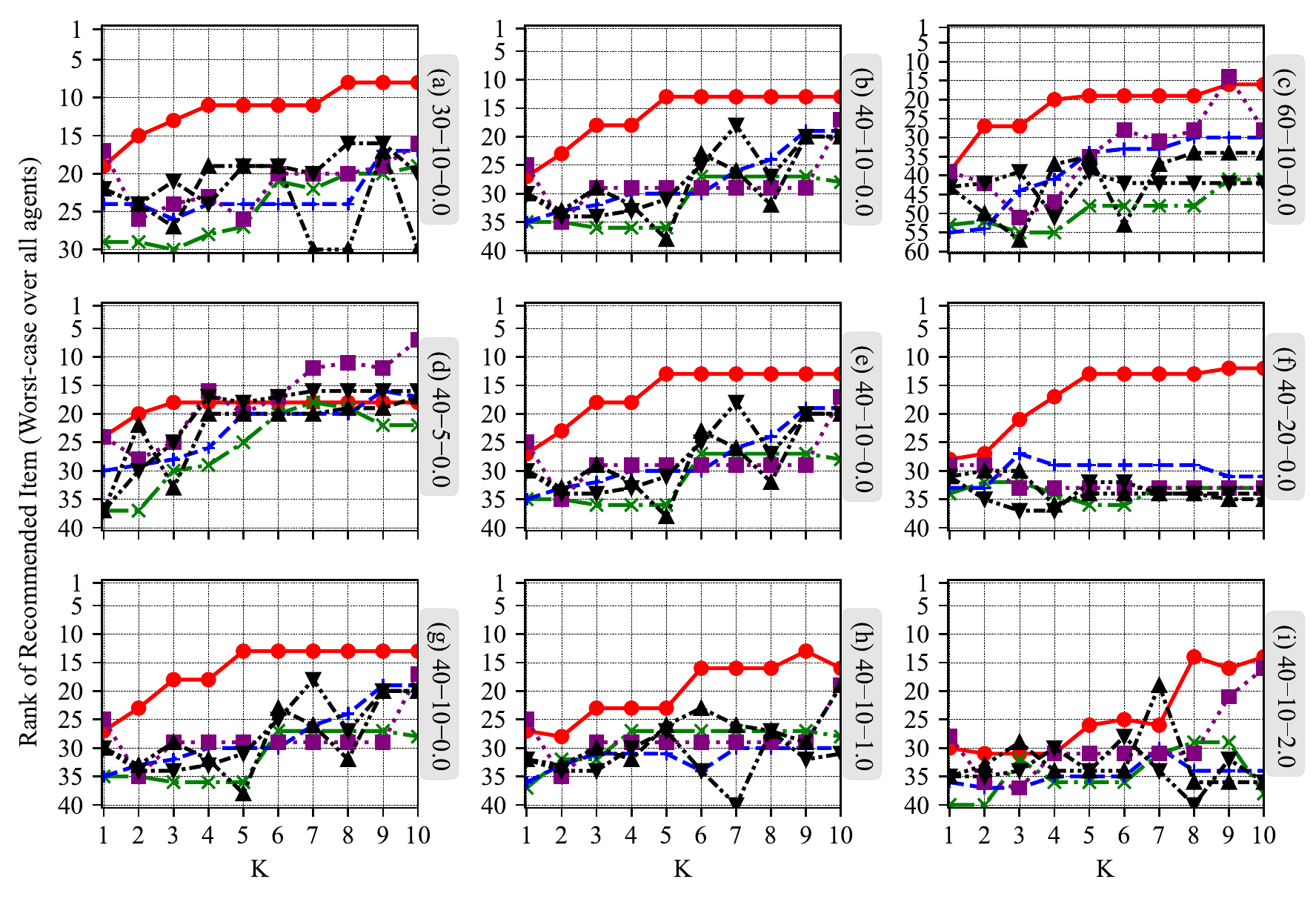}
    \caption{\newpv{Worst-case utility (top) and worst-case rank (bottom)} for the online \newpv{max-min utility} preference elicitation problem on synthetic data. The facet labels have the same interpretation as in Figure~\ref{fig:offline_mmu}. \newpv{The median performance of \texttt{MMU} (resp.\ \texttt{RAND}, \texttt{POLY}, \texttt{PROB}, \texttt{ROB}, \texttt{ELL}) across 50 random utility vectors ${\bm u}$ and inconsistencies ${\bm \epsilon}$ is shown with red dots {(resp.\ blue plus, black triangles, black inverted triangles, green crosses, and purple squares).}}} 
    \label{fig:online_mmu}
\end{figure}


\begin{figure}[h!]
\centering
    \includegraphics[width=0.8\textwidth]{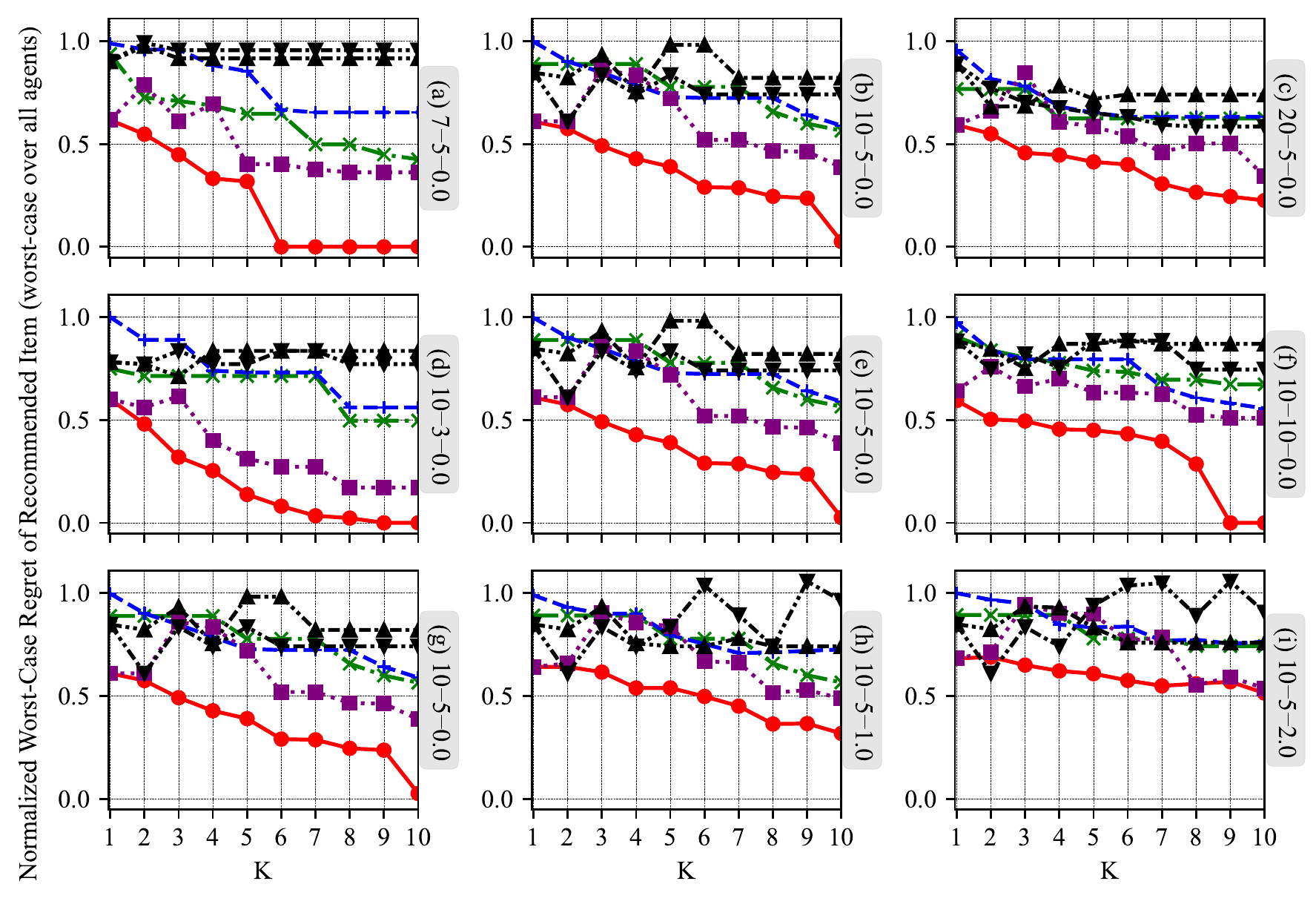}
    \includegraphics[width=0.8\textwidth]{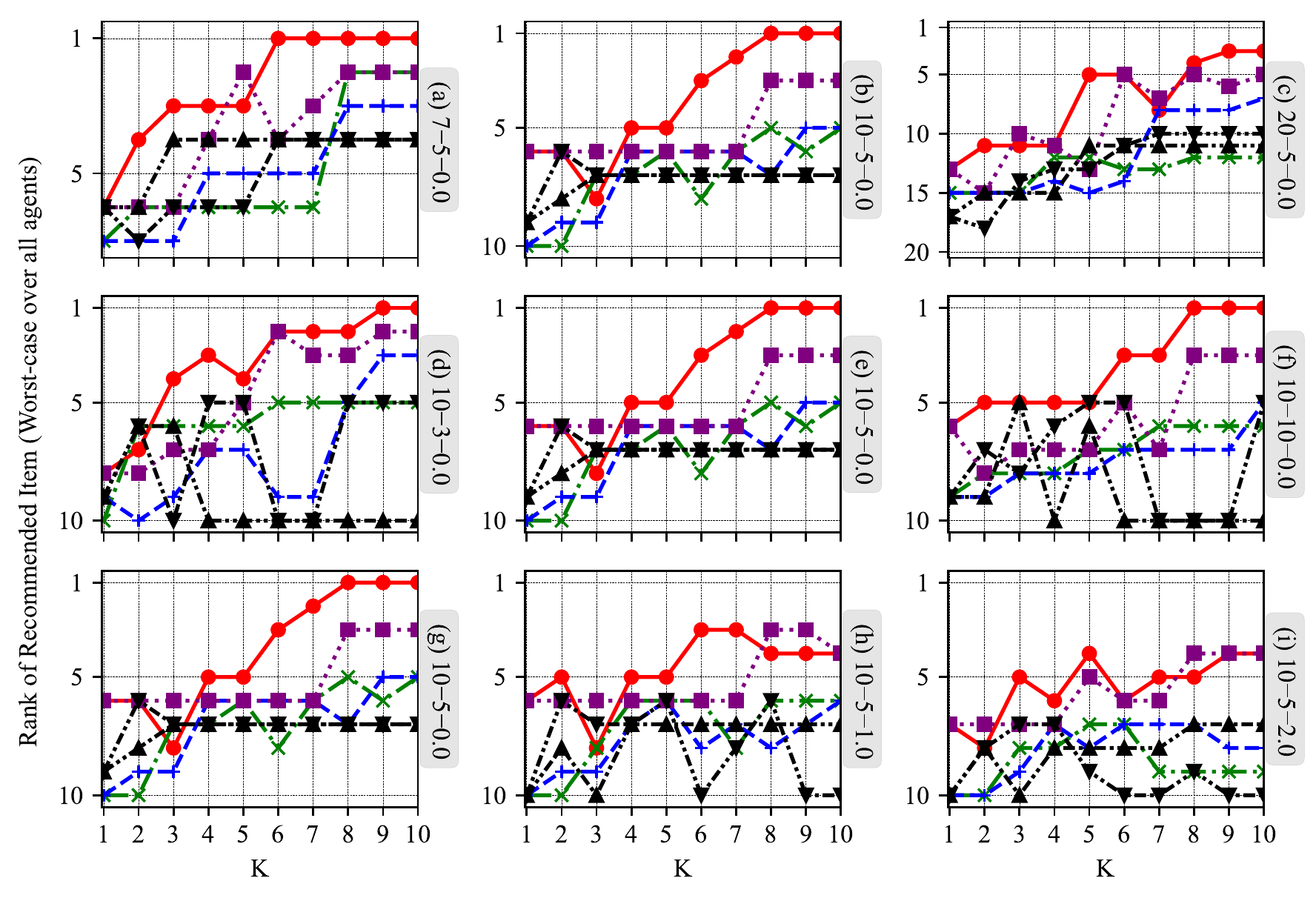}
    \caption{\newpv{Worst-case utility (top) and worst-case rank (bottom)} for the online \newpv{min-max-regret} preference elicitation problem on synthetic data. \newpv{The facet labels and colors have the same interpretation as in Figure~\ref{fig:online_mmu} with the exception that red dots correspond to \texttt{MMR}}.}
    \label{fig:online_mmr}
\end{figure}

}

\newpv{\paragraph{Comparison Between Max-Min Utility \& Min-Max Regret.} In the first set of experiments, we observed that the max-min utility problem is more scalable than its min-max regret counterpart. In our third set of experiments, we investigate whether there are benefits in employing the min-max regret solution relative the max-min utility solution. For this reason, we study the worst-case utility and worst-case regret of solutions to Problems~\eqref{eq:offline_mmu} and~\eqref{eq:offline_mmr}, respectively, on a synthetic dataset with $I = 10$ items and $J = 10$ features ($\Gamma = 0$). The results are summarized in Table~\ref{tab:maxminutility_vs_minmaxregret}.

\begin{table}[t]
\centering 
\caption{Comparison between max-min utility and min-max regret based queries on a synthetic dataset with $I = 10$ items and $J = 10$ features ($\Gamma = 0$). The max-min utility queries ${\bm \iota}_{\rm u}^\star$ and min-max regret queries ${\bm \iota}_{\rm r}^\star$ are computed by solving the MBLP reformulations of Problems~\eqref{eq:offline_mmu} and~\eqref{eq:offline_mmr}, respectively, using the column-and-constraint generation algorithm. The decrease in normalized worst-case utility refers to the drop in normalized worst-case utility experienced by employing the min-max regret rather than max-min utility solution, computed as $(u_{\rm{nwc}}({\bm \iota}^\star_{\rm u})-u_{\rm{nwc}}({\bm \iota}^\star_{\rm r}))$. Similarly, the decrease in normalized worst-case regret refers to the drop in normalized worst-case regret experienced by employing the min-max regret rather than max-min utility solution, computed as $(r_{\rm{nwc}}({\bm \iota}^\star_{\rm u})-r_{\rm{nwc}}({\bm \iota}^\star_{\rm r}))$. All approaches were given a one hour time limit.}
\begin{small}
    \renewcommand{\arraystretch}{1.1}
    \begin{tabular}{@{}c*{8}{c}@{}} 
     \toprule
 &	\multicolumn{2}{c}{Max-Min Utility Solution} &	\multicolumn{2}{c}{Min-Max Regret Solution} & \multirow{4}{*}{ \begin{tabular}[c]{@{}c@{}} Decrease in \\ Normalized \\ Worst-Case \\ Utility (p.p.) \end{tabular}} &	\multirow{4}{*}{\begin{tabular}[c]{@{}c@{}}  Decrease in \\ Normalized \\ Worst-Case \\ Regret (p.p.) \end{tabular}} \\
\cline{2-5} 
\multirow{3}{*}{$K$} & \multirow{3}{*}{\begin{tabular}[c]{@{}c@{}}  Normalized \\ Worst-Case \\ Utility \end{tabular}} &		\multirow{3}{*}{\begin{tabular}[c]{@{}c@{}}  Normalized \\  Worst-Case \\ Regret \end{tabular}} &	\multirow{3}{*}{\begin{tabular}[c]{@{}c@{}}  Normalized \\  Worst-Case \\ Utility \end{tabular}} &		\multirow{3}{*}{\begin{tabular}[c]{@{}c@{}}  Normalized \\  Worst-Case \\ Regret \end{tabular}}	&  & \\ \\ \\ \bottomrule
2& {0.780} &  {0.642} & {0.689} & {0.594} & 9.1 & {4.8} \\
4& {0.783} & {0.620} & {0.691} & {0.532} & {9.2} & {8.8} \\
6& {0.833} & {0.566} & {0.712} & {0.522} & {12.1} & {4.4} \\
8& {0.879} & {0.550} & {0.816} & {0.499} &{6.3} & {5.1} \\
10& {0.912} & {0.510} & {0.844} & {0.495} & {6.8} & {1.5} \\
\bottomrule
    \end{tabular}
    \end{small}
    \label{tab:maxminutility_vs_minmaxregret}
\end{table}


\begin{table}
\caption{Evaluation results of symmetry breaking constraints, CCG algorithm, and CCG-based greedy heuristic approach on three synthetic datasets with $I \in \{5,10,20\}$ items and $J=10$ features ($\Gamma=0$); A dash indicates that the optimal solution was not found within the allotted 3 hour time limit.}
\renewcommand{\arraystretch}{1}
\begin{footnotesize}
\begin{center}
\begin{tabular}{@{}|c|c|c|c|c|c|c|c|c|c|@{}}
\toprule
\multirow{2}{*}{$I$} & \multirow{2}{*}{$K$} & \multicolumn{2}{c|}{MBLP}  & \multicolumn{2}{c|}{MBLP + Symm. Break.} & \multicolumn{2}{c|}{CCG + Symm. Break.} & \multicolumn{2}{c|}{Heuristic}  \\ \cline{3-10} 
 & & \begin{tabular}[c]{@{}c@{}} Normalized \\ Objective \\ Value \end{tabular} & \begin{tabular}[c]{@{}c@{}}Solve \\ Time\\ (sec)\end{tabular} & \begin{tabular}[c]{@{}c@{}} Normalized \\  Objective \\ Value\end{tabular} & \begin{tabular}[c]{@{}c@{}}Solve \\ Time\\ (sec)\end{tabular} & \begin{tabular}[c]{@{}c@{}} Normalized \\ Objective \\ Value\end{tabular} & \begin{tabular}[c]{@{}c@{}}Solve \\ Time\\ (sec)\end{tabular} & \begin{tabular}[c]{@{}c@{}} Normalized \\ Objective \\ Value\end{tabular} & \begin{tabular}[c]{@{}c@{}}Solve \\ Time\\ (sec)\end{tabular} \\ \bottomrule
5 & 2 & {0.800} & 0.34 &  {0.800} & 0.25 & {0.800} & 0.77 & {0.800} & 0.51 \\
5 & 4 & {0.932} & 622.64 & {0.932} & 2.45 & {0.932} & 5.48 & {0.932} & 1.40 \\
5 & 6 & {0.978} & 6094.79 & {0.978} & 76.21 & {0.978} & 10.80 & {0.932} & 2.25 \\
5 & 8 & -- & -- & {0.978} & 1147.51 & {0.978} & 40.91 & {0.972} & 3.53 \\
5 & 10 & -- & -- & {1.000} & 5082.32 & {1.000} & 0.86 & {1.000} & 4.24 \\
\midrule
10 & 2 & {0.780} & 2.66 & {0.780} & 2.61 & {0.780} & 5.54 & {0.780} & 2.13 \\
10 & 4 & -- & -- & {0.879} & 1141.58 & {0.879} & 6054.11 & {0.783} & 5.67 \\
10 & 6 & --  & -- & -- & -- & -- & -- & {0.812}  & 12.53 \\
10 & 8 & -- & -- & -- & -- & -- & -- & {0.879}  & 21.44 \\
10 & 10 & -- & -- & -- & -- & -- & -- & {0.879} & 31.79 \\
\midrule
20 & 2 & {0.646} & 334.33 & {0.646} & 202.75 & {0.646} & 338.17 & {0.646} & {5.21} \\
20 & 4 & -- & -- & -- & -- & -- & -- & {0.658} & {23.24} \\
20 & 6 & -- & -- & -- & -- & -- & -- & {0.765} & {51.16} \\
20 & 8 & -- & -- & -- & -- & -- & -- & {0.786} & {84.35} \\
20 & 10 & -- & -- & -- & -- & -- & -- & {0.812} & {162.32} \\
\bottomrule
\end{tabular}
\end{center}
\end{footnotesize}
\label{tab:speedup_strategies}
\end{table}


\paragraph{Evaluation of Column-and-Constraint Generation, Symmetry Breaking, \& Greedy Heuristic.} For our fourth set of experiments, we evaluate the performance of our speed-up strategies. To this end, we solve the offline risk averse active preference elicitation problem~\eqref{eq:offline_mmu} based on three synthetic datasets using four different approaches: the MBLP problem~\eqref{eq:offline_mmu_MILP}, the MBLP problem~\eqref{eq:offline_mmu_MILP} augmented with symmetry breaking constraints (see Section~\ref{sec:symmetry}), the column-and-constraint generation (CCG) approach from Algorithm~\ref{alg:offline_mmu_ccg} augmented with symmetry breaking constraints, and the CCG-based conservative heuristic approach, see Section~\ref{sec:EC_greedy}. The results are provided in Table~\ref{tab:speedup_strategies}. From the table, it can be seen that the symmetry breaking constraints and the CCG algorithm augmented with symmetry breaking significantly speed-up computation and enable the solution of problems that could not be solved with the MBLP alone. From the table, it also becomes apparent that the CCG-based conservative heuristic returns near-optimal solutions in a fraction of the time, thus scaling to far larger instances.

}



\phantom{sf}

\newpage

\section{\newpv{Detail of Numerical Experiments on Real Data}}
\label{sec:EC_performance_real}

In Tables~\ref{tab:offline_sigma_comparison_detail} and~\ref{tab:online_sigma_comparison_detail}, we provide the full detail of our results on evaluating robustness of our model to the choice of $\Gamma$ in the real dataset.


\begin{table}[ht]
\caption{\newpv{Performance of the offline \texttt{MMR} approach on the real data when the standard deviation~$\sigma$ of the inconsistencies used to select $\Gamma$ in the optimization is different than the true standard deviation $\sigma^\star$ of the inconsistencies.}}
\centering
\renewcommand{\arraystretch}{1}
\begin{footnotesize}
\begin{tabular}{|C{1.2cm}|C{0.5cm}||C{1.4cm}|C{1.4cm}||C{1.4cm}|C{1.4cm}|C{1.4cm}||C{1.4cm}|C{1.4cm}|C{1.4cm}|}
  \hline
\multirow{2}{*}{Metric} & \multirow{2}{*}{K} & \multicolumn{2}{c||}{$\sigma=0.0\%$} & \multicolumn{3}{c|}{$\sigma=2.5\%$} & \multicolumn{3}{c|}{$\sigma=5.0\%$} \\
\cline{3-10}
& & $\sigma^\star=0.0\%$ & $\sigma^\star=0.5\%$ & $\sigma^\star=2.0\%$ & $\sigma^\star=2.5\%$ & $\sigma^\star=3.0\%$ & $\sigma^\star=4.5\%$ & $\sigma^\star=5.0\%$ & $\sigma^\star=5.5\%$ \\
  \hline
\multirow{10}{*}{\rotatebox[origin=c]{90}{Worst-Case Regret}} & 1 & 0.804 & 0.804 & 0.898 & 0.898 & 0.898 & 0.937 & 0.937 & 0.937 \\
 & 2 & 0.716 & 0.716 & 0.810 & 0.810 & 0.810 & 0.934 & 0.934 & 0.934 \\
 & 3 & 0.671 & 0.671 & 0.780 & 0.780 & 0.780 & 0.934 & 0.934 & 0.934 \\
 & 4 & 0.667 & 0.642 & 0.757 & 0.709 & 0.757 & 0.929 & 0.929 & 0.934 \\
 & 5 & 0.658 & 0.641 & 0.721 & 0.630 & 0.712 & 0.925 & 0.816 & 0.926 \\
 & 6 & 0.640 & 0.632 & 0.699 & 0.608 & 0.686 & 0.814 & 0.804 & 0.804 \\
 & 7 & 0.597 & 0.599 & 0.689 & 0.607 & 0.670 & 0.804 & 0.797 & 0.779 \\
 & 8 & 0.547 & 0.551 & 0.670 & 0.581 & 0.651 & 0.787 & 0.769 & 0.746 \\
 & 9 & 0.539 & 0.538 & 0.670 & 0.581 & 0.635 & 0.773 & 0.727 & 0.686 \\
 & 10 & 0.538 & 0.533 & 0.652 & 0.578 & 0.612 & 0.758 & 0.689 & 0.667 \\

  \hline
\multirow{10}{*}{\rotatebox[origin=c]{90}{Worst-Case True Regret}} & 1 & 0.305 & 0.305 & 0.367 & 0.402 & 0.367 & 0.359 & 0.359 & 0.323 \\
 & 2 & 0.305 & 0.305 & 0.305 & 0.305 & 0.305 & 0.360 & 0.321 & 0.311 \\
 & 3 & 0.305 & 0.305 & 0.340 & 0.310 & 0.340 & 0.360 & 0.305 & 0.334 \\
 & 4 & 0.305 & 0.305 & 0.391 & 0.294 & 0.340 & 0.360 & 0.360 & 0.323 \\
 & 5 & 0.305 & 0.305 & 0.391 & 0.244 & 0.340 & 0.347 & 0.351 & 0.315 \\
 & 6 & 0.305 & 0.351 & 0.341 & 0.340 & 0.334 & 0.347 & 0.351 & 0.318 \\
 & 7 & 0.341 & 0.351 & 0.341 & 0.319 & 0.340 & 0.351 & 0.351 & 0.318 \\
 & 8 & 0.341 & 0.217 & 0.329 & 0.301 & 0.334 & 0.358 & 0.351 & 0.315 \\
 & 9 & 0.341 & 0.217 & 0.329 & 0.301 & 0.334 & 0.365 & 0.351 & 0.223 \\
 & 10 & 0.341 & 0.224 & 0.329 & 0.340 & 0.334 & 0.351 & 0.351 & 0.264 \\

  \hline
\multirow{10}{*}{\rotatebox[origin=c]{90}{Worst-Case True Rank}} & 1 & 11 & 11 & 21 & 23 & 21 & 20 & 18 & 18 \\
 & 2 & 11 & 11 & 11 & 11 & 11 & 15 & 15 & 15 \\
 & 3 & 11 & 11 & 18 & 18 & 19 & 14 & 16 & 14 \\
 & 4 & 11 & 11 & 20 & 16 & 13 & 16 & 20 & 20 \\
 & 5 & 11 & 9 & 20 & 8 & 11 & 17 & 18 & 17 \\
 & 6 & 11 & 17 & 19 & 12 & 15 & 15 & 17 & 20 \\
 & 7 & 13 & 13 & 19 & 12 & 15 & 15 & 20 & 20 \\
 & 8 & 13 & 8 & 15 & 12 & 22 & 20 & 17 & 17 \\
 & 9 & 13 & 8 & 15 & 12 & 15 & 15 & 17 & 19 \\
 & 10 & 12 & 8 & 14 & 11 & 12 & 15 & 15 & 16 \\
  \hline 
\end{tabular}
\end{footnotesize}
\label{tab:offline_sigma_comparison_detail}
\end{table}


\begin{table}[ht]
\caption{\newpv{Performance of the online \texttt{MMR} approach on the real data when the standard deviation~$\sigma$ of the inconsistencies used to select $\Gamma$ in the optimization is different than the true standard deviation $\sigma^\star$ of the inconsistencies.}}
\centering
\renewcommand{\arraystretch}{1}
\begin{footnotesize}
\begin{tabular}{|C{1.2cm}|C{0.5cm}||C{1.4cm}|C{1.4cm}||C{1.4cm}|C{1.4cm}|C{1.4cm}||C{1.4cm}|C{1.4cm}|C{1.4cm}|}
  \hline
\multirow{2}{*}{Metric} & \multirow{2}{*}{K} & \multicolumn{2}{c||}{$\sigma=0.0\%$} & \multicolumn{3}{c|}{$\sigma=2.5\%$} & \multicolumn{3}{c|}{$\sigma=5.0\%$} \\
\cline{3-10}
& & $\sigma^\star=0.0\%$ & $\sigma^\star=0.5\%$ & $\sigma^\star=2.0\%$ & $\sigma^\star=2.5\%$ & $\sigma^\star=3.0\%$ & $\sigma^\star=4.5\%$ & $\sigma^\star=5.0\%$ & $\sigma^\star=5.5\%$ \\
  \hline
\multirow{10}{*}{\rotatebox[origin=c]{90}{Worst-Case Regret}} & 1 &0.804 & 0.804 & 0.816 & 0.816 & 0.816 & 0.816 & 0.816 & 0.853 \\
    & 2 & 0.658 & 0.658 & 0.558 & 0.558 & 0.558 & 0.635 & 0.635 & 0.635 \\
    & 3 &0.537 & 0.537 & 0.547 & 0.547 & 0.547 & 0.607 & 0.607 & 0.607 \\
    & 4 &0.487 & 0.487 & 0.509 & 0.509 & 0.509 & 0.616 & 0.616 & 0.616 \\
    & 5 &0.414 & 0.414 & 0.487 & 0.487 & 0.487 & 0.635 & 0.635 & 0.635 \\
    & 6 &0.383 & 0.383 & 0.496 & 0.496 & 0.496 & 0.652 & 0.652 & 0.652 \\
    & 7 &0.264 & 0.264 & 0.497 & 0.497 & 0.497 & 0.661 & 0.661 & 0.661 \\
    & 8 &0.238 & 0.238 & 0.501 & 0.501 & 0.501 & 0.622 & 0.622 & 0.622 \\
    & 9 &0.231 & 0.231 & 0.483 & 0.483 & 0.483 & 0.558 & 0.558 & 0.558 \\
    & 10 &0.229 & 0.229 & 0.489 & 0.489 & 0.489 & 0.566 & 0.566 & 0.566 \\
  \hline
\multirow{10}{*}{\rotatebox[origin=c]{90}{Worst-Case True Regret}} & 1 & 0.268 & 0.268 & 0.268 & 0.268 & 0.268 & 0.268 & 0.268 & 0.245 \\
    & 2 & 0.219 & 0.219 & 0.135 & 0.135 & 0.135 & 0.135 & 0.135 & 0.121 \\
    & 3 & 0.159 & 0.159 & 0.135 & 0.135 & 0.135 & 0.135 & 0.184 & 0.165 \\
    & 4 & 0.135 & 0.135 & 0.135 & 0.135 & 0.135 & 0.135 & 0.159 & 0.143 \\
    & 5 & 0.135 & 0.135 & 0.135 & 0.135 & 0.135 & 0.135 & 0.135 & 0.121 \\
    & 6 & 0.135 & 0.135 & 0.135 & 0.135 & 0.135 & 0.176 & 0.176 & 0.158 \\
    & 7 & 0 & 0.008 & 0.107 & 0.107 & 0.107 & 0.135 & 0.135 & 0.121 \\
    & 8 & 0 & 0.008 & 0.107 & 0.107 & 0.107 & 0.135 & 0.135 & 0.121 \\
    & 9 & 0 & 0.008 & 0.107 & 0.107 & 0.107 & 0.135 & 0.135 & 0.121 \\
    & 10 & 0 & 0.008 & 0.107 & 0.107 & 0.107 & 0.135 & 0.135 & 0.121 \\
  \hline
\multirow{10}{*}{\rotatebox[origin=c]{90}{Worst-Case True Rank}} & 1 & 6 & 6 & 7 & 7 & 7 & 7 & 7 & 15 \\
    & 2 & 4 & 4 & 5 & 5 & 5 & 5 & 5 & 5 \\
    & 3 & 3 & 3 & 5 & 5 & 5 & 5 & 6 & 6 \\
    & 4 & 3 & 3 & 4 & 4 & 4 & 5 & 5 & 5 \\
    & 5 & 4 & 6 & 4 & 4 & 4 & 5 & 5 & 5 \\
    & 6 & 2 & 3 & 4 & 4 & 4 & 5 & 5 & 5 \\
    & 7 & 1 & 3 & 4 & 4 & 4 & 3 & 3 & 3 \\
    & 8 & 1 & 3 & 4 & 6 & 6 & 5 & 5 & 5 \\
    & 9 & 1 & 3 & 4 & 6 & 6 & 7 & 7 & 7 \\
    & 10 & 1 & 3 & 4 & 4 & 4 & 7 & 7 & 7 \\
  \hline 
\end{tabular}
\end{footnotesize}
\label{tab:online_sigma_comparison_detail}
\end{table}

\newpage

\phantom{sf}

\newpage


\section{Proofs of Statements in Section~\ref{sec:maxminutility}}
\label{sec:EC_maxminutility}


The following proposition is needed in the proof of Lemma~\ref{lem:offline_mmu_simplification}. It shows that we can eliminate the dependence of $\newpv{\sets S_\Gamma}({\bm \iota})$ on ${\bm \iota}$ in Problem~\eqref{eq:offline_mmu}.
\begin{proposition}
Problem~\eqref{eq:offline_mmu} is equivalent to
\begin{equation}
\maximize_{{\bm \iota} \in \sets C^K} \quad \min_{{\bm s}\in \sets S^K} \quad \max_{{\bm x} \in \sets R} \quad \min_{ {\bm u} \in  \newpv{\sets U_\Gamma}({\bm \iota} ,{\bm s}) } \quad  {\bm u}^\top {\bm x}
\label{eq:offline_mmu_3}
\end{equation}
in the sense that the two problems have the same optimal objective value and the same sets of optimal solutions.
\label{prop:offline_mmu_simplification_1}
\end{proposition}


\proof{Proof}
Fix ${\bm \iota} \in \sets C^K$. Fix ${\bm s} \in \sets S^K \backslash \{ \newpv{ \sets S_\Gamma({\bm \iota}) } \}$. Then, \newpv{$\sets U_\Gamma({\bm \iota},{\bm s}) = \emptyset$}. Since the choice of ${\bm s} \in \sets S^K \backslash \{ \newpv{ \sets S_\Gamma({\bm \iota})} \}$ was arbitrary, it follows from $\sets R \neq \emptyset$ that
$$
\min\limits_{{\bm s} \in \sets S^K \backslash \{ \newpv{\sets S_\Gamma({\bm \iota}) } \} } \quad \max\limits_{{\bm x} \in \sets R} \quad \min\limits_{ {\bm u} \in \newpv{\sets U_\Gamma}({\bm \iota} ,{\bm s})} \quad  {\bm u}^\top {\bm x} \;\; = \;\; + \infty.
$$
Therefore,
$$
\min\limits_{{\bm s} \in \sets S^K \backslash \{ \newpv{\sets S_\Gamma({\bm \iota}) } \} } \quad \max\limits_{{\bm x} \in \sets R} \quad \min\limits_{ {\bm u} \in \newpv{\sets U_\Gamma}({\bm \iota} ,{\bm s})} \quad  {\bm u}^\top {\bm x}
\quad \geq \quad
\min\limits_{{\bm s} \in \newpv{\sets S_\Gamma}({\bm \iota}) } \quad \max\limits_{{\bm x} \in \sets R} \quad \min\limits_{ {\bm u} \in \newpv{\sets U_\Gamma}({\bm \iota} ,{\bm s})} \quad {\bm u}^\top {\bm x}.
$$
This implies that
$$
\min\limits_{{\bm s} \in \{ \newpv{\sets S_\Gamma}({\bm \iota})\} \cup \sets S^K \backslash \{ \newpv{\sets S_\Gamma}({\bm \iota}) \} } \quad \max\limits_{{\bm x} \in \sets R} \quad \min\limits_{ {\bm u} \in \newpv{\sets U_\Gamma}({\bm \iota} ,{\bm s})} \quad  {\bm u}^\top {\bm x} 
\quad = \quad
\min\limits_{{\bm s} \in \newpv{\sets S_\Gamma}({\bm \iota}) } \quad \max\limits_{{\bm x} \in \sets R} \quad \min\limits_{ {\bm u} \in \newpv{\sets U_\Gamma}({\bm \iota} ,{\bm s})} \quad {\bm u}^\top {\bm x}.
$$
Since the choice of ${\bm \iota} \in \sets C^K$ was arbitrary, this concludes the proof.
\Halmos
\endproof


\proof{Proof of Lemma~\ref{lem:offline_mmu_simplification}} 
\newpv{We begin by showing} that Problem~\eqref{eq:offline_mmu_3} is equivalent to
\begin{equation}
\maximize_{{\bm \iota} \in \sets C^K} \quad \min_{{\bm s}\in {\sets S}^K} \quad \max_{{\bm x} \in \sets R} \quad \min_{ {\bm u} \in \overline{\sets U}_\Gamma({\bm \iota} ,{\bm s}) } \quad  {\bm u}^\top {\bm x},
\label{eq:offline_mmu_2}
\end{equation}
\newpv{where
\begin{equation*}
\overline{\sets U}_\Gamma({\bm \iota},{\bm s}) := 
\left\{
\begin{array}{ccll}  
{\bm u} \in \sets U^0 & : &  \exists {\bm \epsilon} \in \sets E_\Gamma \text{ such that} \\
&& {\bm u}^\top ({\bm x}^{{\bm \iota}^\kappa_1} - {\bm x}^{{\bm \iota}^\kappa_2}) \; \geq \; -{\bm \epsilon}_\kappa &  \quad  \forall \kappa \in \sets K \; : \; {\bm s}_\kappa=1 \\
&& |{\bm u}^\top ( {\bm x}^{{\bm \iota}^\kappa_1} - {\bm x}^{{\bm \iota}^\kappa_2} )| \; \leq \; {\bm \epsilon}_\kappa  & \quad  \forall \kappa \in \sets K \; : \; {\bm s}_\kappa=0 \\
&& {\bm u}^\top ( {\bm x}^{{\bm \iota}^\kappa_1} - {\bm x}^{{\bm \iota}^\kappa_2} ) \; \leq \; {\bm \epsilon}_\kappa  &  \quad  \forall \kappa \in \sets K \; : \; {\bm s}_\kappa=-1
\end{array}
\right\}.
\end{equation*}
}
To this end, fix any ${\bm \iota} \in \sets C^K$. We show that
\begin{equation}
\min\limits_{{\bm s}\in \sets S^K} \quad \max\limits_{{\bm x} \in \sets R} \quad  \min\limits_{ {\bm u} \in \newpv{\sets U_\Gamma}({\bm \iota} ,{\bm s})} \quad  {\bm u}^\top {\bm x}
\label{eq:relax_1}
\end{equation}
and
\begin{equation}
\min\limits_{{\bm s}\in \sets S^K} \quad  \max\limits_{{\bm x} \in \sets R} \quad  \min\limits_{ {\bm u} \in \newpv{\overline{\sets U}_\Gamma}({\bm \iota} ,{\bm s}) } \quad   {\bm u}^\top {\bm x}
\label{eq:relax_2}
\end{equation}
attain the same optimal value. 

Since $\newpv{\sets U_\Gamma}({\bm \iota},{\bm s}) \subseteq \newpv{\overline{\sets U}_\Gamma}({\bm \iota},{\bm s})$, the optimal objective value of Problem~\eqref{eq:relax_2} constitutes a lower bound on the optimal objective value of Problem~\eqref{eq:relax_1}.

For the converse part, observe that, by virtue of the linearity of ${\bm u}^\top {\bm x}$ and the boundedness of $\newpv{\sets U_\Gamma}({\bm \iota},{\bm s})$, Problems~\eqref{eq:relax_1} and~\eqref{eq:relax_2} have the same objective value for all ${\bm s} \in \sets S^K$ such that $\newpv{\sets U_\Gamma}({\bm \iota},{\bm s}) \neq \emptyset$ since, in such cases, $\newpv{\overline{\sets U}_\Gamma}({\bm \iota},{\bm s}) = {\cl}(\newpv{\sets U_\Gamma}({\bm \iota},{\bm s}))$. Let ${\bm s}^\star$ be optimal in Problem~\eqref{eq:relax_2}. If $\newpv{\sets U_\Gamma}({\bm \iota},{\bm s}^\star) \neq \emptyset$, the proof is complete. Suppose instead that $\newpv{\sets U_\Gamma}({\bm \iota},{\bm s}^\star) = \emptyset$. Since $\newpv{\sets U_\Gamma}({\bm \iota},{\bm s}^\star) = \emptyset$ while $\newpv{\overline{\sets U}_\Gamma}({\bm \iota},{\bm s}^\star) \neq \emptyset$, this implies that there exist implied equalities in the set $\newpv{\overline{\sets U}_\Gamma}({\bm \iota},{\bm s}^\star)$\newpv{, which may only happen if $\Gamma=0$}. In particular, since $\sets U^0$ is full-dimensional, see Assumption~\ref{ass:linear_utility}, the implied equalities are all associated with response constraints. Let us collect the indices $\kappa$ of all those implied equalities in the set $\newpv{{\sets K}_{=}}$, i.e., let
$$
\newpv{{\sets K}_{=}} \; := \; \left\{ \kappa \in \sets K \; : \; {\bm s}_\kappa \in \{-1,1\} , \; {\bm u}^\top ( {\bm x}^{{\bm \iota}^\kappa_1} - {\bm x}^{{\bm \iota}^\kappa_2} ) = 0 \quad \forall {\bm u} \in \newpv{\overline{\sets U}_\Gamma}({\bm \iota},{\bm s}^\star) \right\}. 
$$
%
Define
$$
{\bm s}_\kappa' \; := \; 
\begin{cases}
    1 & \text{if } {\bm s}_\kappa = 1 \text{ and } \kappa \notin \newpv{{\sets K}_{=}}, \\
    0 & \text{if } {\bm s}_\kappa = 0 \text{ or } \kappa \in \newpv{{\sets K}_{=}}, \\
    -1 & \text{if } {\bm s}_\kappa = -1 \text{ and } \kappa \notin \newpv{{\sets K}_{=}}.
\end{cases}
$$
This definition of~${\bm s}'$ ensures that equality constraints that are implicit in the set $\newpv{\overline{\sets U}_\Gamma}({\bm \iota},{\bm s}^\star)$ are made explicit in the set $\newpv{\overline{\sets U}_\Gamma}({\bm \iota},{\bm s}')$. By construction, $\newpv{\overline{\sets U}_\Gamma}({\bm \iota},{\bm s}') = \newpv{\overline{\sets U}_\Gamma}({\bm \iota},{\bm s}^\star)$ and thus ${\bm s}'$ is optimal in Problem~\eqref{eq:relax_2}. Moreover, since $\newpv{\overline{\sets U}_\Gamma}({\bm \iota},{\bm s}')$ is non-empty, it follows by Proposition~2.3 in~\cite{NemhauserWolsey_IPbook} that $\newpv{{\sets U}_\Gamma}({\bm \iota},{\bm s}')$ is non-empty. Thus, $\newpv{\overline{\sets U}_\Gamma}({\bm \iota},{\bm s}') = {\cl}(\newpv{\sets U_\Gamma}({\bm \iota},{\bm s}'))$, implying that the objective value attained by ${\bm s}'$ in Problem~\eqref{eq:relax_1} equals the optimal objective value of Problem~\eqref{eq:relax_2}. This in turn yields the required result that the optimal objective value of Problem~\eqref{eq:relax_2} constitutes an upper bound on the optimal objective value of Problem~\eqref{eq:relax_1}.

Thus, the optimal values of Problems~\eqref{eq:relax_1} and~\eqref{eq:relax_2} are equal implying that Problems~\eqref{eq:offline_mmu_3} and~\eqref{eq:offline_mmu_2} are equivalent.
\newpv{It now suffices to show that Problems~\eqref{eq:offline_mmu_2} and~\eqref{eq:offline_mmu_4} are equivalent}. Fix any ${\bm \iota} \in \sets C^K$. Let ${\bm s}^\star$ be optimal in
\begin{equation}
\min\limits_{{\bm s}\in \sets S^K} \quad \max\limits_{{\bm x} \in \sets R} \quad \min\limits_{ {\bm u} \in \newpv{\overline{\sets U}_\Gamma}({\bm \iota} ,{\bm s})} \quad  {\bm u}^\top {\bm x}
\label{eq:indifference_notgood_1}
\end{equation}
and define ${\bm s}' \in \widetilde{\sets S}^K \subset \sets S^K$ through
$$
{\bm s}_\kappa' := 
\begin{cases}
{\bm s}^\star_\kappa & \text{if } {\bm s}^\star_\kappa \neq 0 \\
1 & \text{else.}
\end{cases}
$$
Then, $\newpv{\overline{\sets U}_\Gamma}({\bm \iota},{\bm s}') \supseteq \newpv{\overline{\sets U}_\Gamma}({\bm \iota},{\bm s}^\star)$, implying that the objective value attained by ${\bm s}'$ in Problem~\eqref{eq:indifference_notgood_1} is no greater than the objective value attained by ${\bm s}^\star$. Thus, ${\bm s}'$ is also optimal in~\eqref{eq:indifference_notgood_1}. Since ${\bm s}' \in \widetilde{\sets S}^K$, this implies that we can always restrict our search for an optimal solution to Problem~\eqref{eq:indifference_notgood_1} to the set $\widetilde{\sets S}^K$. As the choice of ${\bm \iota} \in \sets C^K$ was arbitrary \newpv{and noting that $\overline{\sets U}_\Gamma({\bm \iota},{\bm s}) = \widetilde{\sets U}_\Gamma({\bm \iota},{\bm s})$ for all ${\bm s} \in \widetilde{\sets S}^K$ and ${\bm \iota} \in \sets C^K$}, this concludes the proof.
\Halmos
\endproof


\proof{Proof of Theorem~\ref{thm:offline_mmu_complexity}} \emph{(a)} Fix any \newpv{$\Gamma \in \reals_+$}, $K \in \naturals$, ${\bm \iota} \in \sets C^K$, and ${\bm s} \in \widetilde{\sets S}^K$ such that $\newpv{\widetilde{\sets U}_\Gamma}({\bm \iota},{\bm s}) \neq \emptyset$. Then, $\newpv{\widetilde{\sets U}_\Gamma}({\bm \iota},{\bm s})$ is convex. If $\sets R$ is convex, it readily follows from the minimax theorem (see~\cite{Neumann_1928}) that
$$
\min_{ {\bm u} \in \newpv{\widetilde{\sets U}_\Gamma}({\bm \iota} ,{\bm s})} \quad \max_{{\bm x} \in \sets R}  \quad  {\bm u}^\top {\bm x} 
\quad = \quad 
\max_{{\bm x} \in \sets R}  \quad \min_{ {\bm u} \in \newpv{ \widetilde{\sets U}_\Gamma}({\bm \iota} ,{\bm s})} \quad {\bm u}^\top {\bm x}.
$$
Since the choices of \newpv{$\Gamma \in \reals_+$}, $K \in \naturals$, ${\bm \iota} \in \sets C^K$, and ${\bm s} \in \widetilde{\sets S}^K$ such that $\newpv{\widetilde{\sets U}_\Gamma}({\bm \iota},{\bm s}) \neq \emptyset$ were arbitrary, it holds that (under the premise of the theorem statement) Problem~\eqref{eq:offline_mmu_4} is equivalent to
\begin{equation}
    \maximize_{{\bm \iota} \in \sets C^K} \quad \min_{{\bm s}\in \widetilde{\sets S}^K} \quad \min_{ {\bm u} \in \newpv{\widetilde{\sets U}_\Gamma}({\bm \iota} ,{\bm s})} \quad \max_{{\bm x} \in \sets R}  \quad  {\bm u}^\top {\bm x}.
\end{equation}
Fix \newpv{$\Gamma \in \reals_+$} and $K \in \naturals$. Then,
\begin{equation}
\renewcommand{\arraystretch}{1.5}
    \begin{array}{cl}
    & \quad \max\limits_{{\bm \iota} \in \sets C^K} \quad \min\limits_{{\bm s}\in \widetilde{\sets S}^K} \quad  \min\limits_{ {\bm u} \in \newpv{\widetilde{\sets U}_\Gamma}({\bm \iota} ,{\bm s})} \quad \max\limits_{{\bm x} \in \sets R} \quad {\bm u}^\top {\bm x} \\
    = & \quad \max\limits_{{\bm \iota} \in \sets C^K} \quad \min\limits_{ {\bm u} \in \bigcup \limits_{\bm s \in \sets S^K} \newpv{\widetilde{\sets U}_\Gamma}({\bm \iota} ,{\bm s})} \quad \max\limits_{{\bm x} \in \sets R} \quad {\bm u}^\top {\bm x} \\
    = & \quad \max\limits_{{\bm \iota} \in \sets C^K} \quad \min\limits_{ {\bm u} \in \sets U^0} \quad \max\limits_{{\bm x} \in \sets R} \quad {\bm u}^\top {\bm x} \\
    = & \quad \min\limits_{ {\bm u} \in \sets U^0} \quad \max\limits_{{\bm x} \in \sets R} \quad {\bm u}^\top {\bm x} \\
    = & \quad \max\limits_{{\bm x} \in \sets R} \quad \min\limits_{ {\bm u} \in \sets U^0} \quad  {\bm u}^\top {\bm x},
\end{array}
\end{equation}
where the first equality follows since, for any $\bm \iota \in \sets C^K$, we have
$$
\left\{\bm u \in \newpv{\widetilde{\sets U}_\Gamma}({\bm \iota} ,{\bm s}): \bm s \in \widetilde{\sets S}^K \right\} \quad = \quad \bigcup \limits_{\bm s \in \widetilde{\sets S}^K} \newpv{\widetilde{\sets U}_\Gamma}({\bm \iota} ,{\bm s}) ,
$$
and the second equality follows by definitions of $\widetilde{\sets S}^K$ and $\newpv{\widetilde{\sets U}_\Gamma}({\bm \iota} ,{\bm s})$.
If the recommendation set~$\sets R$ is a polyhedron, then
\begin{equation*}
\max\limits_{{\bm x} \in \sets R} \quad \min\limits_{ {\bm u} \in \sets U^0} \quad  {\bm u}^\top {\bm x}
\end{equation*}
can be reformulated equivalently as a linear program of size polynomial in the size of the input parameters, see e.g., \cite{BenTal_Book}, and is thus polynomially solvable. Since the choice of $K \in \naturals$ above was arbitrary, this concludes the first proof of the first item.




\emph{(b)} We use a reduction from the following decision problem that is known to be $\mathcal{NP}$-complete, see~\cite{GJ79}.

\framebox{
\parbox[t][3cm]{15cm}{

\begin{center}
    \textsc{Partition.}
\end{center}

\textbf{Instance.} Given a set $\sets A$ of elements $\sets A := \{ 1 , \ldots, n \}$ with associated positive integer weights ${\bm w}_i \in \naturals_+$, $i\in \sets A$, such that $\displaystyle \sum_{i\in \sets A} {\bm w}_i = 2W$.

\textbf{Question.} Does there exist a partition of $\sets A$ into two subsets, $\sets X$ and $\sets A \backslash \sets X$, such that $\displaystyle \sum_{i \in \sets X} {\bm w}_i = \sum_{i \in \sets A \backslash \sets X} {\bm w}_i = W$? 
}}

\vspace{0.7cm}
We aim to reduce the partition problem to evaluating the objective function of an instance of Problem~\eqref{eq:offline_mmu_4}. To this end, fix an instance $(n,{\bm w},W)$ of the partition problem. Set \newpv{$\Gamma=0$,} $J=n + 1$, $K=n$, and $\sets U^0 = \{ {\bm u} \in [0,1]^J \; : \; {\bm u}_{J} = 0.5\}$. Also, for each $\kappa \in \sets K$, let ${\bm \iota}^\kappa_1 =  \1_\kappa$ and ${\bm \iota}^\kappa_2  = \1_{J}$. Then, given a choice ${\bm s} \in \widetilde{\sets S}^K$, we have
\begin{equation*}
\widetilde{\sets U}({\bm \iota},{\bm s}) = 
\left\{
\begin{array}{ccll}  
{\bm u} \in [0,1]^J & : & {\bm u}_J \; = \; 0.5 \\
&& {\bm u}_\kappa \; \geq \; 0.5 &  \quad  \forall \kappa \in \sets K \; : \; {\bm s}_\kappa=1 \\
&& {\bm u}_\kappa \; \leq \; 0.5  &  \quad  \forall \kappa \in \sets K \; : \; {\bm s}_\kappa=-1
\end{array}
\right\}.
\end{equation*}
Finally, we define the recommendation set $\sets R \subseteq \reals^J$ through $\sets R = \{ {\bm x}^1, {\bm x}^2 \}$, where ${\bm x}^1 = ( 2{\bm w}^\top, -2W )^\top$ and ${\bm x}^2 = ( ( - 2 {\bm w})^\top, 6W )^\top$. For any fixed ${\bm s}\in \widetilde{\sets S}^K$, define
$$
Z({\bm s}) \quad := \quad \max_{{\bm x} \in \sets R} \quad \min_{ {\bm u} \in \widetilde{\sets U}({\bm \iota} ,{\bm s}) } \quad  {\bm u}^\top {\bm x}.
$$
Then, we have
$$
\begin{array}{ccl}
Z({\bm s}) &  \quad  =  & \quad
\displaystyle \max \; \left\{ \left( \min_{{\bm u} \in \widetilde{\sets U}({\bm \iota},{\bm s})} \;\; \sum_{i \in \sets A} 2 {\bm u}_i {\bm w}_i  - 2W {\bm u}_J \right) , \; \left( \min_{{\bm u} \in \widetilde{\sets U}({\bm \iota},{\bm s})} \;\; \sum_{i \in \sets A} -2 {\bm u}_i {\bm w}_i + 6W {\bm u}_J  \right)  \right\} \\
&  \quad  =  & \quad
\displaystyle \max \; \left\{ \left( \min_{{\bm u} \in \widetilde{\sets U}({\bm \iota},{\bm s})} \;\; \sum_{i \in \sets A} 2 {\bm u}_i {\bm w}_i \right) - W , \; \left( \min_{{\bm u} \in \widetilde{\sets U}({\bm \iota},{\bm s})} \;\; \sum_{i \in \sets A} -2 {\bm u}_i {\bm w}_i + 2W   \right) + W  \right\} \\
&  \quad  =  & \quad
\displaystyle \max \; \left\{ \left( \min_{{\bm u} \in \widetilde{\sets U}({\bm \iota},{\bm s})} \;\; \sum_{i \in \sets A} 2 {\bm u}_i {\bm w}_i \right) - W , \; \left( \min_{{\bm u} \in \widetilde{\sets U}({\bm \iota},{\bm s})} \;\; \sum_{i \in \sets A} -2 {\bm u}_i {\bm w}_i + \sum_{i \in \sets A} {\bm w}_i   \right) + W  \right\} \\
&  \quad  =  & \quad
\displaystyle \max \; \left\{ \left( \min_{{\bm u} \in \widetilde{\sets U}({\bm \iota},{\bm s})} \;\; \sum_{i \in \sets A} 2 {\bm u}_i {\bm w}_i \right) - W , \; \left( \min_{{\bm u} \in \widetilde{\sets U}({\bm \iota},{\bm s})} \;\; \sum_{i \in \sets A} -2 ({\bm u}_i - 0.5) {\bm w}_i  \right) + W  \right\} \\
&  \quad  =  & \quad 
\displaystyle \max \; \left\{ \sum_{ \kappa \in \sets K : {\bm s}_\kappa =1  } {\bm w}_\kappa - W , \; \sum_{\kappa \in \sets K : {\bm s}_\kappa =1} - {\bm w}_\kappa + W   \right\} \\
&  \quad  \geq  & \quad 0.
\end{array}
$$
Now, we claim that we are given a ``yes'' instance of \textsc{Partition} if and only the objective value of~${\bm \iota}$ in the constructed instance of Problem~\eqref{eq:offline_mmu_4} is 0. To this end, note that the objective value of ${\bm \iota}$ is given by
\begin{equation}
\displaystyle \minimize_{{\bm s} \in \widetilde{\sets S}^K} \quad Z ({\bm s})
\label{eq:obj_eval_hard}
\end{equation}
and is lower bounded by 0. If there exists $\sets X \subset \sets A$ such that $\sum_{i \in \sets X} {\bm w}_i = W$, then the solution ${\bm s} \in \widetilde{\sets S}^K$ defined through ${\bm s}_\kappa =1$ if $\kappa \in \sets X$, $=-1$ else, $\kappa \in \sets K$, attains an objective value of zero in Problem~\eqref{eq:obj_eval_hard}. Conversely, if the optimal objective value of Problem~\eqref{eq:obj_eval_hard} is zero, then the set $\sets X := \{ \kappa \in \sets A : {\bm s}_\kappa = 1 \}$ is such that $\sum_{i \in \sets X} {\bm w}_i = W$ and the claim follows. This concludes the proof of the second item.

Both items are thus proved.
\Halmos
\endproof

\newpv{The following observation is needed in the proof of Lemma~\ref{lem:offline_mmu_4_finite_program}.
\begin{observation}
Problem~\eqref{eq:offline_mmu_4} is equivalent to the following max-min problem
\begin{equation}
\maximize_{{\bm \iota} \in \sets C^K} \quad \max_{\begin{smallmatrix} {\bm x}^{\bm s} \in \sets R : \\ {\bm s} \in \widetilde{\sets S}^K \end{smallmatrix}} \quad \min_{{\bm s}\in \widetilde{\sets S}^K} \quad \min_{ {\bm u} \in \newpv{\widetilde{\sets U}_\Gamma}({\bm \iota} ,{\bm s})} \quad  {\bm u}^\top {\bm x}^{\bm s},
\label{eq:offline_mmu_4_reorder}
\end{equation}
where ${\bm x}^{\bm s}$ denotes the item to recommend in response scenario~${\bm s}$, ${\bm s} \in \widetilde{\sets S}^K$.
\label{obs:offline_mmu_4_reorder}
\end{observation}}

\proof{Proof of Observation~\ref{obs:offline_mmu_4_reorder}}
Fix ${\bm \iota} \in \sets C^K$. Then,
\begin{equation}
    \begin{array}{cl}
    & \quad \displaystyle \min_{ {\bm s}\in \widetilde{\sets S}^K } \quad \max\limits_{{\bm x} \in \sets R} \quad \min\limits_{ {\bm u} \in \newpv{\widetilde{\sets U}_\Gamma}({\bm \iota} ,{\bm s})} \quad  {\bm u}^\top {\bm x} \\
    = & \quad \min  \;\; \left\{  \;\; \max\limits_{{\bm x}^{\bm s} \in \sets R}  \quad \min\limits_{ {\bm u} \in \newpv{\widetilde{\sets U}_\Gamma}({\bm \iota} ,{\bm s})}  \quad  {\bm u}^\top {\bm x}^{\bm s} \;\; \right\}_{{{\bm s}\in \widetilde{\sets S}}^K} \\
    = & \quad \displaystyle \max_{ \begin{smallmatrix} {\bm x}^{\bm s} \in \sets R : \\ {\bm s} \in \widetilde{\sets S}^K \end{smallmatrix} } \quad \min  \;\; \left\{ \;\;  \min\limits_{ {\bm u} \in \newpv{\widetilde{\sets U}_\Gamma}({\bm \iota} ,{\bm s})} \;\;  {\bm u}^\top {\bm x}^{\bm s} \;\;  \right\}_{{{\bm s}\in \widetilde{\sets S}}^K} \\
    = & \quad \displaystyle \max_{ \begin{smallmatrix} {\bm x}^{\bm s} \in \sets R : \\ {\bm s} \in \widetilde{\sets S}^K  \end{smallmatrix} } \quad \min_{ {\bm s}\in \widetilde{\sets S}^K} \quad  \min\limits_{ {\bm u} \in \newpv{\widetilde{\sets U}_\Gamma}({\bm \iota} ,{\bm s})} \;\;  {\bm u}^\top {\bm x}^{\bm s},  
    \end{array}
\end{equation}
where the second equality follows from the fact that each term in the minimum involves a different set of variables in the maximum, which can be optimized separately. Since the choice of ${\bm \iota} \in \sets C^K$ was arbitrary, the claim follows. 
\Halmos
\endproof


\proof{Proof of Lemma~\ref{lem:offline_mmu_4_finite_program}}
Problem~\eqref{eq:offline_mmu_4_reorder} can be written in epigraph form equivalently as
\begin{equation*}
\begin{array}{cl}
\maximize & \quad \tau \\
\subjectto  & \quad {\bm \iota} \in \sets C^K, \; {\bm x}^{\bm s} \in \sets R, \; {\bm s} \in \widetilde{\sets S}^K \\
& \quad \tau  \;\; \leq \;\; \displaystyle \min_{{\bm s}\in \widetilde{\sets S}^K} \;\;\; \min_{ {\bm u} \in \newpv{\widetilde{\sets U}_\Gamma}({\bm \iota} ,{\bm s})} \;\;\;  {\bm u}^\top {\bm x}^{\bm s}.
\end{array}
\end{equation*}
The problem above is in turn equivalent to
\begin{equation}
\begin{array}{cl}
\maximize & \quad \tau \\
\subjectto  & \quad \tau \in \reals, \; {\bm \iota} \in \sets C^K, \; {\bm x}^{\bm s} \in \sets R ,\; {\bm s} \in \widetilde{\sets S}^K \\
& \quad \tau  \;\; \leq \;\; \min\limits_{ {\bm u} \in \newpv{\widetilde{\sets U}_\Gamma}({\bm \iota} ,{\bm s})} \;\; {\bm u}^\top {\bm x}^{\bm s} \quad \forall {\bm s}\in \widetilde{\sets S}^K.
\end{array}
\label{eq:offline_mmu_4_reorder_2}
\end{equation}
To reformulate this robust problem which involves infinitely many constraints as a finite program, we employ techniques from robust optimization. Fix ${\bm s} \in \widetilde{\sets S}^K$, ${\bm \iota} \in \sets C^K$, and ${\bm x}^{\bm s} \in \sets R$, and consider the minimization subproblem associated with ${\bm s}$ in the epigraph constraint. This subproblem reads
\newpv{$$
\begin{array}{cllllll}
     \minimize & \quad  {\bm u}^\top {\bm x}^{\bm s} \\
     \subjectto  &  \quad {\bm u} \in \reals^J , \; {\bm \epsilon} \in \reals_+^K \\
     & \quad {\bm u}^\top [  {\bm s}_\kappa \; ({\bm x}^{{\bm \iota}^\kappa_1} - {\bm x}^{{\bm \iota}^\kappa_2}) ] + {\bm \epsilon}_\kappa \geq 0 &  \quad  \forall \kappa \in \sets K \\
     & \quad{\bm B} {\bm u} \geq {\bm b} \\
     & \quad \textbf{e}^\top {\bm \epsilon} \leq \Gamma.
\end{array}
$$}
Its dual is expressible as
%
%
\newpv{$$
\begin{array}{cll}
\maximize & \quad{\bm b}^\top {\bm \beta}^{\bm s} + \Gamma \mu^{\bm s} \\
\subjectto &\quad {\bm \alpha}^{\bm s} \in \reals^K_+, \; {\bm \beta}^{\bm s} \in \reals_+^M, \; \mu^{\bm s} \in \reals_- \\
& \quad \displaystyle \sum_{\kappa \in \sets K} \; {\bm s}_\kappa \; ({\bm x}^{{\bm \iota}^\kappa_1} - {\bm x}^{{\bm \iota}^\kappa_2}) \; {\bm \alpha}_\kappa^{\bm s} + {\bm B}^\top {\bm \beta}^{\bm s} \; = \; {\bm x}^{\bm s} \\
& \quad {\bm \alpha}^{\bm s} + \mu^{\bm s} \textbf{e} \leq {\bm 0}.
\end{array}
$$}

\newpv{Next, we claim that the primal-dual pair above satisfies strong duality. If ${\bm s} \in \newpv{\sets S_\Gamma}({\bm \iota})$, then $\newpv{\sets U_\Gamma}({\bm \iota},{\bm s})$ and thus also $\newpv{\widetilde{\sets U}_\Gamma}({\bm \iota},{\bm s})$ are non-empty. Since $\newpv{\widetilde{\sets U}_\Gamma}({\bm \iota},{\bm s})$ is non-empty and compact, by Assumption~\ref{ass:linear_utility}, it follows that the primal problem above is feasible and bounded, so that it is solvable. We show that the statement also holds if ${\bm s} \notin \newpv{ \sets S_\Gamma}({\bm \iota})$. If ${\bm s} \notin \newpv{\sets S_\Gamma}({\bm \iota})$, then the primal problem is infeasible. This implies that the dual is either infeasible or unbounded. We show that it cannot be infeasible. Fix ${\bm \alpha}={\bm 0}$ \newpv{and $\mu=0$} in the dual. Then, the dual reduces to
$$
\begin{array}{cll}
\maximize & \quad {\bm b}^\top {\bm \beta} \\
\subjectto & \quad {\bm \beta} \in \reals_+^M \\
& \quad {\bm B}^\top {\bm \beta}  = {\bm x}.
\end{array}
$$
By Farkas' lemma, exactly one of the following alternatives must hold: \emph{a)} There exists ${\bm \beta} \in \reals_+^M$ such that ${\bm B}^\top {\bm \beta} = {\bm x}$; or \emph{b)} There exists ${\bm u} \in \reals^J$ such that ${\bm B} {\bm u} \geq \bm 0$ and ${\bm u}^\top {\bm x} <0$. Since the set $\sets U^0$ is bounded, by Assumption~\ref{ass:linear_utility}, its recession cone coincides with the origin implying that assertion \emph{b)} cannot hold so that the dual must be feasible. We conclude that the dual must be unbounded so that the optimal objective values of the primal and dual problems coincide for all ${\bm s} \in \widetilde{\sets S}^K$.}

Replacing each minimization subproblem in Problem~\eqref{eq:offline_mmu_4_reorder_2} with its dual, we obtain the equivalent reformulation
\newpv{\begin{equation*}
\begin{array}{cl}
\maximize & \quad \tau \\
\subjectto  & \quad \tau \in \reals, \; {\bm \iota} \in \sets C^K \\
& \quad \!\!\left. \begin{array}{l}
{\bm x}^{\bm s} \in \sets R, \; {\bm \alpha}^{\bm s} \in \reals^K_+, \; {\bm \beta}^{\bm s} \in \reals_+^M ,\;  \mu^{\bm s} \in \reals_- \\
\tau  \;\; \leq \;\; {\bm b}^\top {\bm \beta}^{\bm s} + \Gamma \mu^{\bm s}  \\
\displaystyle \sum_{\kappa \in \sets K} \; {\bm s}_\kappa \; ({\bm x}^{{\bm \iota}^\kappa_1} - {\bm x}^{{\bm \iota}^\kappa_2}) \; {\bm \alpha}_\kappa^{\bm s} + {\bm B}^\top {\bm \beta}^{\bm s} \; = \; {\bm x}^{\bm s}  \\
 {\bm \alpha}^{\bm s} + \mu^{\bm s} \1 \leq {\bm 0}
\end{array} \quad \quad \right\} \quad \forall {\bm s}\in \widetilde{\sets S}^K,
\end{array}
\end{equation*}}
which concludes the proof. 
\Halmos
\endproof


\proof{Proof of Theorem~\ref{thm:offline_mmu_MILP}}
We begin by showing that Problem~\eqref{eq:offline_mmu_4_finite_program} is equivalent to
\newpv{
\begin{equation}
\begin{array}{cl}
\maximize & \quad \tau \\
\subjectto  & \quad \tau \in \reals, \; {\bm v}^\kappa, \; {\bm w}^\kappa \in \{0,1\}^I , \; \kappa \in \sets K  \\
& \quad {\bm \alpha}^{\bm s} \in \reals^K_+, \; {\bm \beta}^{\bm s} \in \reals_+^M ,\;  \mu^{\bm s} \in \reals_- ,\; {\bm x}^{\bm s} \in \sets R, \; {\bm s}\in \widetilde{\sets S}^K \\
& \quad \!\!\left. \begin{array}{l}
\tau  \;\; \leq \;\; {\bm b}^\top {\bm \beta}^{\bm s} + \Gamma \mu^{\bm s}  \\
\displaystyle \sum_{\kappa \in \sets K} \; {\bm s}_\kappa \; \sum_{i \in \sets I} {\bm x}^i \;  ( {\bm v}_i^{\kappa} - {\bm w}_i^{\kappa} ) \; {\bm \alpha}^{\bm s}_\kappa + {\bm B}^\top {\bm \beta}^{\bm s} \; = \; {\bm x}^{\bm s} \\
{\bm \alpha}^{\bm s} + \mu^{\bm s} \1 \leq {\bm 0}
\end{array} \quad \quad \right\} \quad \forall {\bm s}\in \widetilde{\sets S}^K \\
& \quad \!\!\left. \begin{array}{l}
 \1^\top {\bm v}^\kappa = 1, \; \1^\top {\bm w}^\kappa = 1 \quad \quad \quad \quad \\
\displaystyle  1-{\bm w}_i^\kappa \; \geq \; \sum_{i':i'\geq i} {\bm v}_{i'}^\kappa \quad \forall  i \in \sets I \quad \quad
\end{array} \quad  \right\} \quad \forall \kappa \in \sets K.
\end{array}
\label{eq:offline_mmu_MILP_2}
\end{equation}}

For the first direction, let~\newpv{$(\tau,{\bm \iota},\{ {\bm \alpha}^{\bm s}, {\bm \beta}^{\bm s}, \mu^{\bm s},{\bm x}^{\bm s} \}_{{\bm s}\in \widetilde{\sets S}^K})$} be feasible in Problem~\eqref{eq:offline_mmu_4_finite_program}. For each $\kappa\in \sets K$, define ${\bm v}^\kappa$ and ${\bm w}^\kappa \in \{0,1\}^I$ through
$$
{\bm v}_i^\kappa :=
\begin{cases}
1 & \text{if } {\bm \iota}^\kappa_1 = i \\
0 & \text{else,}
\end{cases} \qquad \text{ and } \qquad
{\bm w}_i^\kappa :=
\begin{cases}
1 & \text{if } {\bm \iota}^\kappa_2 = i \\
0 & \text{else,}
\end{cases}
$$
for each $i \in \sets I$. Fix $\kappa \in \sets K$. Since ${\bm \iota}^\kappa \in \sets C$, it readily follows that ${\bm \iota}^\kappa_1$ and ${\bm \iota}^\kappa_2$ are both in $\sets I$, implying that $\1^\top{\bm v}^\kappa=\1^\top {\bm w}^\kappa = 1$. Fix $i \in \sets I$. If ${\bm w}^\kappa_i = 0$, then it holds that $1-{\bm w}_i^\kappa \; \geq \; \sum_{i':i'\geq i} {\bm v}_{i'}^\kappa$. On the other hand, if ${\bm w}^\kappa_i=1$, then it holds that ${\bm \iota}^\kappa_2 = i$. But since ${\bm \iota}^\kappa \in \sets C$, it holds that ${\bm \iota}^\kappa_1 < {\bm \iota}^\kappa_2 = i$, i.e., there must exist $i' < i$ such that ${\bm \iota}^\kappa_1 = i'$ and ${\bm v}_{i'} = 1$. This in turn implies that $\sum_{i':i' \geq i} {\bm v}^\kappa_{i'} = 0$, and therefore $1-{\bm w}_i^\kappa \; \geq \; \sum_{i':i'\geq i} {\bm v}_{i'}^\kappa$ holds. Since the choice of $\kappa \in \sets K$ and $i \in \sets I$ was arbitrary, it holds that
$$
\left. \begin{array}{l}
 \1^\top {\bm v}^\kappa = 1, \; \1^\top {\bm w}^\kappa = 1 \quad \quad \quad \quad \\
\displaystyle  1-{\bm w}_i^\kappa \; \geq \; \sum_{i':i'\geq i} {\bm v}_{i'}^\kappa \quad \forall  i \in \sets I \quad \quad
\end{array} \quad  \right\} \quad \forall \kappa \in \sets K.
$$
Fix ${\bm s} \in \widetilde{\sets S}^K$ and $\kappa \in \sets K$. Then, 
$$
\sum_{i \in \sets I} {\bm x}^i \; ({\bm v}^\kappa_i - {\bm w}^\kappa_i) \; {\bm \alpha}^{\bm s}_\kappa \; = \; ( {\bm x}^{ {\bm \iota}^\kappa_1 } - {\bm x}^{ {\bm \iota}^\kappa_1 }) {\bm \alpha}^{\bm s}_\kappa
$$
holds trivially by definition of ${\bm v}^\kappa$ and ${\bm w}^\kappa$. Thus, \newpv{$(\tau,\{{\bm v}^\kappa,{\bm w}^\kappa\}_{\kappa \in \sets K},\{ {\bm \alpha}^{\bm s}, {\bm \beta}^{\bm s}, \mu^{\bm s},{\bm x}^{\bm s} \}_{{\bm s}\in \widetilde{\sets S}^K})$} is feasible in Problem~\eqref{eq:offline_mmu_MILP_2} with objective value equal to $\tau$.

For the other direction, let \newpv{$(\tau,\{{\bm v}^\kappa,{\bm w}^\kappa\}_{\kappa \in \sets K},\{ {\bm \alpha}^{\bm s}, {\bm \beta}^{\bm s},\mu^{\bm s}, {\bm x}^{\bm s} \}_{{\bm s}\in \widetilde{\sets S}^K})$} be feasible in Problem~\eqref{eq:offline_mmu_MILP_2}. Define
$$
{\bm \iota}_1^\kappa = \sum_{i \in \sets I} i \cdot \I{{\bm v}^\kappa_i = 1} \quad \text{ and } \quad {\bm \iota}_2^\kappa = \sum_{i \in \sets I} i \cdot \I{{\bm w}^\kappa_i = 1}.
$$
Fix $\kappa \in \sets K$. Then, ${\bm \iota}^\kappa_1$ and ${\bm \iota}^\kappa_2$ are both in the set $\sets I$. It follows from
$$
1-{\bm w}_i^\kappa \; \geq \; \sum_{i':i'\geq i} {\bm v}_{i'}^\kappa \quad \forall  i \in \sets I  
$$
that if ${\bm w}_i^\kappa = 1$, then ${\bm \iota}^\kappa_2 = i$ and $\sum_{i':i'\geq i} {\bm v}_{i'}^\kappa = 0$, implying that $\sum_{i':i' < i} {\bm v}_{i'}^\kappa = 1$, i.e., ${\bm \iota}^\kappa_1 < i = {\bm \iota}^\kappa_2$ and it holds that ${\bm \iota}^\kappa \in \sets C$. Since the choice of $\kappa$ was arbitrary, ${\bm \iota} \in \sets C^K$. Fix ${\bm s} \in \widetilde{\sets S}^K$ and $\kappa \in \sets K$. Then, 
$$
( {\bm x}^{ {\bm \iota}^\kappa_1 } - {\bm x}^{ {\bm \iota}^\kappa_1 }) \; {\bm \alpha}^{\bm s}_\kappa \; = \;    \sum_{i \in \sets I} {\bm x}^i \; ({\bm v}^\kappa_i - {\bm w}^\kappa_i) \; {\bm \alpha}^{\bm s}_\kappa 
$$
holds trivially by definition of ${\bm \iota}^\kappa$. Therefore, \newpv{$(\tau,{\bm \iota},\{ {\bm \alpha}^{\bm s}, {\bm \beta}^{\bm s},\mu^{\bm s}, {\bm x}^{\bm s} \}_{{\bm s}\in \widetilde{\sets S}^K})$} is feasible in Problem~\eqref{eq:offline_mmu_4_finite_program} with objective value equal to $\tau$. We have thus shown that Problems~\eqref{eq:offline_mmu_4_finite_program} and~\eqref{eq:offline_mmu_MILP_2} are equivalent.

Equivalence of the bilinear problem~\eqref{eq:offline_mmu_MILP_2} and the MILP~\eqref{eq:offline_mmu_MILP} follows directly by using standard linearization techniques which apply since all bilinear terms involve products of binary and real valued variables, see e.g., \cite{hillier2012introduction}. \Halmos
\endproof


\proof{Proof of Proposition~\ref{prop:offline_mmu_ccg_algo_correct_1}}
Since ${\bm \iota}$ is feasible in Problem~\eqref{eq:offline_mmu_ccg_rmp}, it follows that ${\bm \iota} \in \sets C^K$. By following a proof strategy similar to that in the Proof of \newpv{Lemma~\ref{lem:offline_mmu_4_finite_program}} (based on Farkas' lemma), it follows that~${\bm \iota}$ is feasible in Problem~\eqref{eq:offline_mmu_4}. Thus, it remains to show that for any given ${\bm \iota} \in \sets C^K$, Problem~\eqref{eq:offline_mmu_ccg_feas} is equivalent to
\begin{equation}
\minimize_{{\bm s}\in \widetilde{\sets S}^K} \quad \max\limits_{{\bm x} \in \sets R} \quad \min\limits_{ {\bm u} \in \newpv{\widetilde{\sets U}_\Gamma}({\bm \iota} ,{\bm s})} \quad  {\bm u}^\top {\bm x}
\label{eq:inner_offline_1}
\end{equation}
in the sense that the two problems have the same optimal objective value. To this end, fix ${\bm \iota}\in\sets C^K$. Using an epigraph reformulation, we can write  Problem~\eqref{eq:inner_offline_1} equivalently as
\begin{equation}
\begin{array}{cl}
    \minimize & \quad \theta \\
    \subjectto & \quad \theta\in \reals,\; {\bm s} \in \widetilde{\sets S}^K  \\
    & \quad \theta \; \geq \; \max\limits_{{\bm x} \in \sets R} \;\; \min\limits_{ {\bm u} \in \newpv{\widetilde{\sets U}_\Gamma}({\bm \iota} ,{\bm s})} \;\;  {\bm u}^\top {\bm x}.
\end{array}
\label{eq:epi_inner}
\end{equation}
Since $\sets R$ has fixed finite cardinality, we can equivalently express the above problem as 
\begin{equation*}
\begin{array}{cl}
    \minimize & \quad \theta \\
    \subjectto & \quad \theta\in \reals,\; {\bm s} \in \widetilde{\sets S}^K \\
    & \quad {\bm u}^{\bm x} \in \newpv{\widetilde{\sets U}_\Gamma}({\bm \iota} ,{\bm s})\quad \forall {\bm x} \in \sets R\; \\
    & \quad \theta \; \geq \; ({\bm u}^{\bm x})^\top {\bm x}  \quad  \forall {\bm x} \in \sets R.
\end{array}
\end{equation*}
From the definition of $\newpv{\widetilde{\sets U}_\Gamma}({\bm \iota} ,{\bm s})$, we can rewrite the preceding problem equivalently as
\newpv{\begin{equation*}
\begin{array}{cl}
    \minimize & \quad \theta \\
    \subjectto & \quad \theta\in \reals,\; {\bm u}^{\bm x} \in \sets U^0, \; {\bm \epsilon}^{\bm x} \in \sets E_\Gamma \;\;\; \forall {\bm x}\in \sets R, \;\;\;  {\bm s} \in \widetilde{\sets S}^K  \\
    & \quad \theta \; \geq \; ({\bm u}^{\bm x})^\top {\bm x} \quad \forall {\bm x} \in \sets R \\
    & \!\! \left.\begin{array}{l}
    \quad ({\bm u}^{\bm x})^\top ( {\bm x}^{{\bm \iota}_k} - {\bm x}^{{\bm \iota}_k'}) \; \geq \; -{\bm \epsilon}_\kappa^{\bm x} \quad  \forall \kappa \in \sets K \; : \; {\bm s}_\kappa=1\\
    \quad ({\bm u}^{\bm x})^\top ( {\bm x}^{{\bm \iota}_k} - {\bm x}^{{\bm \iota}_k'}) \; \leq \; {\bm \epsilon}_\kappa^{\bm x}  \quad  \forall \kappa \in \sets K \; : \; {\bm s}_\kappa=-1 \\ 
    \end{array} \quad \right\} \quad \forall \bm x \in \sets R,
\end{array}
\end{equation*}}
which, for $M$ sufficiently large, is equivalent to Problem~\eqref{eq:offline_mmu_ccg_feas}. Thus, Problems~\eqref{eq:inner_offline_1} and~\eqref{eq:offline_mmu_ccg_feas} are equivalent. Since the choice of ${\bm \iota}$ was arbitrary, the claim follows.
\Halmos
\endproof


\proof{Proof of Lemma~\ref{lem:offline_mmu_ccg_algo_correct_2}}
\begin{enumerate}[label=\emph{(\roman*)}]
    \item By virtue of Proposition~\ref{prop:offline_mmu_ccg_algo_correct_1}, $\theta$ yields a lower bound to the optimal value of Problem~\eqref{eq:offline_mmu_4}. At the same time, since ${\sets S}'\subseteq \widetilde{\sets S}^K$, it is evident that $\tau$ is an upper bound to the optimal objective value of Problem~\eqref{eq:offline_mmu_4}. Therefore, $\theta \leq \tau$.
    \item Suppose that $\theta = \tau$ and that there exists $\bm s \in \widetilde{\sets S}^K$ such that Problem~\eqref{eq:offline_mmu_ccg_sp} is infeasible. This implies that Problem~\eqref{eq:offline_mmu_ccg_rmp} is solvable and that there exists ${\bm s} \in \widetilde{\sets S}^K$ such that $\tau$ is strictly larger than the optimal objective value of
    \newpv{\begin{equation}
    \begin{array}{cl}
    \maximize & \quad {\bm b}^\top {\bm \beta} +\Gamma \mu \\
    \subjectto &\quad  {\bm x} \in \sets R, \; {\bm \alpha} \in \reals^K, \; {\bm \beta} \in \reals_+^M, \; \mu \in \reals_- \\
    & \quad\displaystyle \sum_{\kappa \in \sets K} {\bm s}_\kappa \;  \; ({\bm x}^{{\bm \iota}^\kappa_1} - {\bm x}^{{\bm \iota}^\kappa_2}) \; {\bm \alpha}_\kappa + {\bm B}^\top {\bm \beta}  = {\bm x} \\
    & \quad {\bm \alpha} + \mu \1 \leq {\bm 0}.
    \end{array}
    \label{eq:upper_bound_of_tau}
    \end{equation}}
    An inspection of the proof of Proposition~\ref{prop:offline_mmu_ccg_algo_correct_1} reveals that Problems~\eqref{eq:offline_mmu_ccg_feas} and~\eqref{eq:inner_offline_1} are equivalent. An inspection of the proof of Lemma~\ref{lem:offline_mmu_4_finite_program} shows that Problem~\eqref{eq:inner_offline_1} is equivalent to
    \newpv{\begin{equation}
    \begin{array}{ccl}
    \minimize\limits_{{\bm s}\in \widetilde{\sets S}^K} \phantom{xx} &\max & \quad {\bm b}^\top {\bm \beta}+\Gamma \mu \\
    &\st &\quad  {\bm x} \in \sets R, \; {\bm \alpha} \in \reals^K, \; {\bm \beta} \in \reals_+^M, \; \mu \in \reals_- \\
    && \quad\displaystyle \sum_{\kappa \in \sets K} {\bm s}_\kappa \; ({\bm x}^{{\bm \iota}^\kappa_1} - {\bm x}^{{\bm \iota}^\kappa_2}) \;  {\bm \alpha}_\kappa + {\bm B}^\top {\bm \beta}  = {\bm x} \\
    && \quad {\bm \alpha} + \mu \1 \leq {\bm 0}.
    \end{array}
    \label{eq:dual_of_inner_of_feas}
    \end{equation}}
    Thus, Problems~\eqref{eq:dual_of_inner_of_feas} and~\eqref{eq:offline_mmu_ccg_feas} are equivalent and Problem~\eqref{eq:dual_of_inner_of_feas} has an optimal objective value of~$\theta$. This implies that $\theta < \tau$, a contradiction.
    \item Suppose $\theta < \tau$, and let $\bm s$ be defined as in the premise of the lemma. Then, the proof of item \emph{(ii)} reveals that~$\bm s$ must be optimal in Problem~\eqref{eq:dual_of_inner_of_feas} with associated optimal value~$\theta$. This in turn implies that ${\bm s}$ is such that the optimal objective value of Problem~\eqref{eq:upper_bound_of_tau} is strictly less than $\tau$, implying that the ${\bm s}$th subproblem~\eqref{eq:offline_mmu_ccg_sp} is infeasible, which concludes the proof.
\end{enumerate}
We have thus proved all claims. \Halmos
\endproof


\proof{Proof of Theorem~\ref{thm:offline_mmu_ccg_algo_converges}}
First, note that finite termination is guaranteed since at each iteration, either ${\rm{UB}}-{\rm{LB}} \leq \delta$ (in which case the algorithm terminates) or a new set of constraints (indexed by the infeasible index~$\bm s$) is added to the master problem \ref{eq:offline_mmu_ccg_rmp}, see Lemma~\ref{lem:offline_mmu_ccg_algo_correct_2}. Since the set of all indices $\widetilde{\sets S}^K$ is finite, the algorithm will terminate in a finite number of steps. Second, by construction, at any iteration of the algorithm, $\tau$ (i.e., \rm{UB}) provides a upper bound on the optimal objective value of the problem. On the other hand, the returned (feasible) solution has as objective value $\theta$ (i.e., \rm{LB}). Since the algorithm only terminates if ${\rm{UB}}-{\rm{LB}} \leq \delta$, we are guaranteed that, at termination, the returned solution will have an objective value
that is within $\delta$ of the optimal objective value of the problem. This concludes the proof.
\Halmos
\endproof

\section{Proofs of Statements in Sections~\ref{sec:minimaxregret}}
\label{sec:EC_minimaxregret}


\newpv{The following lemma shows that, in a manner paralleling the risk averse case, Problem~\eqref{eq:offline_mmr} can be considerably simplified by eliminating the dependence of~$\sets S({\bm \iota})$ on~${\bm \iota}$, by dropping the ``indifferent'' scenarios, and by replacing the strict inequalities in the set~$\sets U({\bm \iota},{\bm s})$ by their loose counterparts. It is needed in the proof of Theorem~\ref{thm:offline_mmr_complexity}.
\begin{lemma}
Problem~\eqref{eq:offline_mmr} is equivalent to
\begin{equation}
\tag{$\widetildeofflineregret{K}$}
\minimize_{{\bm \iota} \in \sets C^K} \quad \max_{{\bm s}\in \widetilde{\sets S}^K} \quad \min_{{\bm x} \in \sets R} \quad \max_{ {\bm u} \in \newpv{\widetilde{\sets U}_\Gamma}({\bm \iota} ,{\bm s})} \quad \left\{ \;\; \max_{  {\bm x}' \in \sets R }  \;\; {\bm u}^\top {\bm x}' - {\bm u}^\top {\bm x} \;\; \right\},
\label{eq:offline_mmr_2}
\end{equation}
in the sense that the two problems have the same optimal objective value and the same sets of optimal solutions.
\label{lem:offline_mmr_simplification}
\end{lemma}}

\proof{Proof of Lemma~\ref{lem:offline_mmr_simplification}}
The proof parallels exactly the proof of Lemma~\ref{lem:offline_mmu_simplification} and can thus be omitted. \Halmos 
\endproof


\proof{Proof of Theorem~\ref{thm:offline_mmr_complexity}} We \newpv{work with the equivalent reformulation~\eqref{eq:offline_mmr_2} of~\eqref{eq:offline_mmr}, see Lemma~\ref{lem:offline_mmr_simplification} and} use a reduction from \textsc{Partition}, as in the proof of Theorem~\ref{thm:offline_mmu_complexity}. We aim to reduce the partition problem to evaluating the objective function of an instance of Problem~\eqref{eq:offline_mmr_2}. For convenience, we work with the negative of the objective function of Problem~\eqref{eq:offline_mmr_2} given by
$$
\minimize_{{\bm s}\in \widetilde{\sets S}^K} \quad \max_{{\bm x} \in \sets R} \quad \min_{ {\bm u} \in \newpv{\widetilde{\sets U}_\Gamma}({\bm \iota} ,{\bm s})} \quad \left\{ \;\; \min_{  {\bm x}' \in \sets R }  \;\; {\bm u}^\top {\bm x} - {\bm u}^\top {\bm x}' \;\; \right\}.
$$ 
To this end, fix an instance $(n,{\bm w},W)$ of the partition problem. Set \newpv{$\Gamma=0$}, $J=2n + 3$, $K=n$, and
$$
\sets U^0 \; = \; \left\{ 
\begin{array}{lcl}
{\bm u} \in \reals^{2n + 3} & : &
    {\bm u}_i \in [0,1] \quad \forall i \in \{1,\ldots, n\}   \\
  &&  {\bm u}_{n+i} \; = \; {\bm u}_i -  0.5 \quad \forall i \in \{1,\ldots, n\} \\
  &&  {\bm u}_{2n+1} \in [0,1] \\
  &&  {\bm u}_{2n+2} \in [-1,1] \\
  &&  {\bm u}_{2n+3} = 0.5
\end{array} \right\}.
$$
Also, for each $\kappa \in \sets K$, let ${\bm \iota}^\kappa_1 =  \1_\kappa$ and ${\bm \iota}^\kappa_2  = \1_{2n+3}$. Then, given a choice ${\bm s} \in \widetilde{\sets S}^K$, we have
$$
\newpv{\widetilde{\sets U}_\Gamma}({\bm \iota},{\bm s}) = 
\left\{
\begin{array}{lcl}  
{\bm u} \in \sets U^0 & : &
%
 {\bm u}_\kappa \; \geq \; 0.5  \quad  \forall \kappa \in \sets K \; : \; {\bm s}_\kappa=1 \\
&& {\bm u}_\kappa \; \leq \; 0.5   \quad  \forall \kappa \in \sets K \; : \; {\bm s}_\kappa=-1
\end{array}
\right\}.
$$
Finally, we define the recommendation set $\sets R \subseteq \reals^J$ through $\sets R = \{ {\bm x}^1, {\bm x}^2 \}$, where ${\bm x}^1 = ( {\bm w}^\top, {\bm 0}^\top, -W, 0, 0)^\top$ and ${\bm x}^2= ( {\bm 0}^\top, - {\bm w}^\top, W, 2W, 0 )^\top$. For any fixed ${\bm s}\in \widetilde{\sets S}^K$, define
$$
Z({\bm s}) \quad := \quad \max_{{\bm x} \in \sets R} \quad \min_{ {\bm u} \in \widetilde{\sets U}({\bm \iota} ,{\bm s})} \quad \min_{  {\bm x}' \in \sets R }  \;\; \left\{ \;  {\bm u}^\top {\bm x} - {\bm u}^\top {\bm x}' \; \right\}.
$$
Then, we have
$$
\begin{array}{ccl}
Z({\bm s}) &  \quad  =  & \quad
\displaystyle \max_{\bm x \in \sets R} \; \min \left\{ \min_{{\bm u} \in \widetilde{\sets U}({\bm \iota},{\bm s})} \;\; {\bm u}^\top ({\bm x}-{\bm x}^1), \;\min_{{\bm u} \in \widetilde{\sets U}({\bm \iota},{\bm s})} \;\; {\bm u}^\top ({\bm x}-{\bm x}^2) \right\} \\
&  \quad  =   & \quad
\displaystyle \max \;\left[ \min \left\{ \min_{{\bm u} \in \widetilde{\sets U}({\bm \iota},{\bm s})} \;\; {\bm u}^\top ({\bm x}^1-{\bm x}^1), \;\min_{{\bm u} \in \widetilde{\sets U}({\bm \iota},{\bm s})} \;\; {\bm u}^\top ({\bm x}^1-{\bm x}^2) \right\},  \right. \\
&& \qquad \qquad \displaystyle \left. \min \left\{ \min_{{\bm u} \in \widetilde{\sets U}({\bm \iota},{\bm s})} \;\; {\bm u}^\top ({\bm x}^2-{\bm x}^1), \;\min_{{\bm u} \in \widetilde{\sets U}({\bm \iota},{\bm s})} \;\; {\bm u}^\top ({\bm x}^2-{\bm x}^2) \right\} \right] \\
&  \quad  =   & \quad
 \displaystyle \max \;\left[ \min \left\{ 0, \;\min_{{\bm u} \in \widetilde{\sets U}({\bm \iota},{\bm s})} \;\; {\bm u}^\top ({\bm x}^1-{\bm x}^2) \right\}, \; \min \left\{ \min_{{\bm u} \in \widetilde{\sets U}({\bm \iota},{\bm s})} \;\; {\bm u}^\top ({\bm x}^2-{\bm x}^1), \; 0 \right\} \right].
\end{array}
$$
Next, observe that
$$
\begin{array}{cl}
  \quad   &  \displaystyle \min_{{\bm u} \in \widetilde{\sets U}({\bm \iota},{\bm s})} \;\; {\bm u}^\top ({\bm x}^1-{\bm x}^2)  \\
 = \quad  & \displaystyle  \min_{{\bm u} \in \widetilde{\sets U}({\bm \iota},{\bm s})} \;\; \sum_{i=1}^n {\bm u}_i {\bm w}_i + \sum_{i=1}^n {\bm u}_{n+i} {\bm w}_i - 2{\bm u}_{2n+1} W - 2{\bm u}_{2n+2} W + 0 \cdot {\bm u}_{2n+3}  \\
 = \quad  & \displaystyle  \min_{{\bm u} \in \widetilde{\sets U}({\bm \iota},{\bm s})} \;\; \sum_{i=1}^n ( 2{\bm u}_i - 0.5 ) {\bm w}_i -2W - 2 W  \\
 = \quad  & \displaystyle  \min_{{\bm u} \in \widetilde{\sets U}({\bm \iota},{\bm s})} \;\; \sum_{i=1}^n  2{\bm u}_i  {\bm w}_i - W -2W - 2 W  \\
 = \quad  & \displaystyle  \sum_{\kappa \in \sets K : {\bm s}_\kappa = 1}  {\bm w}_\kappa - W  - 4 W  \\
 \leq \quad  & \quad 0.
\end{array}
$$
Similarly, it holds that
$$
\begin{array}{cl}
  \quad   &  \displaystyle \min_{{\bm u} \in \widetilde{\sets U}({\bm \iota},{\bm s})} \;\; {\bm u}^\top ({\bm x}^2-{\bm x}^1)  \\
 = \quad  & \displaystyle  \min_{{\bm u} \in \widetilde{\sets U}({\bm \iota},{\bm s})} \;\; \sum_{i=1}^n - {\bm u}_i {\bm w}_i + \sum_{i=1}^n - {\bm u}_{n+i} {\bm w}_i + 2{\bm u}_{2n+1} W + 2{\bm u}_{2n+2} W  -  0 \cdot {\bm u}_{2n+3} \\
 = \quad  & \displaystyle  \min_{{\bm u} \in \widetilde{\sets U}({\bm \iota},{\bm s})} \;\; \sum_{i=1}^n -(2 {\bm u}_i - 0.5 ) {\bm w}_i  + 0  - 2 W  \\
 = \quad  & \displaystyle  \min_{{\bm u} \in \widetilde{\sets U}({\bm \iota},{\bm s})} \;\; \sum_{i=1}^n -2( {\bm u}_i - 0.5 ) {\bm w}_i - W + 0  - 2 W -2W +2W  \\
 = \quad  & \displaystyle    \sum_{\kappa \in \sets K : {\bm s}_\kappa = 1}  -{\bm w}_\kappa + W - 4 W  \\
 \leq \quad  & \quad 0.
\end{array}
$$
Therefore,
$$
\begin{array}{ccl}
Z({\bm s}) &  \quad  =  & \quad 
 \displaystyle \max \;\left\{ \displaystyle  \sum_{\kappa \in \sets K : {\bm s}_\kappa = 1}  {\bm w}_\kappa - W  , \;  \sum_{\kappa \in \sets K : {\bm s}_\kappa = 1}  -{\bm w}_\kappa + W   \right\} - 4 W \\
& \quad \geq  & \quad -4W.
\end{array}
$$
Now, we claim that we are given a ``yes'' instance of \textsc{Partition} if and only the objective value of ${\bm \iota}$ in the constructed instance of Problem~\eqref{eq:offline_mmu_4} is $-4W$. To this end, note that the objective value of ${\bm \iota}$ is given by
\begin{equation}
\displaystyle \minimize_{{\bm s} \in \widetilde{\sets S}^K} \quad Z ({\bm s})
\label{eq:obj_eval_hard_2}
\end{equation}
and is lower bounded by $-4W$. If there exists $\sets X \subset \sets A$ such that $\sum_{i \in \sets X} {\bm w}_i = W$, then the solution ${\bm s} \in \widetilde{\sets S}^K$ defined through ${\bm s}_\kappa =1$ if $\kappa \in \sets X$, $=-1$ else, $\kappa \in \sets K$, attains an objective value of $-4W$ in Problem~\eqref{eq:obj_eval_hard_2}. Conversely, if the optimal objective value of Problem~\eqref{eq:obj_eval_hard_2} is $-4W$, then the set $\sets X = \{ \kappa \in \sets A : {\bm s}_\kappa = 1 \}$ is such that $\sum_{i \in \sets X} {\bm w}_i = W$ and the claim follows. This concludes the proof of the second item. \Halmos
\endproof

\newpv{The following Lemma is needed in the proof of Theorem~\ref{thm:offline_mmr_MILP}.

\begin{lemma}
Problem~\eqref{eq:offline_mmr_2} is equivalent to the following finite program
\newpv{\begin{equation}
\begin{array}{cl}
\minimize & \quad \tau \\
\subjectto  & \quad \tau \in \reals, \; {\bm \iota} \in \sets C^K ,\; {\bm x}^{\bm s} \in \sets R, \; {\bm s} \in \widetilde{\sets S}^K \\
& \quad \!\!\left. \begin{array}{l}
{\bm \alpha}^{(\bm x', \bm s)} \in \reals_-^K, \; {\bm \beta}^{(\bm x', \bm s)} \in \reals_-^M  , \; \mu^{(\bm x', \bm s)} \in \reals_+ \\
\tau  \;\; \geq \;\; {\bm b}^\top {\bm \beta}^{(\bm x', \bm s)} + \Gamma \mu^{(\bm x', \bm s)} \\
\displaystyle \sum_{\kappa \in \sets K} {\bm s}_\kappa \; ({\bm x}^{{\bm \iota}^\kappa_1} - {\bm x}^{{\bm \iota}^\kappa_2}) \; {\bm \alpha}_\kappa^{(\bm x', \bm s)} + {\bm B}^\top {\bm \beta}^{(\bm x', \bm s)}  = {\bm x}'  - {\bm x}^{\bm s} \\
{\bm \alpha}^{(\bm x', \bm s)} + \mu^{(\bm x', \bm s)} \1 \geq {\bm 0}
\end{array} \quad \quad \right\} \quad \forall {\bm s}\in \widetilde{\sets S}^K, \;\;\forall \bm x' \in \sets R,
\end{array}
\label{eq:offline_mmr_finite_program}
\end{equation}}
where ${\bm \iota}$ denotes the queries to make and ${\bm x}^{\bm s}$ the items to recommend in response scenario ${\bm s}\in \widetilde{\sets S}^K$.
\label{lem:offline_mmr_finite_program}
\end{lemma}

\proof{Proof of Lemma~\ref{lem:offline_mmr_finite_program}}
In a way that parallels the proof of Observation~\ref{obs:offline_mmu_4_reorder}, it can be readily shown that Problem~\eqref{eq:offline_mmr_2} is equivalent to
\begin{equation*}
\minimize_{{\bm \iota} \in \sets C^K}  \quad \min_{\begin{smallmatrix} {\bm x}^{\bm s} \in \sets R : \\ {\bm s} \in \widetilde{\sets S}^K \end{smallmatrix}}  \quad \max_{{\bm s}\in \widetilde{\sets S}^K}\quad \max_{ {\bm u} \in \newpv{\widetilde{\sets U}_\Gamma}({\bm \iota} ,{\bm s})} \quad \left\{ \;\; \max_{  {\bm x}' \in \sets R }  \;\; {\bm u}^\top {\bm x}' - {\bm u}^\top {\bm x}^{\bm s} \;\; \right\}.
\end{equation*}
The above problem can be written in epigraph form equivalently as
\begin{equation*}
\begin{array}{cl}
\minimize & \quad \tau \\
\subjectto  & \quad \tau \in \reals, \; {\bm \iota} \in \sets C^K, \; {\bm x}^{\bm s} \in \sets R ,\; {\bm s} \in \widetilde{\sets S}^K \\
& \quad \tau  \;\; \geq \;\; \max\limits_{{\bm s}\in \widetilde{\sets S}^K} \;\; \max\limits_{ {\bm u} \in \newpv{\widetilde{\sets U}_\Gamma}({\bm \iota} ,{\bm s})} \;\;  \max\limits_{{\bm x'}\in \sets R} \;\; {\bm u}^\top ({\bm x'} - {\bm x}^{\bm s}),
\end{array}
\end{equation*}
which is in turn equivalent to
\begin{equation}
\begin{array}{cl}
\minimize & \quad \tau \\
\subjectto  & \quad \tau \in \reals, \; {\bm \iota} \in \sets C^K, \; {\bm x}^{\bm s} \in \sets R ,\; {\bm s} \in \widetilde{\sets S}^K \\
& \quad \tau  \;\; \geq \;\; \max\limits_{ {\bm u} \in \newpv{\widetilde{\sets U}_\Gamma}({\bm \iota} ,{\bm s})} \;\; {\bm u}^\top ({\bm x'} - {\bm x}^{\bm s}) \quad \forall {\bm s}\in \widetilde{\sets S}^K, \; \forall {\bm x'}\in \sets R .
\end{array}
\label{eq:regretaverse_recommendation_equiv_epigraph}
\end{equation}
Fix ${\bm s} \in \widetilde{\sets S}^K$, $\tau \in \reals$, ${\bm \iota} \in \sets C^K$, ${\bm x}^{\bm s} \in \sets R$, and $\bm x' \in \sets R$ and consider the associated maximization subproblem in the constraints of the above problem. This reads
\newpv{$$
\begin{array}{cllllll}
     \maximize & \quad  {\bm u}^\top ({\bm x'} - {\bm x}^{\bm s}) \\
     \subjectto  &  \quad {\bm u} \in \reals^J , \; {\bm \epsilon} \in \reals_+^K \\
     & \quad {\bm u}^\top [  {\bm s}_\kappa \; ({\bm x}^{{\bm \iota}^\kappa_1} - {\bm x}^{{\bm \iota}^\kappa_2}) ] + {\bm \epsilon}_\kappa \geq 0 &  \quad  \forall \kappa \in \sets K \\
     & \quad{\bm B} {\bm u} \geq {\bm b} \\
     & \quad \textbf{e}^\top {\bm \epsilon} \leq \Gamma.
\end{array}
$$}
Its dual is expressible as
\newpv{$$
\begin{array}{cll}
\minimize & \quad {\bm b}^\top {\bm \beta}^{(\bm x', \bm s)} + \Gamma \mu^{(\bm x', \bm s)} \\
\subjectto &\quad {\bm \alpha}^{(\bm x', \bm s)} \in \reals_-^K, \; {\bm \beta}^{(\bm x', \bm s)} \in \reals_-^M , \; \mu^{(\bm x', \bm s)} \in \reals_+\\
& \quad\displaystyle \sum_{\kappa \in \sets K}  {\bm s}_\kappa \; ({\bm x}^{{\bm \iota}^\kappa_1} - {\bm x}^{{\bm \iota}^\kappa_2}) \; {\bm \alpha}_\kappa^{(\bm x', \bm s)} + {\bm B}^\top {\bm \beta}^{\bm s}  = {\bm x}'-{\bm x}^{\bm s} \\
& \quad {\bm \alpha}^{(\bm x', \bm s)} + \mu^{(\bm x', \bm s)} \textbf{e} \geq {\bm 0}.
\end{array}
$$}
Following a proof strategy similar to that taken in the proof of \newpv{Lemma~\ref{lem:offline_mmu_4_finite_program}}, it can be shown that the optimal objective values of the primal-dual pair above are always equal (even when the primal is infeasible). Replacing each maximization subproblem in Problem~\eqref{eq:regretaverse_recommendation_equiv_epigraph} with its dual yields the equivalent formulation
\newpv{\begin{equation}
\begin{array}{cl}
\minimize & \quad \tau \\
\subjectto  & \quad \tau \in \reals, \; {\bm \iota} \in \sets C^K ,\; {\bm x}^{\bm s} \in \sets R, \; {\bm s} \in \widetilde{\sets S}^K \\
& \quad \!\!\left. \begin{array}{l}
{\bm \alpha}^{(\bm x', \bm s)} \in \reals_-^K, \; {\bm \beta}^{(\bm x', \bm s)} \in \reals_-^M  , \; \mu^{(\bm x', \bm s)} \in \reals_+\\
\tau  \;\; \geq \;\; {\bm b}^\top {\bm \beta}^{(\bm x', \bm s)} + \Gamma \mu^{(\bm x', \bm s)} \\
\displaystyle \sum_{\kappa \in \sets K} {\bm s}_\kappa \; ({\bm x}^{{\bm \iota}^\kappa_1} - {\bm x}^{{\bm \iota}^\kappa_2}) \; {\bm \alpha}_\kappa^{(\bm x', \bm s)} + {\bm B}^\top {\bm \beta}^{(\bm x', \bm s)}  = {\bm x}'  - {\bm x}^{\bm s} \\
{\bm \alpha}^{(\bm x', \bm s)} + \mu^{(\bm x', \bm s)} \textbf{e} \geq {\bm 0}
\end{array} \quad \quad \right\} \quad \forall {\bm s}\in \widetilde{\sets S}^K, \;\;\forall \bm x' \in \sets R.
\end{array}
\end{equation}}
This concludes the proof.  \Halmos
\endproof}


\proof{Proof of Theorem~\ref{thm:offline_mmr_MILP}}
The proof follows directly \newpv{from Lemma~\ref{lem:offline_mmr_finite_program}} by following an approach similar to that taken in the Proof of Theorem~\ref{thm:offline_mmu_MILP} and is thus omitted. \Halmos
\endproof


\newpv{
We introduce the following problem which is needed in the proof of Theorem~\ref{thm:offline_mmr_ccg_algo_converges}. Given variables $(\tau,{\bm \iota})$ feasible in the main problem, we define the $(\bm x', {\bm s})$th subproblem, $\bm x' \in \sets R$, ${\bm s} \in \widetilde{\sets S}^K$, through
\begin{equation}
\tag{$\mathcal{CCG}^{{\rm sub},(\bm x', {\bm s})}_{\text{mmr}}(\tau, {\bm \iota})$}
\begin{array}{cl}
\minimize & \quad 0 \\
\subjectto &\quad {\bm \alpha} \in \reals_-^K, \; {\bm \beta} \in \reals_-^M , \; \mu \in \reals_+ , \; {\bm x}\in \sets R \\
&\quad \tau  \;\; \geq \;\; {\bm b}^\top {\bm \beta} + \Gamma \mu \\
&\quad \displaystyle \sum_{\kappa \in \sets K} {\bm s}_\kappa \; ({\bm x}^{{\bm \iota}^\kappa_1} - {\bm x}^{{\bm \iota}^\kappa_2}) \; {\bm \alpha}_\kappa + {\bm B}^\top {\bm \beta} \; = \; {\bm x}' - {\bm x} \\
& \quad {\bm \alpha} + \mu \1 \geq {\bm 0}.
\end{array}
\label{eq:offline_mmr_ccg_sp}
\end{equation}

The following proposition and lemma are needed in the proof of Theorem~\ref{thm:offline_mmr_ccg_algo_converges}.
\begin{proposition}
Let ${\bm \iota}$ be feasible in the relaxed master problem~\eqref{eq:offline_mmr_ccg_rmp}. Then, ${\bm \iota}$ is feasible in Problem~\eqref{eq:offline_mmr_2} and the objective value of ${\bm \iota}$ in Problem~\eqref{eq:offline_mmr_2} is given by the optimal objective value of Problem~\eqref{eq:offline_mmr_ccg_feas}. Moreover, if the set $\sets R$ has fixed finite cardinality, then Problem~\eqref{eq:offline_mmr_ccg_feas} is a mixed-binary linear program of size polynomial in the size of the input.
\label{prop:offline_mmr_ccg_algo_correct_1}
\end{proposition}}

\proof{Proof of Proposition~\ref{prop:offline_mmr_ccg_algo_correct_1}}
Since ${\bm \iota}$ is feasible in Problem~\eqref{eq:offline_mmr_ccg_rmp}, it follows that ${\bm \iota} \in \sets C^K$. Thus,~${\bm \iota}$ is feasible in Problem~\eqref{eq:offline_mmr_2}. It remains to show that for any given ${\bm \iota} \in \sets C^K$, Problems~\eqref{eq:offline_mmr_ccg_feas} and
\begin{equation}
\max\limits_{{\bm s}\in \widetilde{\sets S}^K} \quad \min\limits_{{\bm x} \in \sets R} \quad \max\limits_{ {\bm u} \in \newpv{\widetilde{\sets U}_\Gamma}({\bm \iota} ,{\bm s})} \quad \left\{ \; \max\limits_{{\bm x'} \in \sets R} \quad  {\bm u}^\top {\bm x'} - {\bm u}^\top {\bm x} \; \right\}
\label{eq:inner_offline_1_regret}
\end{equation}
have the same optimal objective value. To this end, fix ${\bm \iota}\in\sets C^K$. Using an epigraph reformulation, we can write  Problem~\eqref{eq:inner_offline_1_regret} equivalently as
\begin{equation}
\begin{array}{cl}
    \maximize & \quad \theta \\
    \subjectto & \quad \theta\in \reals,\; {\bm s} \in \widetilde{\sets S}^K  \\
    & \quad \theta \;\; \leq \;\; \min\limits_{{\bm x} \in \sets R} \;\; \max\limits_{ {\bm u} \in \newpv{ \widetilde{\sets U}_\Gamma}({\bm \iota} ,{\bm s})} \;\; \max\limits_{{\bm x'} \in \sets R} \;\;  {\bm u}^\top {\bm x'} - {\bm u}^\top {\bm x},
\end{array}
\label{eq:epi_inner_regret}
\end{equation}
and thus as
\begin{equation*}
\begin{array}{cl}
    \maximize & \quad \theta \\
    \subjectto & \quad \theta\in \reals,\; {\bm s} \in \widetilde{\sets S}^K,\\
    & \quad {{\bm x}'}^{,\bm x} \in \sets R,\; {\bm u}^{\bm x} \in \newpv{\widetilde{\sets U}_\Gamma}({\bm \iota} ,{\bm s}) \quad \forall {\bm x} \in \sets R\; \\
    & \quad \theta \; \leq \; ({\bm u}^{\bm x})^\top  ( {{\bm x}'}^{,\bm x} - {\bm x} )  \;\;\; \quad \forall {\bm x} \in \sets R.
\end{array}
\end{equation*}
By definition of $\newpv{\widetilde{\sets U}_\Gamma}({\bm \iota} ,{\bm s})$, the above problem can be written as
\newpv{\begin{equation*}
\begin{array}{cl}
    \maximize & \quad \theta \\
    \subjectto & \quad \theta\in \reals,\; {\bm s} \in \widetilde{\sets S}^K \\
    & \quad {{\bm x}'}^{,\bm x} \in \sets R,\; {\bm u}^{\bm x} \in \sets U^0, \; {\bm \epsilon}^{\bm x} \in \sets E_\Gamma \quad \forall {\bm x} \in \sets R \\
    & \quad \theta \; \leq \; ({\bm u}^{\bm x})^\top ( {{\bm x}'}^{,\bm x} - {\bm x}) \quad \forall {\bm x} \in \sets R \\
    & \!\! \left.\begin{array}{l}
    \quad ({\bm u}^{\bm x})^\top ( {\bm x}^{{\bm \iota}_k} - {\bm x}^{{\bm \iota}_k'}) \; \geq \; -{\bm \epsilon}_\kappa^{\bm x} \quad  \forall \kappa \in \sets K \; : \; {\bm s}_\kappa=1\\
    \quad ({\bm u}^{\bm x})^\top ( {\bm x}^{{\bm \iota}_k} - {\bm x}^{{\bm \iota}_k'}) \; \leq \; {\bm \epsilon}_\kappa^{\bm x}  \quad  \forall \kappa \in \sets K \; : \; {\bm s}_\kappa=-1 \\ 
    \end{array} \quad \right\} \quad \forall \bm x \in \sets R.
\end{array}
\end{equation*}}
The claim then follows by rewriting the logical constraints above as linear constraints using a ``big-$M$'' constant.
\Halmos
\endproof


\newpv{\begin{lemma}
The relaxed master problem~\eqref{eq:offline_mmr_ccg_rmp} is always feasible. If~\eqref{eq:offline_mmr_ccg_rmp} is solvable, let $(\tau,{\bm \iota}, \{{{\bm \alpha}}^{(\bm x', \bm s)}, {\bm \beta}^{(\bm x', \bm s)} \}_{\bm x' \in {\sets R}, {\bm s} \in {\sets S}'}, \{ {\bm x}^{\bm s} \}_{{\bm s} \in {\sets S}'} )$ be an optimal solution. Else, if~\eqref{eq:offline_mmr_ccg_rmp} is unbounded, set $\tau = -\infty$ and let ${\bm \iota} \in \sets C^K$ be such that~\eqref{eq:offline_mmr_ccg_rmp} is unbounded when ${\bm \iota}$ is fixed to that value. Moreover, let $(\theta,\{ {\bm u}^{{\bm x}}, {{\bm x}'}^{,{\bm x}} \}_{{\bm x}\in \sets R}, {\bm s} )$ be optimal in Problem~\eqref{eq:offline_mmr_ccg_feas}. Then, the following hold:
\begin{enumerate}[label=(\roman*)]
    \item $\theta \geq \tau$;
    \item If $\theta = \tau$, then Problem~\eqref{eq:offline_mmr_ccg_sp} is feasible for all $\bm x' \in \sets R$ and ${\bm s} \in \widetilde{\sets S}^K$;
    \item If $\theta > \tau$, then there exists ${\bm x}'\in \sets R$ such that the pair~$(\bm{x'}, {\bm s})$ corresponds to an infeasible subproblem, i.e., Problem~\eqref{eq:offline_mmr_ccg_sp} is infeasible.
\end{enumerate}
\label{lem:offline_mmr_ccg_algo_correct_2}
\end{lemma}}

\proof{Proof of Lemma~\ref{lem:offline_mmr_ccg_algo_correct_2}}
\begin{enumerate}[label=\emph{(\roman*)}]
    \item By virtue of Proposition~\ref{prop:offline_mmr_ccg_algo_correct_1}, $\theta$  yields an upper bound to the optimal value of the Problem~\eqref{eq:offline_mmr_2}. At the same time, since $\widetilde{\sets R} \subseteq \sets R$ and ${\sets S}'\subseteq \widetilde{\sets S}^K$, it is evident that $\tau$ is a lower bound on the optimal objective value of Problem~\eqref{eq:offline_mmr_2}. Therefore, $\theta \geq \tau$.
    \item Suppose that $\theta = \tau$ and that there exists $\bm x' \in \sets R$ and $\bm s \in \widetilde{\sets S}^K$ such that Problem~\eqref{eq:offline_mmr_ccg_sp} is infeasible. This implies that Problem~\eqref{eq:offline_mmu_ccg_rmp} is solvable and that $\tau$ is strictly smaller than the optimal objective value of
    \newpv{\begin{equation}
    \begin{array}{cl}
    \minimize & \quad {\bm b}^\top {\bm \beta} + \Gamma \mu  \\
    \subjectto &\quad  {\bm \alpha} \in \reals_-^K, \; {\bm \beta} \in \reals_-^M ,\; \mu \in \reals_+,\; {\bm x} \in \sets R \\
    & \quad \displaystyle \sum_{\kappa \in \sets K} {\bm s}_\kappa \; ({\bm x}^{{\bm \iota}^\kappa_1} - {\bm x}^{{\bm \iota}^\kappa_2}) \; {\bm \alpha}_\kappa + {\bm B}^\top {\bm \beta} \; = \; {\bm x}' - {\bm x} \\
    & \quad {\bm \alpha} + \mu \textbf{e} \geq {\bm 0}.
    \end{array}
    \label{eq:upper_bound_of_tau_regret}
    \end{equation}}
    An inspection of the proof of Proposition~\ref{prop:offline_mmr_ccg_algo_correct_1} reveals that Problems~\eqref{eq:offline_mmr_ccg_feas} and~\eqref{eq:inner_offline_1_regret} are equivalent. An inspection of the proof of Lemma~\ref{lem:offline_mmr_finite_program} shows that Problem~\eqref{eq:inner_offline_1_regret} is equivalent to
    \newpv{\begin{equation}
    \begin{array}{ccl}
    \maximize \limits_{\bm x' \in \sets R,\; {\bm s}\in \widetilde{\sets S}^K} \;\; & \quad \min & \quad {\bm b}^\top {\bm \beta} + \Gamma \mu \\
    &\quad \st &\quad   {\bm \alpha} \in \reals_-^K, \; {\bm \beta} \in \reals_-^M ,\;\mu \in \reals_+, \; {\bm x} \in \sets R\\
    && \quad\displaystyle \sum_{\kappa \in \sets K} {\bm s}_\kappa \; ({\bm x}^{{\bm \iota}^\kappa_1} - {\bm x}^{{\bm \iota}^\kappa_2}) \;  {\bm \alpha}_\kappa + {\bm B}^\top {\bm \beta}  = {\bm x}' - {\bm x}\\
    && \quad {\bm \alpha} + \mu \textbf{e} \geq {\bm 0}.
    \end{array}
    \label{eq:dual_of_inner_of_feas_regret}
    \end{equation}}
    Thus, Problems~\eqref{eq:offline_mmr_ccg_feas} and~\eqref{eq:dual_of_inner_of_feas_regret} are equivalent and Problem~\eqref{eq:dual_of_inner_of_feas_regret} has an optimal value of $\theta$. This implies that $\tau < \theta$, a contradiction.
    \item Suppose $\theta > \tau$ and let $(\bm x', \bm s)$ be defined as in the premise of the lemma. Then, the proof of item \emph{(ii)} reveals that ${\bm s}$ is optimal in~\eqref{eq:inner_offline_1_regret} with associated optimal value $\theta$. This implies that ${\bm s}$ is such that the optimal objective value of Problem~\eqref{eq:upper_bound_of_tau_regret} is strictly greater than $\tau$, so that the $(\bm x', \bm s)$th subproblem~\eqref{eq:offline_mmr_ccg_sp} is infeasible, which concludes the proof. 
\end{enumerate}
All claims are thus proved. \flushright \Halmos
\endproof


\proof{Proof of Theorem~\ref{thm:offline_mmr_ccg_algo_converges}}
The proof mirrors the proof of Theorem~\ref{thm:offline_mmu_ccg_algo_converges} and is thus omitted in the interest of space. \Halmos
\endproof

\end{document}